%% file: arxiv.tex
\documentclass[
 onecolumn,
 superscriptaddress,
 amsmath,amssymb,
 aps
]{revtex4-2}

\usepackage{graphicx}
\usepackage{dcolumn}
\usepackage{bm}
\usepackage{hyperref}

\input{my_commands}

\newcommand{\SM}{Supplemental Material}

\begin{document}

\title{\textbf{An Energy-Based Mechanism for Compositional Behavior} 
}

\author{Francesca Rossi}
\email{Contact author: f.rossi@ssmeridionale.it}
\affiliation{Scuola Superiore Meridionale, Naples, Italy.}

\author{Veronica Centorrino}
\affiliation{Automatic Control Laboratory (IfA), ETH Z\"urich, Switzerland.}

\author{Francesco Bullo}
\affiliation{Center for Control, Dynamical Systems, and Computation, University of California, Santa Barbara, California, USA.}

\author{Giovanni Russo}
\affiliation{Department of Information and Electrical Engineering and Applied Mathematics, University of Salerno, Italy.}

\date{\today}

\begin{abstract}
Flexible intelligence relies on the ability to reuse previously acquired behaviors and combine them differently as circumstances change. In biological and artificial systems, this ability is often attributed to gating mechanisms that determine how much each available behavior should contribute at a given time. Yet these gating rules, the dynamics that compute them, and the neural circuits that may implement them are usually introduced separately, leaving unclear whether they reflect a common underlying principle. Here, we show that they can all be derived from a single variational principle for behavioral composition. The resulting mechanism naturally gives rise to softmax gating, evolves as an energy-based dynamical system with guaranteed convergence, and admits a recurrent neural network instantiation featuring context-dependent and local interactions. Across collective behavior, human decision-making, and layered control, the same mechanism reproduces characteristic behavioral patterns, provides interpretable accounts of how different behaviors are combined, and matches or outperforms established approaches. These results provide a unified account of how behavioral composition can emerge from a common principle, with implications for understanding flexible behavior in natural systems and for designing artificial agents that can adapt by recombining existing capabilities.
\end{abstract}
\maketitle

\section{\label{sec:introduction}Introduction}

A bird in a flock can modulate its social interactions with other birds to participate in distinct collective behaviors; a person can shift among decision strategies as their beliefs change. These examples illustrate a general capacity of biological agents and a hallmark of flexible intelligence: reusing and combining previously acquired behaviors~\cite{BT-KT-HT-TS:20}. Humans and other animals often reuse and recombine behavioral schemas,  continuously adjusting their relative contributions according to context~\cite{BML-TDU-JBT-SJG:17,JAM-EE:23,QM-CS:26}.
Yet, the principle that determines these context-dependent combinations and the mechanisms by which they can be realized in biological or artificial agents remain unclear~\cite{BT-KT-HT-TS:20,TJP-SPW:23,JAM-EE:23,PV-JFQ-SVF-JT:25}.

Experimental evidence in neuroscience suggests that the prefrontal cortex plays a key role in composing behavioral schemas, potentially via a gating mechanism that regulates information flow across neural circuits~\cite{RR-AG-MH:18,KJ-HL-MK-SE:07,BT-KT-HT-TS:20,LG-HC-YC-YN-SY:26}. These empirical observations have motivated computational models that combine recurrent neural networks with mixture-of-experts 
frameworks to assign context-dependent weights to reusable modules or policies~\cite{JW-ZKN-TD-SD-JL-DH-MB:18,BT-KT-HT-TS:20,WC-JJ-FW-JT-SK-JH:25,SM-RE:14}. Such models often implement gating through a softmax or another prespecified mapping. 
To both advance a more general understanding of policy composition mechanisms and enable their deployment in autonomous and adaptive systems, researchers have increasingly developed and evaluated gating-based computational models, often using decision-making and motor-control benchmarks~\cite{TJP-SPW:23}.
Although effective, these methods are typically tied to particular architectures and tasks, leaving open why a particular gating rule should arise or how its computation could be mechanistically implemented.

Despite remarkable cross-disciplinary efforts, the gating rules, their dynamics, and the neural circuits that implement them are typically introduced as separate modeling or design choices, leaving unresolved whether they can be derived from a common principle.
Deriving these elements from a common principle, while remaining consistent with existing conceptual frameworks, would explain why a particular gating mechanism emerges instead of merely demonstrating that a prescribed architecture works. Addressing this interdisciplinary challenge would both advance a unified understanding of compositional behavior in natural systems and provide principled foundations for designing artificial agents that can flexibly reuse and combine learned capabilities.

To address this challenge, this work shows that the gating rule, its dynamics, and its neural instantiation can all be derived from a single variational principle. 
Through this central idea, the gating rule is not postulated but instead emerges as the necessary consequence of the evolution of a dynamical system featuring a neural instantiation. More in detail, we cast the composition of behaviors as the minimization of a variational objective over mixtures of reusable policies, balancing divergence from a generative model against the entropy of the composition weights. 
This objective gives rise to a continuous-time proximal-gradient dynamics in which softmax gating emerges from the underlying optimization, and whose equilibrium encodes the optimal policy composition. The resulting dynamics is an energy model having as energy function the variational objective. The energy decreases monotonically along the trajectories of the dynamics, which remain on the probability simplex and converge exponentially to the optimum with explicit rates. Moreover, the same equations admit a recurrent neural network instantiation with context-dependent and local interactions, thereby linking the variational objective, its dynamical computation, and a neural network instantiation within a unified framework.

We evaluate the resulting mechanism in three domains that probe distinct aspects of policy composition: collective behavior, human decision-making, and layered control. In collective systems, the model dynamically coordinates social forces to reproduce polarization, milling, and goal-directed flocking. In human decision-making, it provides trial-by-trial estimates of how participants combine exploitation, uncertainty seeking, and risk aversion, while improving model evidence relative to established accounts. In control, it enables a drone to compose individually insufficient controllers to achieve trajectory tracking in a noisy environment and outperform standard baselines.
Collectively, the results advance our understanding of flexible behavior beyond current explanations and may guide the design of artificial agents capable of more adaptive behavioral composition.

\section{Results}\label{sec:results}

\begin{figure*}[b]
\centering
\includegraphics[width=1\linewidth]{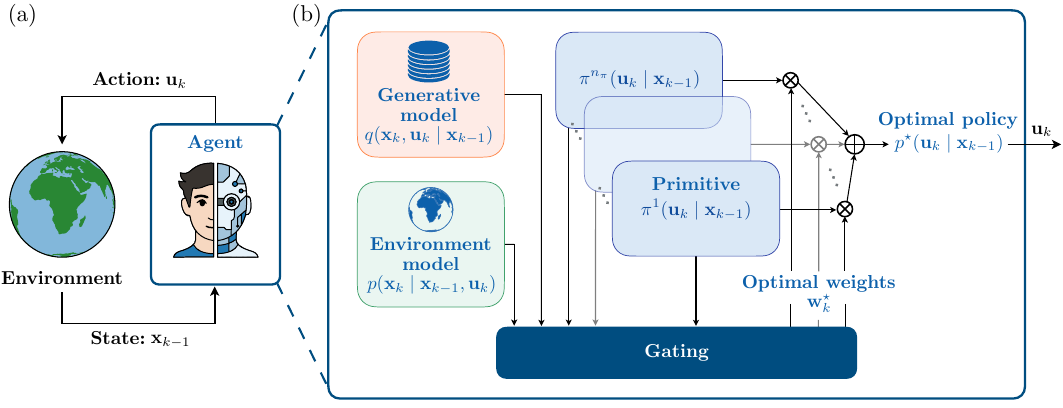}
\caption{Set-up. (a) At a given time step, an agent (e.g., a boid in a flock, a person in a multi-armed bandit task, a drone tracking a trajectory) receives the state $\bv{x}_{k-1}$ from the environment and determines the action $\bv{u}_k$. Bold symbols denote vector variables. For a boid or a drone, the state may be its position/velocity, and the action its acceleration or rotor thrust, respectively; for a person in a multi-armed bandit task, the state may represent the belief associated with the reward history, and the action may be the selected arm.
Both $\bv{x}_{k-1}$ and $\bv{u}_k$ are realizations of random variables, $\bv{X}_{k-1}$ and $\bv{U}_k$, respectively.
(b) At each time step, the agent computes the optimal policy $\optimalpolicy{k}{k-1}$ by composing a set of  primitives $\primitive{k}{k-1}{1}$, \dots, $\primitive{k}{k-1}{\np}$ via a gating mechanism. For boids in a flock, primitives may be social forces, for a person in a multi-armed bandit task these could be behavioral schemas, for a drone simple control policies. The optimal weights vector $\optimalweights{k}$ is obtained by solving an optimization problem that involves minimizing the statistical complexity (discrepancy, captured via a divergence) from a  generative model, $\refjointxu{k}{k-1}$. 
We use the term ``generative model'' broadly to denote a time-series model, specify a target or reference, encode a cost, or combine these elements.}
\label{fig:figure1}
\end{figure*}

In Fig.~\ref{fig:figure1}(a) an agent interacts with a stochastic environment. 
The environment transitions from state $\states{k-1}$ to $\states{k}$ in response to an action $\bv{u}_k$ sampled from policy $\policy{k}{k-1}$. The environment model available to the agent is $\plant{k}{k-1}$. 
The agent computes (Fig.~\ref{fig:figure1}(b)) the optimal policy $\optimalpolicy{k}{k-1}$ from $\np$ primitives (reusable  policies).

\subsection{Variational Objective for Policy Composition}

\begin{figure*}[b]
\centering
\includegraphics[width=1\linewidth]{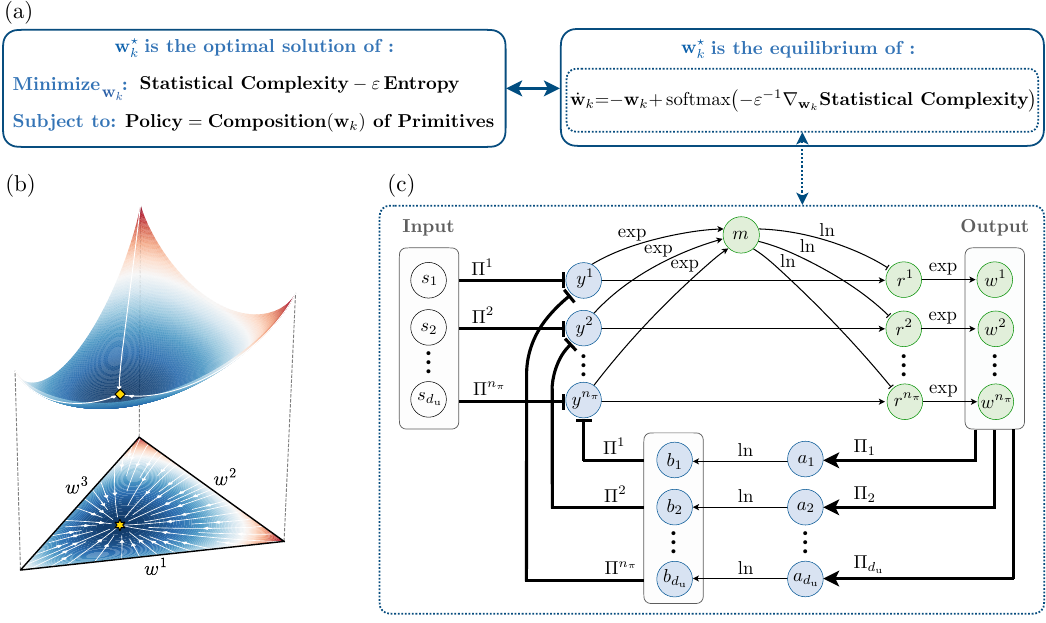}
\caption{Computational model overview.
(a) At each $k$, the agent computes the optimal weights $\bv{w}_k^\star$ for policy composition by solving an entropy-regularized optimization {problem}, finding a trade-off between statistical complexity (discrepancy from generative model) and entropy (left). The optimal weights are equal to the equilibrium of a continuous-time dynamical system defined by a softmax (proximal) gradient flow (right). 
(b) The softmax gradient flow is an energy model featuring highly ordered behavior, with guaranteed and explicit exponential convergence rates to the solution of the optimization problem in panel~(a). 
The optimization cost is an energy function (top) for the proximal gradient flow (for illustration purposes, a case with three primitives at a step $k$ is shown). The energy decreases along the trajectories of the flow toward the optimal solution (diamond).
Convergence is global (bottom): from any initial conditions in the simplex, the trajectories of the flow converge exponentially to the unique equilibrium, $\optimalweights{k}$, which is also the optimal solution of the normative optimization (star).
(c) Mechanistic neural instantiation of the flow at each time $k$. The circuit consists of two units operating at different time scales. The symbols $\Pi^\alpha$ and $\Pi_{i}$ denote the $\alpha$-th column and $i$-th row, respectively, of the synaptic matrix $\Pi$, dependent on state $\bv{x}$; $\Pi$ is the $\du \times \np$-dimensional matrix whose $\alpha$-th column is the $\du$-dimensional probability vector $\pi^\alpha$.
The fast unit (blue neurons) computes the gradient of the objective through a feedforward motif featuring contextual synaptic projections (associated to the $\Pi$ synapses, highlighted with thicker lines) and a hidden layer with logarithmic activations. The slow unit (green neurons) implements the softmax. Jointly, the units realize divisive normalization~\cite{DH-KZ:20} applied to exponentiated inputs determined by the gradient. 
The input is a cost combining a mismatch from the generative model and a negative log-likelihood. Arrowheads denote excitatory connections; T-bar heads denote inhibitory connections.
}
\label{fig:figure2}
\end{figure*}

We derive a top-down normative variational framework for gating that applies to a wide range of decision-making problems. 
The underlying optimization is provably solvable via a continuous-time dynamics yielding a neural circuit. In a different context, neural circuits derived from the similarity-matching normative framework have been described in the literature~\cite{CP-AMS-DBC:17,NMC-CP-DBC:23}.

The weights $\weights{k}$ for composing the primitives $\primitive{k}{k-1}{1},\dots,\primitive{k}{k-1}{\np}$ 
are computed via an optimization problem (Fig.~\ref{fig:figure2}(a), left). 
The cost comprises two terms: (i) a statistical complexity penalty, specifically the Kullback-Leibler (KL) divergence between the agent–environment dynamics $\jointxu{k}{k-1} = \plant{k}{k-1}\policy{k}{k-1}$ and the generative model $\refjointxu{k}{k-1}$, and 
(ii) an entropic regularizer scaled by a temperature-like parameter $\eps$. Formally,  
\begin{align}
\min_{\weights{k}\in \simplex{\np}} & \!\!\!\! \overbrace{{\DKLb{\jointxu{k}{k-1}\!}{\!\refjointxu{k}{k-1}}}}^{\text{Statistical Complexity}} {-} \eps \overbrace{\entropy{(\weights{k})}}^{\text{Entropy}} \nonumber\\
\text{s.t. }& \policy{k}{k-1} = \underbrace{\sum_{\alpha=1}^{\np}\weight{k}{\alpha} \primitive{k}{k-1}{\alpha}}_{\text{Composition of Primitives}},
\label{eq:problem1}
\end{align}
where $\simplex{\np}$ is the probability simplex. Embedding the constraint in the cost makes the dependence on the decision variables explicit and reveals that the problem is strictly convex (see Sec.~\ref{sec:methods}). To stress the dependence of the statistical complexity on $\weights{k}$, we use the notation $\F(\weights{k})$ so that the cost in Eq.~\eqref{eq:problem1} is written compactly as
$\F(\weights{k}) - \eps\entropy{(\weights{k})}$. 

Eq.~\eqref{eq:problem1} is an entropy-regularized optimization problem that arises in the context of, e.g.,  inference, maximum entropy and free energy principle. To elucidate this point, note that the generative model can be written as $\frac{1}{Z}\refjointxutwisted{k}{k-1}\e^{-c(\bv{x}_k,\bv{u}_k)}$.
Here, $Z$ is a normalization constant, $\refjointxutwisted{k}{k-1}$ is a probability, and $c(\bv{x}_k,\bv{u}_k)$ encodes a state/action cost. Applying the logarithm product rule, the complexity term in Eq.~\eqref{eq:problem1} becomes $\DKLb{\jointxu{k}{k-1}}{\refjointxutwisted{k}{k-1}} + \E_{\jointxu{k}{k-1}}\!\left[c(\bv{X}_k,\bv{U}_k)\right] + \ln Z$ and the constant $\ln Z$ can be dropped from the optimization. The resulting functional is minimized to find optimal policies under, e.g., the free energy principle~\cite{AS-HJ-KF-GR:25_new} and arises across active inference, KL control, control-as-inference, maximum entropy learning frameworks as well as entropy-regularized optimal transport~\cite{ET:09,MB-MT:12,HK-VC-MO:12,BDZ-AM-JAB-AKD:08,PM-TV-BD:21,TP-GP-KF:22,YK-HN:23}.
The optimization is also relevant to combinatorial optimization~\cite{ZSS-FP-YW-YDM-WBX-MHY-PZ:25} and variational inference, where the weights serve as variational parameters~\cite[Chapter $10$]{KM:23}.

\subsection{Softmax Gating from the Variational Objective}
The optimization in Eq.~\eqref{eq:problem1} can be solved via a continuous-time dynamics, the softmax proximal gradient flow (Fig.~\ref{fig:figure2}(a), right). We derive the dynamics here, and refer to the next section for its convergence properties. For an optimization problem with cost $f + g$, where $f$ is convex and $g$ is a regularizer, the proximal gradient flow generalizes gradient descent by handling possibly poorly behaved regularizers through a proximal correction step. Just like the dynamics of gradient descent, the proximal gradient flow is determined by the optimization cost.
We investigate the structure of the proximal gradient flow, revealing three key results. First, the dynamics converges to an equilibrium that is also the solution of Eq.~\eqref{eq:problem1}.  {Second, the equilibrium defines a softmax gating rule, arising from the entropic regularizer}.
Third, the statistical complexity in Eq.~\eqref{eq:problem1} determines how much each primitive contributes to the optimal policy. Details and definitions in Sec.~\ref{sec:methods}.

To establish the results, we reformulate the constrained optimization in Eq.~\eqref{eq:problem1} into the unconstrained problem
\begin{equation}
\label{eq:entropy_regularized_problem_simplex_weights}
\min_{\weights{k} \in \R^{\np}} \F(\weights{k}) + \eps \entropybar(\weights{k}).
\end{equation}
Here, $\F(\weights{k})$ is defined as in the previous section, and the {\em entropic barrier} $\entropybar(\weights{k})$ is equal to the negative entropy on the simplex and infinity outside of this set. Although Eq.~\eqref{eq:entropy_regularized_problem_simplex_weights} is unconstrained, the entropic barrier forces the optimal solution to be in the simplex, so that the optimal solution of Eq.~\eqref{eq:entropy_regularized_problem_simplex_weights} coincides with that of Eq.~\eqref{eq:problem1}. From this reformulation, we derive the associated proximal gradient {flow}~\cite{SHM-MRJ:21, VC-AG-AD-GR-FB:23a, AD-VC-AG-GR-FB:23f}, that we term softmax gradient flow
\begin{equation}
\label{eq::FreeNet}
\tau \dot{\bv{w}}_k = - \weights{k} + \softmax\bigl(- \eps^{-1} \nabla \F( \weights{k})\bigr).
\end{equation}
The parameter $\tau >0$ is a time-scale and $\nabla \F( \weights{k})$ is the gradient of $\F(\weights{k})$ with respect to $\weights{k}$. The $\softmax$ emerges from the proximal operator of the entropic regularization (\SM, Sec.~3). Each state variable (primitive  weight) evolves according to the corresponding component of the softmax applied to the scaled negative gradient of the statistical complexity. Therefore, the dynamics downweights primitives inducing a larger statistical complexity. 

The optimal solution of Eq.~\eqref{eq:entropy_regularized_problem_simplex_weights}, $\optimalweights{k} = \softmax\bigl(- \eps^{-1} \nabla \F( \optimalweights{k})\bigr)$, is an equilibrium of the softmax gradient flow in Eq.~\eqref{eq::FreeNet} and defines a softmax gating rule.
This broadly used mechanism is typical of architectures where all primitives contribute to the policy. As $\eps$ decreases, the softmax approaches the argmax so that the primitive with lowest gradient is selected with increasing probability~\cite{FR-EG-GR:25}.
When a weight bias $\pweights{k}$ is included in our normative framework, the entropy regularizer is replaced by a KL divergence from the bias, and Eq.~\eqref{eq::FreeNet} returns a prior-weighted softmax gating rule (\SM, Sec.~4).

\subsection{The Gating Dynamics Is an Energy Model}
To solve the problem in Eq.~\eqref{eq:problem1}, the softmax gradient flow Eq.~\eqref{eq::FreeNet} must admit the optimal solution $\optimalweights{k}$ as its equilibrium and guarantee convergence to it. Here we show that the softmax gradient flow is an energy model and that its trajectories $\weights{k}(t)$ converge exponentially to $\optimalweights{k}$ from any initial condition $\weights{k}(0)$ in the simplex. Derivations in Sec.~\ref{sec:methods}.

The energy function of the softmax gradient flow is the cost of Eq.~\eqref{eq:entropy_regularized_problem_simplex_weights}, which we denote by $V(\weights{k})$.
For all $t >0$, along the trajectories of the softmax gradient flow, it holds that 
\begin{equation}
\label{eq:dissipation_inequality_main}
\begin{aligned}
\frac{d}{dt}V(\weights{k}(t)) = -\frac{\eps}{\tau} \Bigl( \DKL{\weights{k}(t)\!}{\!\softmax(-\eps^{-1}\nabla \F(\weights{k}))} + \DKL{\softmax(-\eps^{-1}\nabla \F(\weights{k}))\!}{\!\weights{k}(t)} \Bigr).
\end{aligned}
\end{equation}
The symmetric Kullback--Leibler divergence on the right-hand side of Eq.~\eqref{eq:dissipation_inequality_main} is nonnegative and vanishes only at the equilibrium $\optimalweights{k}$ of Eq.~\eqref{eq::FreeNet}. Therefore, $V$ is strictly decreasing along every solution of the softmax gradient flow that is not the equilibrium. 
Moreover, for any initial condition belonging to the simplex, $\weights{k}(t)$ remains in the simplex for all $t\ge0$. Hence, the trajectories of the softmax gradient flow are always feasible for the problem in Eq.~\eqref{eq:problem1}.

The dissipation identity in Eq.~\eqref{eq:dissipation_inequality_main} also reveals two key convergence results with guaranteed convergence rates. First, the energy exponentially converges to the minimum $V(\optimalweights{k})$ with rate ${\tau}^{-1}$. Second, the convergence of $\weights{k}(t)$ to the optimal solution $\optimalweights{k}$ is also exponential with rate $(2\tau)^{-1}$.

These results show that the softmax gradient flow is a continuous-time energy model provably solving the optimization postulated under our normative framework. Remarkably, the softmax function also arises in the context of best response maps~\cite{BG-LP:17}, thus making our results potentially relevant to game-theoretical settings~\cite{SP-NEL:25,NM-MW-MB:23} (\SM, Sec.~3).


\subsection{The Flow Admits a Neural Instantiation}
\label{sec:bio_implementation}
We instantiate the softmax gradient flow in a neural circuit (Fig.~\ref{fig:figure2}(c)). This circuit is a soft-competitive continuous-time recurrent network with contextual computations~\cite{AM-MG-AS-DE-VP-MW:25}. Each element of the architecture is derived from the normative framework, specifically from Eq.~\eqref{eq::FreeNet}. Our approach therefore provides a theoretical understanding of the computational principles of the circuit.

The circuit consists of two units: 
a fast unit, computing the argument of the $\softmax$ in Eq.~\eqref{eq::FreeNet} and a slow unit implementing the $\softmax$ itself. 
To derive the fast unit, given $\du$ discrete actions, computing the $\softmax$ argument reduces to evaluating the $\np$-dimensional vector  $\bar{\bv{y}}_k {=} - \eps^{-1}\Pi(\bv{x}_{k-1})^\top\bigl(\ln ({\Pi(\bv{x}_{k-1}){\weights{k}}}) + \bv{s}(\bv{x}_{k-1},\bv{u}_{k})\bigr)$.
Given the state $\bv{x}_{k-1}$, $\Pi(\states{k-1})$ is the $\du \times \np$-dimensional matrix having on its $\alpha$-th column, $\Pi^\alpha(\bv{x}_{k-1})$, the primitive $\primitive{k}{k-1}{\alpha}$, and $\bv{s}(\bv{x}_{k-1} , \bv{u}_{k})$ is a $\du$-dimensional input.
This input is the sum of a statistical complexity and a negative log-likelihood, i.e., $\DKLb{\plant{k}{k-1}}{\refplant{k}{k-1}} - \ln\refpolicy{k}{k-1}$.
Thus, given the current state, the neural circuit processes an input associated to the action space to compute the gating weights. 
A growing body of experimental and computational evidence suggests that these quantities, which we obtain through theoretical derivations, may be processed in brain circuits~\cite{AL-QZ-MM-SM-JD-NU:25,ACC-MG-LM-PK-VN-FB-RC-JR-SC-SP-JB:23}.

The fast unit computes $\bar{\bv{y}}_k$. Denoting $\bv{s}(\bv{x}_{k-1}, \bv{u}_k)$ by $\bv{s}$, this computation is implemented via the dynamics
\begin{subequations}
\label{eq:gradient_bio}
\begin{align}
\subscr{\tau}{g} \dot{\bv{a}} &=-\bv{a}+\Pi(\states{k-1}) \bv{w},\label{eq:gradient_bio_a}\\
\subscr{\tau}{g} \dot{\bv{b}} &=-\bv{b}+\ln{(\bv{a})},\label{eq:gradient_bio_b}\\
\subscr{\tilde\tau}{g} \dot{\bv{y}} &= - \eps\bv{y} -\Pi(\states{k-1})^\top (\bv{b} + \bv{s}),\label{eq:gradient_bio_c}
\end{align}
\end{subequations}
where $\subscr{\tau}{g} \ll \subscr{\tilde\tau}{g} $ are time-scale constants. 

Following Eq.~\eqref{eq:gradient_bio_a}, the neuron $a_i$ receives the projected scalar quantity $\Pi_i(\states{k-1})\bv{w}$, where $\Pi_i(\states{k-1})$ is the $i$-th row of $\Pi(\states{k-1})$. The dynamics in Eqs.~\eqref{eq:gradient_bio_a}--\eqref{eq:gradient_bio_c} are exponentially stable and converge to $\Pi(\states{k-1})\bv{w}$, $\ln(\bv{a})$, and $\bar{\bv{y}}$, respectively. The time-scale separation $\subscr{\tau}{g}\ll\subscr{\tilde\tau}{g}$ ensures that the $\bv{a}$ and $\bv{b}$ layers converge before the $\bv{y}$ layer, so that the fast unit computes the desired softmax argument.

The dependence of the synaptic matrix $\Pi$ on the state  shows the presence of context-dependent synaptic couplings in the circuit.
In neural systems, context-dependent modulation might arise from multiplicative (Sigma--Pi-type) interactions between contextual variables and inputs~\cite{BM-CK:89,CK-TP:92,DL-CW-SG-DC:20}. 
At the cellular level, dendritic mechanisms have long been suggested to support context-dependent computations~\cite{ML-MH:05}.
Related forms of state-dependent modulation might also emerge from neuron--glia interactions, whereby astrocytes regulate synaptic efficacy based on local network activity~\cite{GP-MN-AA:09}.

The output of the fast unit feeds the slow unit, implementing the softmax in Eq.~\eqref{eq::FreeNet} and returning the weights. The circuit leverages a neural implementation of the softmax exhibiting elements of biological plausibility~\cite{MS-JO:22}. The slow unit dynamics is
\begin{subequations}
\label{eq:softmax_bio_main}
\begin{align}
\subscr{\tau}{s} \dot m &= - m + \sum_{\alpha=1}^{\np} \e^{y^\alpha},\label{eq:softmax_bio_main_a}\\
\subscr{\tau}{s} \dot{\bv{r}} &= - \bv{r} + \bv{y} - \1_{\np}\ln(m),\label{eq:softmax_bio_main_b}\\
\tau \dot{\bv{w}} &= - \bv{w} + \e^{\bv{r}},\label{eq:softmax_bio_main_c}
\end{align}
\end{subequations}
where $\1_{\np}$ is the all-ones vector and $\subscr{\tilde\tau}{g} \ll \subscr{\tau}{s} \ll \tau$, ensuring that the first two equations are faster than the last, output, dynamics.
Given $\bv{y}$, having components $y^\alpha$, the dynamics in Eqs.~\eqref{eq:softmax_bio_main_a}--\eqref{eq:softmax_bio_main_b} are exponentially stable and converge to the softmax normalizer $\sum_{\alpha=1}^{\np}\e^{y^\alpha}$ and to $\bv{y}-\1_{\np}\ln(m)$, respectively. The output layer in Eq.~\eqref{eq:softmax_bio_main_c} converges to the component-wise exponential of this last quantity. This value is $\softmax(\bv{y})$, which at steady state returns the optimal solution of Eq.~\eqref{eq:problem1}.
Details in Sec.~\ref{sec:methods}.

\subsection{Evaluation}

We evaluate our computational model on three paradigmatic domains: collective behaviors, human decision-making, and control tasks. The applications have been selected to identify and quantitatively measure model performance.
Experimental settings in Sec.~\ref{sec:methods}.

\subsubsection{Emerging Flocking Behaviors}

\begin{figure*}[t]
\centering
\includegraphics[width=1\linewidth]{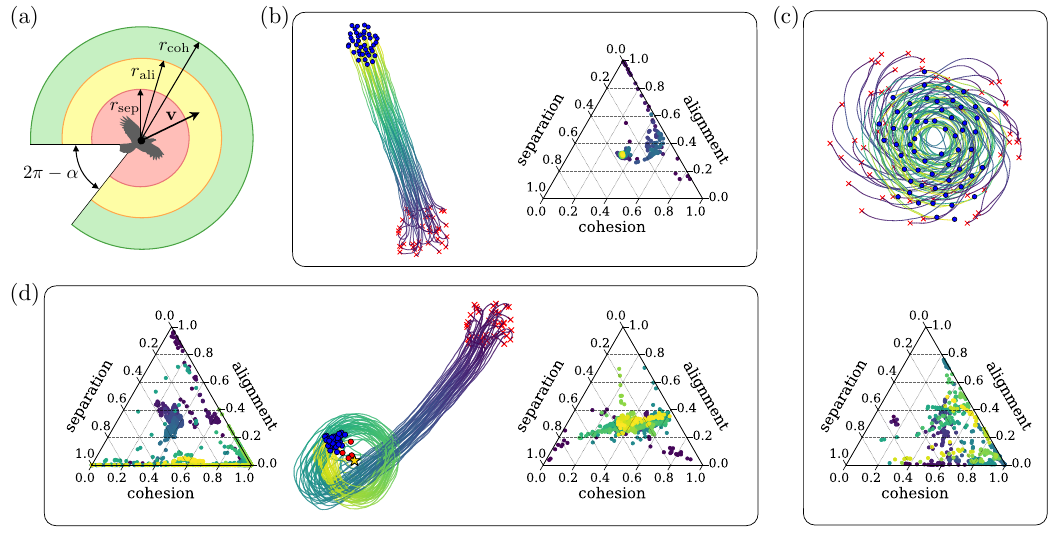}
\caption{Collective behavior experiments. (a) For the $i$-th boid in a flock, position and velocity components form the four-dimensional state $\bv{x}_k^i$. The action $\bv{u}_k^i$ is the acceleration vector, determined based on the state of the boid and of the neighbors that lie within its field of view. The field angle $\alpha$ is $320^\circ$. The radii correspond to three concentric separation, alignment and cohesion zones~\cite{IDC-JK-RJ-GDR-NRF:02,RL-YL-LEK:10}.
The innermost region is the separation zone, where a boid attempts to avoid collisions with neighbors within this area by applying a separation force. The alignment zone surrounds the separation zone. Through an alignment force, the boid attempts to achieve directional consensus with boids within this area. The outer zone is the cohesion zone. Through a cohesion force, the boid is drawn toward the average position of the boids in this area.  
(b) Each line is the trajectory of an individual and color gradients represent time.
Boids (filled circles) start from random initial positions (crosses).
Our model recovers polarization as an emerging collective behavior.
The evolution of the primitives/forces weights under our model is shown in the probability simplex. Over time, the weights tend to become uniformly distributed, linking a balanced use of the forces to the onset of polarization.
(c) Milling behavior, i.e., collective rotation around a common center, can be recovered when a collision avoidance cost is introduced in the generative model. This term prioritizes the separation force when the boid is too close to its neighbors.
(d) Middle: from random initial positions, the flock is soft-controlled~\cite{JH-ML-LG:06} toward the goal (star) when $10$\% of boids (leaders) have a goal-directed behavior. Our model shows that leaders (red-filled circles) use their primitives differently from the uninformed boids (blue-filled circles).
Followers (right) again lean toward a balanced use of social forces over time. Instead, leaders (left) exhibit more variable temporal modulation of social forces, alternating between balanced and more selective use of primitives.
See \SM, Table~S1 for simulation parameters.}
\label{fig:boids}
\end{figure*}

A central view in the study of collective animal behaviors is that these emerge from agent (e.g., boid)-level interactions, frequently expressed via social forces~\cite{CWR:87,IDC-JK-RJ-GDR-NRF:02,TV-AC-EBJ-IC-OS:95,JAC-MF-GT-FV:10,DM-DN-GP:26} that can emerge from surprise minimization~\cite{CH-BM-LDC-RPM-KJF-IDC:24}. Quantitative studies commonly infer the weights of these forces from data~\cite{RL-YL-LEK:10}. This approach yields strong agreement with observations but often limited explanatory insights that rely on interpreting constant fitted parameters.
Here, we interpret social forces as primitives and show that our computational model can dynamically modulate their use. At each time step $k$, the $i$-th boid (Fig.~\ref{fig:boids}(a)) determines its acceleration by sampling from a policy obtained by combining the social forces with weights from our model. We evaluate the ability of our model to recover and explain key behavioral signatures in leader and leaderless settings.

We first consider a leaderless setting. Inspired by previous works~\cite{CH-BM-LDC-RPM-KJF-IDC:24}, for the $i$-th boid the generative model comprises $\refplantboid{k}{k-1}{i}$ -- a multivariate Gaussian centered on the average position and velocity of neighboring boids within its cohesion and alignment zones, respectively -- and a uniform $\refpolicyboid{k}{k-1}{i}$. Our model not only recovers polarization consistently with the literature~\cite{CH-BM-LDC-RPM-KJF-IDC:24}, but also suggests a link between this collective behavior and weights tending to balance over time.
Fig.~\ref{fig:boids}(b) is representative of this phenomenon (see \SM, Fig.~S1 for the evolution of the weights of all boids). Milling is also recovered when a collision-avoidance term is included in the generative model (Fig.~\ref{fig:boids}(c)).
By promoting short-range repulsion over other social interactions, this term reproduces the interaction hierarchy of zonal models in which milling emerges~\cite{IDC-JK-RJ-GDR-NRF:02}.

Next, we consider a setting with leaders. Leaders are boids informed of a goal destination and their generative model encodes a goal-directed behavior. Fig.~\ref{fig:boids}(d) shows the emerging group behavior when a small fraction ($10\%$) of boids is goal-informed and seeks to reach the goal position.
Consistently with the literature~\cite{IDC-JK-NRF-SAL:05,DB-DJTS-JM-TG:06}, the group achieves goal-directed flocking without loss of cohesion. 
Moreover, experiments suggest that leaders tend to be more flexible than followers in their use of social forces, in the sense that the weights of the leaders exhibit higher temporal variability. 
This may indicate a more flexible behavior of leaders as increasingly emphasized in the literature~\cite{HL-GEM-JW-KVDV-RTV-AT-NTO:19,SB-RE-KJP:21,SN-MRM-MC-MP:17}. In \SM, Sec.~5 and Fig.~S2, we quantitatively assess the emerging behaviors in Fig.~\ref{fig:boids} via standard metrics.
Moreover, Fig.~S3 and Fig.~S4 report additional experiments with different numbers of informed boids, temperatures, and models for uninformed boids.
These experiments support our findings.

\subsubsection{Exploration/Exploitation in Human Decision-Making} 

\begin{figure*}[t]
\centering
\includegraphics[width=1\linewidth]{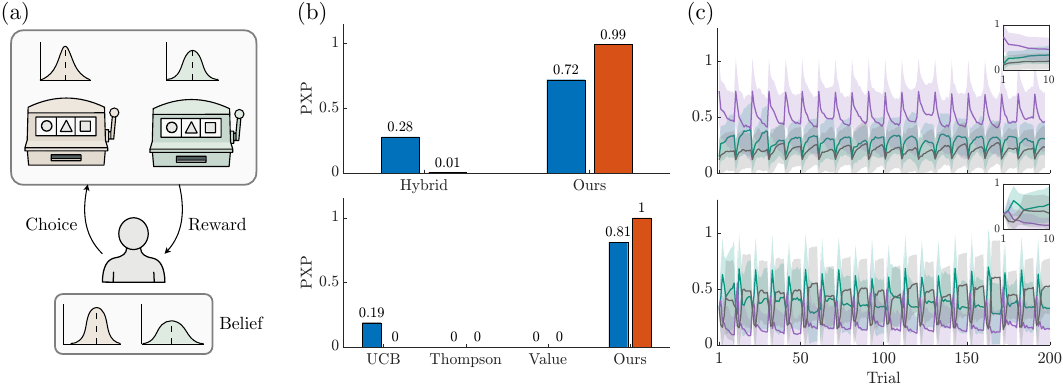}
\caption{Exploration/exploitation experiments.
(a) Participants performing two bandit experiments from~\cite{SJG:18}.
At each trial, a participant chooses one of two arms and receives a reward drawn from a distribution. Participants are instructed to maximize their total reward.
In Experiment~1, one arm returns stochastic rewards (positive or negative) while the other always returns zero. In Experiment~2, both arms return stochastic rewards. Over trials, participant $i$ updates a belief (mean and variance) of the reward for each arm~\cite{SJG:18}, forming the four-dimensional state $\mathbf{x}^i_{k-1}$. Based on this belief, the participant selects one of the two arms; the action $u^i_k=j$ corresponds to selecting arm $j\in\{1,2\}$. 
(b) PXP comparisons for Experiment~1 (blue) and~2 (orange) between: (i) the Hybrid model~\cite{SJG:18} and our model (top); (ii) UCB, Thompson, Value, and our model (bottom). PXP quantifies the probability that each considered model is the most frequent process that generated the data. Higher PXP for our model suggests that it provides better explanations for the data across experiments and benchmark models. 
(c) Our model suggests that participants behave differently across the experiments. These differences may be attributed to the different reward structures. Evolution of the mean weights (solid lines) of exploitation (teal), uncertainty-seeking (purple), risk-aversion (gray) behavioral schemas/primitives. Shaded areas indicate standard deviation across subjects. Trials from consecutive blocks are shown sequentially (see insets for a zoom on the first block).
In Experiment~1 (top), risk-aversion behavior remains low across trials, while uncertainty-seeking has the largest weight. In Experiment~2 (bottom), exploitation dominates early choices, after which the risk-aversion primitive contribution increases, while uncertainty-seeking remains comparatively small after a brief initial peak. See \SM, Table~S2 for parameter values.
}
\label{fig:mab}
\end{figure*}

A growing body of experimental and theoretical evidence suggests that human decision-making is better characterized by a combination of multiple mechanisms rather than a single strategy~\cite{SJG:18,AC:25}. For multi-armed bandit tasks, this view implies that human policies arise from composing decision rules. For example, human decisions are explained by a probabilistic model that combines algorithms based on uncertainty bonus (Upper Confidence Bound, UCB) and posterior sampling (Thompson)~\cite{SJG:18}. Using this probability as generative model and primitives capturing behavioral schemas, we evaluate the ability of our model to explain the data and suggest how schemas are used across  tasks.

We use data from the literature~\cite{SJG:18}  collected during two web-based two-armed bandit experiments (Fig.~\ref{fig:mab}(a)) involving $45$ and $44$ participants. Each experiment consists of $20$ blocks of $10$ trials. 
At the beginning of each block, the mean rewards of the arms are initialized by sampling from zero-mean Gaussian distributions. In Experiment~1 one arm deterministically returns zero reward, while the other produces stochastic rewards around a mean that is updated at each block. In Experiment~2 both arms return stochastic rewards, sampled from Gaussians with means updated between blocks.

In our model, a trial corresponds to a time step.
The state is the belief held by each participant about the mean and variance of the rewards associated with each arm, while the action is the arm pulled by the participant (Fig.~\ref{fig:mab}(a)). The generative model is the hybrid probabilistic model from~\cite{SJG:18}. Primitives encode simple behaviors: (i) exploitation, favoring the arm with the highest expected reward; (ii) uncertainty-seeking exploration, favoring the arm associated with the highest uncertainty; and (iii) risk aversion, favoring the least uncertain arm. At each trial our model returns 
a set of weights that modulate the use of these schemas. To quantify and benchmark the explanatory value of our model we use the Protected Exceedance Probability (PXP) also used in~\cite{SJG:18}.

First, we compare our model with the hybrid model~\cite{SJG:18}. Fig.~\ref{fig:mab}(b) (top) shows that our model achieves higher PXPs than the hybrid model on the same dataset for both experiments, thus suggesting better explanations of the data. The results are confirmed in Fig.~\ref{fig:mab}(b) (bottom), where our model also outperforms standard bandit algorithms.
Fig.~\ref{fig:mab}(c) reveals that individuals use behavioral schemas differently across the experiments and this difference may be attributed to the different structure of the rewards. In Experiment~1, participants tend to favor the most uncertain arm over the alternative that always returns zero. 
Instead, in Experiment~2 early choices in a block are consistently more exploitative, while  risk aversion increases over trials. See \SM, Fig.~S5 for additional experiments (cross-validation and experiments with different values of $\eps$) supporting the findings.

\subsubsection{Layered Control}

\begin{figure*}[t]
\centering
\includegraphics[width=1\linewidth]{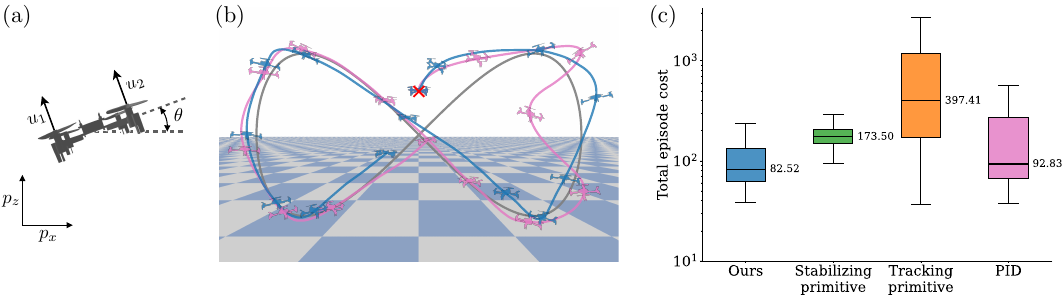}
\caption{Layered control experiments.
(a) For the drone, the state at time step $k$, $\bv{x}_k$, is six-dimensional: this is the vector of horizontal and vertical positions, pitch angle, and the corresponding linear and angular velocities. The control input is two-dimensional and each component represents the thrust applied to a rotor. 
(b) Tracking when our computational model is used (blue). Reference trajectory in grey; trajectory obtained with the baseline PID in pink. The start position is shown with a red cross and the translucent drone renderings show vehicle configurations during the experiment. 
(c) Comparison of the total episode cost between our model, the individual primitives and the baseline PID. The boxes around the median values represent interquartile ranges; whiskers denote the extrema excluding outliers. Statistics obtained from $1000$ experiments with random seeds. Renderings generated using \texttt{Safe-Control-Gym}. See \SM, Table~S3 for simulation parameters.}
\label{fig::quadrotor_pid}
\end{figure*}

Layered control architectures are common in biological~\cite{YN-QL-TS-JD:21} and engineered~\cite{TJP-SPW:23} systems. In these architectures, higher-level mechanisms coordinate the outputs of specialized lower-level controllers (primitives in our model) to produce complex behaviors.
Using a quadrotor tracking task and standard controllers from a benchmark environment~\cite{ZY-AWH-SZ-LB-MG-JP-APS:22} as primitives, we evaluate the ability of our model to orchestrate them, comparing performance against standard baselines and individual primitives.

The drone (Fig.~\ref{fig::quadrotor_pid}(a)) is a Crazyflie available from the \texttt{Safe-Control-Gym}~\cite{ZY-AWH-SZ-LB-MG-JP-APS:22} library. State and control variables are affected by Gaussian noise, capturing the presence of environmental and actuation disturbances. The reference is a lemniscate (a figure-eight-shaped curve) from the same environment, typically used to evaluate curvature tracking and stability.  
Primitives are two proportional-derivative (PD) controllers: 
(i) a stabilizing primitive that regulates the drone angle and damps the velocity, so that the drone stays around its current position without progressing along the reference, and (ii) a tracking primitive with control gains that do not allow for proper curvature following.

When the drone is equipped with our policy composition model, tracking is possible. Leveraging the flexibility offered by our normative framework in the choice of the generative model, we set $\refjointxu{k}{k-1}\propto \refjointxutwisted{k}{k-1} \e^{-c(\bv{x}_k,\bv{u}_k)}$, with $\refjointxutwisted{k}{k-1}$ being the product of $\tilderefplant{k}{k-1}$ set as the drone model available to the agent and a uniform $\tilderefpolicy{k}{k-1}$. The state-action cost $c$ is defined as the squared weighted tracking errors of state and action. Fig.~\ref{fig::quadrotor_pid}(b) shows the results from one experiment in this setting. 
The figure shows the trajectory obtained when the drone samples the control input from the policy returned by our model and the one obtained with a standard Crazyflie proportional-integral–derivative (PID) controller~\cite{JP-HZ-SZ-JX-APS:21}. Our model achieves lower median weighted cumulative error (i.e., cost) than the PID baseline even if the individual primitives behave poorly (Fig.~\ref{fig::quadrotor_pid}(c)). 
Differences in total episode cost are statistically significant under a two-sided Mann–Whitney $U$ test ($p<10^{-7}$).
Our model also achieves reduced variability across runs, indicating more consistent performance than the baseline.
See \SM, Sec.~5 and Fig.~S7 for additional experiments with a different set of primitives, which support the findings.

\section{Discussion}
This work introduces a theoretically grounded mechanism for policy composition. Gating emerges from the minimization of a {variational objective combining divergence from a generative model and entropy regularization.} The optimal weights arise as the equilibrium of an energy model that features explicit convergence rates and admits a neural instantiation.
Comprehensive evaluations across diverse applications suggest that our model provides a general framework for understanding the computational principles of policy composition and the neural circuits implementing this function.
The results also highlight context-dependent synaptic couplings as a principled computational necessity. 
Bridging data-driven and first-principles perspectives, the findings yield insights into policy composition in nature and can pave the way for more biologically grounded artificial agents.

The results also open avenues for further exploration and improvement.
In this work we explain the onset of softmax gating rules, which are widely used but not the only options.
The results hint that it may be possible to obtain other rules by adopting different regularizations in the normative framework. For example, sparsity-promoting regularizers yield biologically plausible dynamics~\cite{VC-AG-AD-GR-FB:23a} for reconstruction problems featuring rectifying linear units in the flow. Although our study focuses on linear primitive combinations, the results extend to broader settings in which the statistical complexity is {still} convex in the weights.
Such settings include nonlinear combinations of primitives, opening a principled route toward product-of-experts–type architectures. Primitives are inputs to our model. An open question is how to learn sets of effective primitives for a given task. Addressing this question may require augmenting our model with learning mechanisms in the spirit of, e.g., hierarchical learning schemes. Our results suggest that a possible way to derive a learning scheme is to exploit the links between our model and active inference~\cite{TP-GP-KF:22}. Inspired by social and network sciences, introducing a feedback reward signal from the gating mechanism to the primitives might be a promising way to implement learning and primitives' innovation~\cite{II_SM_VL:18,GR-QS-LW-JP:25}.
This signal may yield a dynamics for the synaptic weights in our circuit, whose convergence should be assessed~\cite{VC-FB-GR:22k}. In multi-agent settings, where agents may themselves be primitives and evolve their capabilities over time, a central question is whether our model yields a form of collective intelligence~\cite{DC:22} and under what conditions.
These limitations and open questions may inspire future interdisciplinary work.

\section{Materials and Methods}\label{sec:methods}

We recall that the agent has access to: (i) the generative model, $\refjointxu{k}{k-1} = \refplant{k}{k-1}\refpolicy{k}{k-1}$; (ii) the model $\plant{k}{k-1}$; (iii)  $\np$ primitives, $\primitive{k}{k-1}{\alpha}$, $\alpha=1,\ldots, \np$, with support spanning the action space. In \SM~ we unpack the formal treatment and show that the properties derived here also apply when a bias vector for the weights is introduced in our computational model.

\subsection{Optimization Properties}\label{sec:modelframe_methods} In Eq.~\eqref{eq:problem1}, the decision variables $\weights{k}$ belong to the simplex. Therefore, the term $-\entropy{(\weights{k})}$ in the cost is strictly convex and the constraints are convex in the decision variables.
We now show that the statistical complexity term in the cost, $\F(\weights{k})$, is  convex in $\weights{k}$, thus making the overall problem strictly convex. To show this, we embed the constraints into the expression of $\F(\weights{k})$. Then, using the chain rule for the KL divergence reveals that the complexity term can be written as $
\F(\weights{k})=\sum_{\alpha =1}^{\np}\weight{k}{\alpha} \big(\E_{\primitive{k}{k-1}{\alpha}}\big[ \ln\sum_{\beta =1}^{\np}\weight{k}{\beta}\primitive{k}{k-1}{\beta}  - \ln\refpolicy{k}{k-1}  + \DKLb{\plant{k}{k-1}}{\refplant{k}{k-1}}\big] \big)$. Here, $\weight{k}{\alpha}$ is the element of the vector $\weights{k}$ corresponding to primitive $\alpha$. The function $\F(\weights{k})$ is twice differentiable in $\weights{k}$. The Hessian of $\F(\weights{k})$ is a positive semi-definite matrix, proving that the map $\F(\weights{k})$ is convex. 

\subsection{Why the Framework Yields the Softmax Gradient}
Our starting point is the reformulation in Eq.~\eqref{eq:entropy_regularized_problem_simplex_weights}. To  derive the results, we leverage the continuous-time proximal gradient flow~\cite{BA-HA:14, SHM-MRJ:21, AD-VC-AG-GR-FB:23f}. Consider a generic composite optimization problem of the form $\min_{\bv{x} \in \R^n} f(\bv{x}) + g(\bv{x})$, where $f$ is a smooth function and $g$ closed, convex, proper (CCP) and possibly non-smooth. The continuous-time proximal gradient method  updates an estimate of the optimal solution via the dynamics
\begin{equation*}
 \bv{\dot x} = - \bv{x} + \prox{\gamma g}{\bigl(\bv{x} - \gamma \nabla f(\bv{x})\bigr)},
\end{equation*}
where $\prox{\gamma g}(\bv{x}) :=  \argmin_{\bv{z} \in \R^n} g(\bv{z}) + \frac{1}{2 \gamma}\|\bv{x} - \bv{z}\|_2^2$, for all $\bv{x} \in \R^n$ is the proximal operator of $g$ and $\gamma$ is a parameter~\cite{NP-SB:14}. The softmax gradient flow 
is the continuous-time proximal gradient flow associated to Eq.~\eqref{eq:entropy_regularized_problem_simplex_weights}, and hence to our normative framework optimization. Namely, the softmax gradient flow (Eq.~\eqref{eq::FreeNet}) is obtained from  reformulating Eq.~\eqref{eq:entropy_regularized_problem_simplex_weights} as
$$
\min_{\weights{k} \in \R^{\np}} f(\weights{k}) + g(\weights{k}),
$$
where ${f(\weights{k}) = \F(\weights{k}) + \eps \frac{\norm{\weights{k}}^2}{2}}$ is continuously differentiable and $g(\weights{k}) = \eps\bigl(\entropybar(\weights{k}) - \frac{\norm{\weights{k}}^2}{2}\bigr)$ is CCP. The $\softmax$ in the softmax gradient flow is (\SM, Sec.~3) $\prox{\eps^{-1} g}$.

\subsection{Softmax Gradient Flow Properties}\label{sec::softmax:gf} We first derive the forward invariance of the simplex $\simplex{\np}$ and then characterize convergence. 

Intuitively, forward invariance of the simplex means that when the dynamics is initialized at some initial condition $\weights{k}(0) \in \simplex{\np}$, then its trajectories remain in the simplex for all $t$, that is, $\weights{k}(t)\in\simplex{\np}$, for all $t >0$. This desirable property can be shown by proving that: (i) the positive orthant is forward invariant; (ii) the total {\em mass} is conserved, that is, $\sum_{\alpha = 1}^{\np} \frac{d{\weight{k}{\alpha}}}{dt} = 0$ for all $\weights{k} \in \R^{\np}$ such that $\sum_{\alpha = 1}^{\np} \weight{k}{\alpha}=1$.
The first property can be shown via Nagumo's theorem; see, e.g.,~\cite[Exercise~3.12]{FB:24-CTDS}. In particular, each component of the softmax gradient vector field, i.e., the right-hand side of Eq.~\eqref{eq::FreeNet}, is non-negative on the boundary of the positive orthant.
Therefore, trajectories on the boundary are subject to a vector field pointing toward the interior of the positive orthant, ensuring that they remain in the orthant. The second property follows from direct computation. Summing all the components of the softmax gradient yields 
$\sum_{\alpha=1}^{\np}\bigl( \softmax\bigl(- \eps^{-1}\nabla \F(\weights{k})\bigr)_\alpha - \weight{k}{\alpha}\bigr) = 0$,
as desired.

Detailed derivations for the dissipation identity (Eq.~\eqref{eq:dissipation_inequality_main}) are in \SM, Sec.~3. From Eq.~\eqref{eq:dissipation_inequality_main}, exponential convergence of the energy function to the energy minimum $V(\optimalweights{k})$ is established via convexity of $\F$ and Gr\"onwall's inequality. This yields the convergence bound:
\begin{equation}
V(\weights{k}(t))-V(\optimalweights{k}) \le \e^{-t/\tau}\bigl(V(\weights{k}(0))-V(\optimalweights{k})\bigr),
\end{equation}
for all $t\ge 0$. Convergence of $\weights{k}(t)$ to the optimal weights $\optimalweights{k}$ is subsequently established via convexity of $\F$ and Pinsker's inequality, yielding
\begin{equation}
\|\weights{k}(t)-\optimalweights{k}\| \le \sqrt{\frac{2}{\eps}\bigl(V(\weights{k}(0))-V(\optimalweights{k})\bigr)}\,\e^{-t/(2\tau)},
\end{equation}
for all $t \ge0$, with $\norm{\cdot}$ being the Euclidean norm.

\subsection{Deriving the Circuit} 
The fast unit (Fig.~\ref{fig:figure2}(c) and Eq.~\eqref{eq:gradient_bio}) returns the argument of the $\softmax$. This argument is obtained from the gradient of $\F(\weights{k})$, which we recall is the statistical complexity in the optimization. The negative gradient of $\F(\weights{k})$ is a $\np$-dimensional vector. The component $\alpha$ of this vector is
\begin{equation*}
-\E_{\primitive{k}{k-1}{\alpha}}\left[\ln \weights{k}^\top\primitives{k}{k-1} + \bv{s}(\bv{x}_{k-1} , \bv{u}_{k})\right] - 1.
\end{equation*}
From this expression, the first step to obtain the fast unit dynamics in Eq.~\eqref{eq:gradient_bio} is to note that the constant term can be dropped due to the translation invariance property of the softmax. In the discrete-action setting, this expression yields the component $\alpha$ of the vector $\bar{\bv{y}}_k$ from the results.
The second step is to rewrite this vector into a form that is suitable for a neural dynamics representation. To this aim,
$\bar{\bv{y}}_k$ can be conveniently written (dropping subscripts in the network state variables) as $-\eps^{-1}\Pi(\bv{x}_{k-1})^{\top}(\bar{\bv{b}} + \bv{s})$, where we introduced the change of variables $\bar{\bv{b}}:=\ln(\bar{\bv{a}})$ (component-wise), with $\bar{\bv{a}}=\Pi(\states{k-1})\bv{w}$.
The slow unit (Fig.~\ref{fig:figure2}(c) and Eq.~\eqref{eq:softmax_bio_main}) is adapted from the literature~\cite{MS-JO:22} and implements the $\softmax$ by leveraging the identity: $\softmax(\bv{y})_{\alpha} =  \e^{y^{\alpha} - \ln(m)}$.

\subsection{Experimental Settings}
Detailed settings for all experiments are provided in \SM, Sec.~5.
In all experiments, primitive weights are initialized uniformly on the simplex.
In the collective behavior experiments, velocities and accelerations are bounded. As in previous works~\cite{CH-BM-LDC-RPM-KJF-IDC:24}, the boid dynamics $\plantboid{k}{k-1}{i}$ is a Gaussian centered in the state provided by a second-order model~\cite{FC-SS:07, HL-WJR-IC:00}. The expressions for the separation, alignment, and cohesion primitives are in accordance with prior literature. Initial positions are sampled randomly near the origin, initial velocities are random.
Consistent with previous literature~\cite{CH-BM-LDC-RPM-KJF-IDC:24}, social-forces models can yield cases in which, depending on initial conditions and parameter values, cohesion of the group is not maintained;
as common in this literature, in this work we focus on cohesive behaviors.
In the multi-armed-bandit experiments, the state of participant $i$, $\mathbf{x}_{k-1}^i$, is the stack of the estimated posterior means and variances of each arm after observing the reward at trial $k-1$ and before selecting an arm at trial $k$.
Belief/state updating is as in prior work~\cite{SJG:18}, leveraging a Kalman filter and with beliefs being reset at the start of each block.
We use the original code from the literature~\cite{SJG:18} to compute posterior estimates and to implement, via probit regression, the Hybrid, UCB, Thompson and Value policies.
See \SM, Fig.~S6 for additional analyses using an alternative to probit.
In the layered control experiments, the drone model, reference trajectory and controller structures are from the \texttt{Safe-Control-Gym}. The environment enables the injection of Gaussian noise in the dynamics/actions.
Consequently, the primitives are Gaussians centered in the output of the PD controllers.

\section*{Acknowledgments}
Apple Keynote was used for drawing schematics and assembling figure panels.
We acknowledge the use of ChatGPT--5.5 for assistance with improving wording and grammar of this document. All outputs were carefully reviewed and verified by the authors. No AI tools were used to produce scientific findings or perform analysis.

V.C. was at University of Salerno while preparing the theoretical bulk of this work. V.C. and G.R. were supported in part by the European Union-Next Generation EU Mission 4 Component 1 CUP E53D23014640001. F.B. was supported in part by AFOSR grant FA9550-22-1-0059.

F.R. and V.C. contributed equally to this work. F.B. and G.R. contributed equally to this work.
F.R., F.B., and G.R. designed the research; F.R. and V.C. performed the theoretical analysis with inputs from F.B. and G.R.; F.R. developed the code for the experiments and obtained the numerical results; F.R., V.C., F.B., and G.R. analyzed the results; F.R., V.C., F.B., and G.R. wrote the paper.

There are no competing interests to declare.

\section*{Data Availability}

The human multi-armed-bandit data analyzed in this study are publicly available from the repository associated with Ref.~\cite{SJG:18}.
The custom code used for the simulations and data analyses in this article is openly available on GitHub~\cite{FR:26_repo}.
The multi-armed-bandit analysis uses code adapted from the repository associated with Ref.~\cite{SJG:18}. The layered-control experiments use the \texttt{Safe-Control-Gym} environment from Ref.~\cite{ZY-AWH-SZ-LB-MG-JP-APS:22}.

\bibliography{bib_arxiv}

\end{document}


\setcounter{section}{0}
\renewcommand{\thesection}{\arabic{section}}
\renewcommand{\thesubsection}{\thesection.\arabic{subsection}}
\renewcommand{\thesubsubsection}{%
  \thesubsection.\arabic{subsubsection}}

\setcounter{equation}{0}
\renewcommand{\theequation}{S\arabic{equation}}

\setcounter{figure}{0}
\renewcommand{\thefigure}{S\arabic{figure}}

\setcounter{table}{0}
\renewcommand{\thetable}{S\arabic{table}}

\renewcommand{\theHequation}{supp.equation.\arabic{equation}}
\renewcommand{\theHfigure}{supp.figure.\arabic{figure}}
\renewcommand{\theHtable}{supp.table.\arabic{table}}

\title{\textbf{Supplemental Material for\\ An Energy-Based Mechanism for Compositional Behavior} 
}%

\author{Francesca Rossi}
\email{Contact author: f.rossi@ssmeridionale.it}
\affiliation{Scuola Superiore Meridionale, Naples, Italy.}

\author{Veronica Centorrino}
\affiliation{Automatic Control Laboratory (IfA), ETH Z\"urich, Switzerland.}

\author{Francesco Bullo}
\affiliation{Center for Control, Dynamical Systems, and Computation,
University of California, Santa Barbara, California, USA.}

\author{Giovanni Russo}
\affiliation{Department of Information and Electrical Engineering and Applied Mathematics, University of Salerno, Italy.}

\maketitle

\section{Introduction}

After providing some background in Section~\ref{sec:math_preliminaries},
we consider a general class of entropy-regularized minimization problems and obtain a dynamical system that provably converges to the optimal solution of such problems (Section~\ref{sec:general_results}). 
Building on these general findings, we develop (Section~\ref{sec:model_details}) the formal statements for the results in the main text. After reporting the supplemental details of the experiments (Section~\ref{sec:experiments_details})
we provide the proofs of all the statements (Section~\ref{sec:proofs}). Finally, in Section~\ref{sec:figures} and Section~\ref{sec:tables} we report supplemental figures and tables, respectively.

\section{Background}
\label{sec:math_preliminaries}
We introduce notation, definitions, and standard results related to Kullback-Leibler (\KL) divergence and entropy. We also briefly survey key elements of convex optimization and proximal operator theory relevant to our analysis. 

\paragraph*{\textbf{Notation.}}
Sets are in {\em calligraphic} letters and vectors in {\bf bold}. The set $\displaystyle \simplex{n} := \Bigsetdef{\bv{p} \in \R^n}{\sum_{i = 1}^{n} p_i=1, p_i\geq 0, \ \forall i \in \{1,\ldots, n\}}$ is the \emph{probability simplex} in $\R^n$. We denote by $\0_n$ and $\1_n$ the all-zeros and all-ones vector in $\R^n$, respectively. We let $I_n$ denote the $n$-dimensional identity matrix.
For a vector $\bv{x} \in \R^n$, we denote by $\|\bv{x}\|_p := \left(\sum_{i=1}^n |x_i|^p\right)^{1/p}$ the $\ell_p$-norm for $p \ge 1$. When the subscript is omitted, $\|\cdot\|$ denotes the Euclidean norm, i.e., $\|\cdot\|_2$.
Random variables are denoted via capital letters and their realizations by lowercase letters. The {probability density function} (pdf) of $\mathbf{V}$ is denoted by $p(\mathbf{v})$; for discrete variables, $p(\mathbf{v})$ is the {probability mass function} (pmf).
The expectation of a function $h(\cdot)$ with respect to a random variable $\mathbf{V}$ is defined as $\ds \E_{{p}}[h(\mathbf{V})]:=\int_{\mathcal{S}(p)}h(\mathbf{v})p(\mathbf{v})d\bv{v}$, where $\mathcal{S}(p)$ denotes the support of $p(\mathbf{v})$. For discrete variables, the integral is replaced by a sum. When clear from context, we omit the domain of integration or summation.
The joint pdf/pmf of two random variables $\mathbf{V}_1$ and $\mathbf{V}_2$ is  $p(\mathbf{v}_1,\mathbf{v}_2)$ and the conditional pdf/pmf of $\mathbf{V}_1$ with respect to (w.r.t.) $\mathbf{V}_2$ is $p\left( \mathbf{v}_1\mid \mathbf{v}_2 \right)$. Given two pdfs/pmfs $p(\mathbf{v})$ and $q(\mathbf{v})$, we say that $p(\mathbf{v})$ is \emph{absolutely continuous} with respect to $q(\mathbf{v})$ if the support of $p(\bv{v})$ is contained in the support of $q(\bv{v})$. We adopt the standard convention $0 \ln(0) = 0$. The multivariate Gaussian distribution with mean $\boldsymbol{\mu} \in \mathbb{R}^d$ and covariance matrix $\Sigma \in \mathbb{R}^{d \times d}$ is
$$
\mathcal{N}(\boldsymbol{\mu}, \Sigma)  = \frac{1}{(2\pi)^{d/2} |\Sigma|^{1/2}} \exp \Bigl(-\frac{1}{2} (\bv{v} - \boldsymbol{\mu})^\top \Sigma^{-1} (\bv{v} - \boldsymbol{\mu})\Bigr).
$$

\subsection*{The Softmax Function}
The {softmax function} is the map $\map{\softmax}{\R^n}{\R^n}$ defined by $\softmax(\bv{x}) {=} \bigl(\softmax(\bv{x})_i\bigr)$, $i \in \until{n}$, where
$$
\softmax(\bv{x})_i := \ds \frac{\exp(x_i)}{\sum_{j= 1}^n \exp(x_j)}.
$$
This function provides a differentiable and continuous approximation of the argmax operator~\cite{IG-YB-AC:16} and naturally arises in the context of statistical physics (here it is also known as \emph{Boltzmann} or \emph{Gibbs distribution}~\cite{RSS-AGB:98}), economics and game theory (where it is often referred to as \emph{logit response function}, \emph{logit map}, or \emph{perturbed best response}~\cite{RDM-TRP:95, RC-EM-SS:10, PC-BG-PM:15, DSL-EJC:05}),  and  machine learning~\cite{JP-KM-YA:10, BG-LP:17, WHS:10, IG-YB-AC:16, LK-KVK-DK:23}. In neuroscience, this function  emerges, for example, in Winner-Take-All (WTA) networks to promote competition among neurons~\cite{IME-JLW:93, ALY-DG:95}. Finally, we recall the following.
\begin{Lemma}[Softmax translation invariance]
\label{rem:softmax_sum}
Let $\F$ be a function of the form $\F(\bv{x}) = \tilde \F(\bv{x}) + u\1_n^\top \bv{x}$, for some $u \in \R$. Then,
$$
\softmax(\nabla \F(\bv{x})) = \softmax(\nabla \tilde \F(\bv{x}) + u \1_n) = \softmax(\nabla \tilde \F(\bv{x})).
$$
\end{Lemma}

\subsection*{The Kullback-Leibler Divergence and Entropy}
We recall the definition of \KL~divergence, which provides a measure of discrepancy between two probability distributions.
\bd[Kullback-Leibler Divergence~\cite{SK-RAL:51}]
Given two probability distributions $p(\mathbf{v})$ and $q(\mathbf{v})$, the \emph{Kullback-Leibler divergence of $p$ from $q$} is
$$
\DKL{p}{q}:= \int_{\mathcal{S}(p)} p(\bv{v}) \ln\biggl( \frac{p(\bv{v})}{q(\bv{v})}\biggr)d\bv{v}.
$$
\ed
The $\DKL{p}{q}$ is finite only if $p(\bv{v})$ is absolutely continuous w.r.t. $q(\bv{v})$~\cite[Chapter $8$]{TC-JT:06}.
The following result, see, e.g.,~\cite[Theorem 2.5.3]{TC-JT:06}, is used in the main text.
\begin{Lemma}[Chain rule for the \KL~divergence]
\label{lem:chain_rule_KL}
Let $p\left(\mathbf{v}, \mathbf{u}\right)$ and $q\left(\mathbf{v}, \mathbf{u} \right)$ be joint probability density functions. Then, 
\[
\DKL{p\left(\mathbf{v}, \mathbf{u}\right)}{q\left(\mathbf{v}, \mathbf{u}\right)} =
\DKL{p(\mathbf{u})}{q(\mathbf{u})} + \E_{p(\mathbf{u})}\left[\DKL{p\left(\mathbf{v} \mid \mathbf{u}\right)}{q\left(\mathbf{v} \mid \mathbf{u}\right)} \right].
\]
\end{Lemma}
\bd[Entropy and Cross-Entropy]
The \emph{entropy function} is $\map{\entropy}{\simplex{n}}{[0, \ln n]}$, $\displaystyle \entropy(p):= - \sum_{i=1}^{n}{p_i}\ln{p_i}$.
The \emph{cross-entropy of $q$ relative to $p$} 
is  $\map{\entropy}{\simplex{n} \times \simplex{n}}{[0, +\infty]}$, $\displaystyle \entropy(p,q):= - \sum_{i=1}^{n}{p_i}\ln{q_i}$.
\ed
\begin{pro}[Relation between \KL~divergence and entropy~\cite{AU:96}]
\label{prop:cross_entropy}
Let $p(\mathbf{v})$ and $q(\mathbf{v})$ be two probability mass/density functions. The following identity holds: 
$\DKL{p}{q}=\entropy(p,q)-\entropy(p)$.
\end{pro}

\subsection*{Convex Analysis and Proximal Gradient Dynamics}
We recall the following standard definition.
\begin{definition}[$L$-Lipschitz function]
A map $\map{f}{\R^n}{\R}$ is \emph{Lipschitz with constant $L$ ($L$-Lipschitz)} if for all $\bv{x}_1, \bv{x}_2 \in \R^n$, it holds that
$$
\abs{f(\bv{x}_1) - f(\bv{x}_2)} \leq L \norm{\bv{x}_1 - \bv{x}_2}.
$$
\end{definition}
A map $\map{g}{\R^n}{\realextended}:= [{-}\infty,+\infty]$, is (i) \emph{convex} if $\operatorname{epi}(g) := \setdef{(\bv{x},y) \in \R^{n+1}}{g(\bv{x}) \leq y}$ is a convex set (with $\operatorname{epi}$ denoting epigraph); (ii) \emph{proper} if its value is never $-\infty$ and there exists at least one $\bv{x} \in \R^n$ such that $g(\bv{x}) < \infty$; (iii) \emph{closed} if it is proper and $\operatorname{epi}(g)$ is a closed set. We denote by $\partial g$ the \emph{subdifferential} of $g$. Moreover, we say that a map $g$ is 
\begin{enumerate}[label = (\roman*)]
\item \emph{strictly convex} if for all $\bv{x}, \bv{y} \in \R^n$ with $\bv{x} \neq \bv{y}$, and for all $\lambda \in (0,1)$, it holds that $g(\lambda \bv{x} + (1-\lambda)\bv{y}) < \lambda g(\bv{x}) + (1-\lambda) g(\bv{y})$;
\item \emph{strongly convex with parameter $\rho > 0$} if the map $\ds \bv{x} \mapsto g(\bv{x}) - \frac{\rho}{2}\|\bv{x}\|_2^2$ is convex;
\item {\emph{$L$-smooth} if it is differentiable and $\nabla g$ is $L$-Lipschitz.}
\end{enumerate}

Next, we define the proximal operator of $g$, which maps a point $\bv{x} \in \R^n$ into a subset of $\R^n$. This subset can be either empty, contain a single element, or be a set with multiple vectors.
\bd[Proximal Operator]
\label{apx:def:prox_operator}
Let $\map{g}{\R^n}{\realextended}$ be a closed, convex, and proper function, and let $\gamma > 0$. The \emph{proximal operator} of $g$ with parameter $\gamma$ is the map $\map{\prox{\gamma g}}{\R^n}{\R^n}$ defined by
\begin{equation*}
 \prox{\gamma g}(\bv{x}) = \displaystyle \argmin_{\bv{z} \in \R^n} g(\bv{z}) + \frac{1}{2 \gamma}\|\bv{x} - \bv{z}\|_2^2, \quad \textup{ for all } \bv{x} \in \R^n.
\end{equation*}
\ed
For any closed, convex, and proper (CCP) function $g$, the proximal operator $\prox{\gamma g}(\bv{x})$ exists and is unique for all $\bv{x} \in \R^n$~\cite[Theorem 6.3]{AB:17}.
Based on the use of proximal operators, {\em proximal gradient method} (see, e.g.,~\cite{NP-SB:14}) can be devised to iteratively solve a class of composite (possibly non-smooth) convex problems of the form
\begin{equation*}
\min_{\bv{x} \in \R^n} f(\bv{x}) + g(\bv{x}),
\end{equation*}
where $\map{f}{\R^n}{\R}$, $\map{g}{\R^n}{\realextended}$ are CCP functions, and $f$ is differentiable. At its core, the proximal gradient method updates the estimate of the solution of the optimization problem by computing the proximal operator of $\alpha g$, where $\alpha > 0$ is a step size, evaluated at the difference between the current estimate and the gradient of $\alpha f$ computed at the current estimate. This method has been extended and generalized to a continuous-time framework~\cite{BA-HA:14, SHM-MRJ:21}, resulting in the \emph{continuous-time proximal gradient dynamics}
\begin{equation}
\label{eq:prox_gradient}
 \bv{\dot x} = - \bv{x} + \prox{\gamma g}{\bigl(\bv{x} - \gamma \nabla f(\bv{x})\bigr)}, \ \ \ \gamma >0.
\end{equation}

Recent work~\cite{VC-AG-AD-GR-FB:23a} has provided an alternative interpretation of the dynamics in Eq.~\eqref{eq:prox_gradient}, particularly in the case where $f$ is a quadratic function and $g$ is a separable sum of scalar functions.

\section{Entropy-regularized Optimization and Proximal Gradient Flows}
\label{sec:general_results}
We now start developing the formal treatment and statements for the results in the main text. To this end, we consider a general class of entropy-regularized minimization problems of the form
\begin{equation}
\label{eq:entropy_regularized_problem}
\min_{\bv{x} \in \simplex{n}} \F(\bv{x}) - \eps \entropy(\bv{x}),
\end{equation}
with $\map{\F}{\R^n}{\R}$ being convex and continuously differentiable, and $\eps > 0$ being a regularization (\emph{temperature}) parameter. This setting embeds as a special case the optimization problem from our normative framework. Our goal is to obtain a continuous-time dynamics that provably converges to the optimal solution of the problem in Eq.~\eqref{eq:entropy_regularized_problem}. The formal derivations
for the normative framework, the softmax gradient flow, and its neural instantiation in the main text build on the results reported here.

Entropy-regularized optimization problems of the form~\eqref{eq:entropy_regularized_problem} naturally arise across reinforcement learning, game theory, optimal transport, and statistical physics~\cite{PM-WHS:16, VM-APB-MM-AG-TL-TH-DS-KK:16,GP-MC:19,AS-HJ-KF-GR:25_new}. This class of problems also appears in the form of smooth best response maps, which leads to the logit (i.e., softmax) function as a solution~\cite{PC-BG-PM:15, PM-WHS:16, BG-LP:17}.
Finally, optimization over the probability simplex with entropy regularization is also closely related to mirror descent, where the negative entropy acts as the mirror map~\cite{AB-MT:03,AB:17}.

\begin{Remark}[Existence and uniqueness of the minimizer of problem~{\eqref{eq:entropy_regularized_problem}}]
\label{remark:exist_uniq}
Consider the entropy-regularized problem in Eq.~\eqref{eq:entropy_regularized_problem}. Since the set $\simplex{n}$ is compact and the map $x_i\mapsto - x_i\ln x_i$ extends continuously to $x_i=0$ (under the standard convention $0 \ln(0) = 0$), the objective function is continuous on $\simplex{n}$. Hence, by the Weierstrass Theorem, there exists a minimizer.
Regarding uniqueness, note that the gradient of the regularization term $- \eps \entropy(\bv{x})$ satisfies $\partial (\eps \sum_{j=1}^n x_j\ln x_j)/\partial x_i \to -\infty$ as $x_i\to 0^+$.
Therefore, no minimizer can lie on the boundary of the simplex, and every minimizer must belong to the relative interior of $\simplex{n}$ where the map $- \eps \entropy(\bv{x})$ is strictly convex.
Since $\F$ is convex, the cost function $\F(\bv{x}) - \eps \entropy(\bv{x})$ is strictly convex on the relative interior of $\simplex{n}$, implying uniqueness of the minimizer.
\end{Remark}

Equipped with this observation, next we derive a dynamics provably converging to the optimal solution of the problem in Eq.~\eqref{eq:entropy_regularized_problem}. We first propose a reformulation of the optimization problem that enables computing a proximal operator. Then, we obtain the corresponding continuous-time proximal dynamics, and characterize its properties. Namely, we show that the equilibrium of this dynamics is the optimal solution of the problem in Eq.~\eqref{eq:entropy_regularized_problem} and that trajectories provably converge to such equilibrium.

\subsection{Reformulating Problem~\eqref{eq:entropy_regularized_problem}}
We introduce an unconstrained reformulation of the problem in Eq.~\eqref{eq:entropy_regularized_problem}. This is obtained by defining the \emph{entropic barrier function} $\map{\entropybar}{\R^{n}}{]-\infty, + \infty{]}}$ as 
$$
{
\entropybar(\bv{x}) := 
\begin{cases}
- \entropy(\bv{x}) & \text{if } \bv{x} \in \simplex{n},\\
+\infty & \text{otherwise}.
\end{cases}
}
$$
This function is equal to the {negative} entropy on the simplex, the feasibility domain of the cost function in Eq.~\eqref{eq:entropy_regularized_problem}, and infinity outside of this set. Note that $\entropybar$ is CCP.
Problem~\eqref{eq:entropy_regularized_problem} can be conveniently recast as
\begin{equation}
\label{eq:entropy_regularized_problem_simplex}
\min_{\bv{x} \in \R^n} \F(\bv{x}) + \eps \entropybar(\bv{x}).
\end{equation}
To obtain the proximal gradient flow associated with the problem in Eq.~\eqref{eq:entropy_regularized_problem_simplex}, we add and subtract $\eps \frac{\norm{\bv{x}}^2}{2}$ to the objective function. 
This motivates the following result, in which we compute the proximal operator of $\entropybar(\bv{x}) - \frac{\norm{\bv{x}}^2}{2}$. While this result can be proved via an adaptation of~\cite[Example 2.23]{PLC-JCP:20b}, where tools from monotone operator theory are used, in Section~\ref{sec:proofs} 
we give an alternative proof, which is interesting {\em per se}. This self-contained proof is based directly on the definition of the proximal operator.
\begin{Lemma}\label{lem:prox_softmax}
Given the function $\map{\entropybar}{\R^n}{{]-\infty, + \infty]}}$, we have
$
\displaystyle \prox{ \entropybar - \frac{\norm{\cdot}^2}{2}}(\bv{x}) = \softmax(\bv{x}).
$
\end{Lemma}
Beyond being instrumental in obtaining a dynamical system that converges to the optimal solution of the problem in Eq.~\eqref{eq:entropy_regularized_problem}, Lemma~\ref{lem:prox_softmax} also reveals a connection between proximal operators and smooth best response maps.
In particular, Eq.~\eqref{eq:lem_prox_entropy2} in the proof of Lemma~\ref{lem:prox_softmax} shows that the smooth best response map is the proximal operator of the function $\entropybar - \frac{\norm{\cdot}^2}{2}$.
This follows directly by observing that Eq.~\eqref{eq:lem_prox_entropy2} can be equivalently written as
$\ds \argmax_{\bv{z} \in \simplex{n}} \left(\bv{x}^\top \bv{z} - \sum_{i=1}^n z_i \ln z_i \right)
= \argmax_{\bv{z} \in \simplex{n}} \left( \bv{x}^\top \bv{z} + \entropy(\bv{z}) \right)$.
To the best of our knowledge, this observation does not appear explicitly in the existing literature.

This proximal interpretation also aligns naturally with the classical derivation of smooth best response maps in game theory. As discussed in~\cite{PC-BG-PM:15, BG-LP:17}, in the context of game theory, the softmax function can be obtained by considering the \emph{argmax function} under entropy regularization resulting in what is called \emph{smooth best response map}.
Specifically, let $\bv{x} \in \R^n$ and consider the argmax function, or best response map (the agent seeks the strategy with the highest score), given by
$$
\bv{x} \mapsto \argmax_{\bv{z} \in \simplex{n}} \bv{x}^\top \bv{z}.
$$
This best response map carries several challenges (see, e.g.,~\cite{PC-BG-PM:15, BG-LP:17} for more details). For example, if two or more components are equal, it leads to a multi-valued function, while in many applications it is highly desirable to have a single-valued map (i.e., having a unique maximizer). To address this, the most common approach is to introduce a regularizer $h$ that acts as a penalty (or \emph{control cost}) to the maximization objective. This yields the \emph{regularized argmax function} or \emph{smooth best response map}
$$
\bv{x} \mapsto \argmax_{\bv{z} \in \simplex{n}} \left\{\bv{x}^\top \bv{z} - h(\bv{z}) \right\}.
$$
A popular choice for the regularizer is the negative entropy function restricted to the simplex~\cite{BG-LP:17}. This choice leads to the softmax, also known as \emph{logit map} in the context of game theory, that is,
\begin{equation}
\label{eq:best_response_map}
    \ds \argmax_{\bv{z} \in \simplex{n}} \left\{\bv{x}^\top \bv{z} + \entropy(\bv{z}) \right\} = \softmax{(\bv{x})}.
\end{equation}

\subsection{The Proximal Gradient Flow}
\label{sec:general_softmaxflow_properties}
To solve the optimization problem in Eq.~\eqref{eq:entropy_regularized_problem_simplex} we propose the following dynamics:
\beq
\label{eq:softmax_gf}
\tau \bv{\dot x} = - \bv{x} + \softmax\left(- \eps^{-1} \nabla \F(\bv{x})\right) = \fs(\bv{x}).
\eeq
Before characterizing the links between  the optimization problem in Eq.~\eqref{eq:entropy_regularized_problem} and the dynamics in Eq.~\eqref{eq:softmax_gf} we make the following observations:
\begin{enumerate}
\item as the parameter $\eps$ increases, the dynamics in Eq.~\eqref{eq:softmax_gf} promotes solutions with higher entropy. This follows from the fact that as $\eps \to + \infty$ the softmax converges to the uniform distribution.
\item as $\eps$ goes to zero, the dynamics in Eq.~\eqref{eq:softmax_gf} converges to a trajectory that follows the steepest descent direction along the simplex vertices. Specifically, as $\eps \to 0^+$, the softmax function converges to the vector $\bv{e}_{i^\star(\bv{x})}$, where $\bv{e}_i \in \R^n$ denotes the $i$-th canonical basis vector and $i^\star(\bv{x})\in\argmin_{1\le j\le n}\nabla \F(\bv{x})_j$. Thus, the dynamics becomes $\bv{\dot x} \approx - \bv{x} + \bv{e}_{i^\star(\bv{x})}$, which means that the trajectory moves toward the vertex of the simplex corresponding to the coordinate of steepest descent of the function $\F$ on the simplex. This type of dynamics is reminiscent of the Frank–Wolfe algorithm~\cite{MF:PW-56}, where at each step the update direction is toward an extreme point (vertex) of the feasible set that minimizes the linear approximation of the objective function. 
\end{enumerate}

The next result formalizes that the dynamics in Eq.~\eqref{eq:softmax_gf}: (i) is the proximal dynamics associated to problem~\eqref{eq:entropy_regularized_problem}; (ii) encodes the optimal solution of problem~\eqref{eq:entropy_regularized_problem} as an equilibrium. In the next section we characterize convergence to the optimal solution/equilibrium. Remarkably, the convergence analysis also shows that Eq.~\eqref{eq:softmax_gf} is an energy model.

\begin{Lemma}[Relating the optimization in Eq.~{\eqref{eq:entropy_regularized_problem}} and the dynamics in Eq.~{\eqref{eq:softmax_gf}}]\label{lem:softmax_optimization_connection}
Let $\map{\F}{\R^n}{\R}$ be a convex and continuously differentiable function, and let $\eps > 0$. Consider the entropy-regularized optimization problem in Eq.~\eqref{eq:entropy_regularized_problem} and the dynamics in Eq.~\eqref{eq:softmax_gf}. 
Then:
\begin{enumerate} 
\item
\label{item1:op_soft}
the dynamics in Eq.~\eqref{eq:softmax_gf} is the proximal gradient flow associated with the optimization problem in Eq.~\eqref{eq:entropy_regularized_problem};
\item
\label{item2:op_soft}
a vector $\bv{\xstar} \in \R^n$ is an equilibrium point of Eq.~\eqref{eq:softmax_gf} if and only if $\bv{\xstar}$ is an optimal solution of the problem in Eq.~\eqref{eq:entropy_regularized_problem}.
\end{enumerate}
\end{Lemma}
The proof of Lemma~\ref{lem:softmax_optimization_connection} is given in Section~\ref{sec:proofs}.

\subsection{Properties and Convergence of the Softmax Gradient Flow}
We now show that the proximal gradient flow in Eq.~\eqref{eq:softmax_gf} converges exponentially to its equilibrium, which coincides with the optimal solution of the problem in Eq.~\eqref{eq:entropy_regularized_problem}, and admits an explicit convergence rate. Moreover, our convergence analysis also reveals that: (i) the softmax gradient flow Eq.~\eqref{eq:softmax_gf} is an energy model; (ii) the energy function is the cost  of the problem in Eq.~\eqref{eq:entropy_regularized_problem_simplex}. More precisely, the results of this section establish that 
\begin{equation}\label{eqn:energy}
V(\bv{x}) := \F(\bv{x}) + \eps \entropybar(\bv{x})
\end{equation}
is an energy function for Eq.~\eqref{eq:softmax_gf}.  To develop the convergence analysis, we start by showing that any trajectory of Eq.~\eqref{eq:softmax_gf} with initial condition in the simplex never leaves this set, so that the probability simplex is forward invariant for the dynamics. Thus, trajectories starting from a feasible point for problem~\eqref{eq:entropy_regularized_problem} remain feasible for all time.
\begin{Lemma}[Invariance property]
\label{lem:positive_softmax_gf}
The probability simplex $\simplex{n}$ is a forward invariant set for the softmax gradient flow in Eq.~\eqref{eq:softmax_gf}.
\end{Lemma}
The formal proof is in
Section~\ref{sec:proofs}. 
To study convergence, we characterize the evolution of $V$ given in Eq.~\eqref{eqn:energy} along the trajectories of the softmax gradient flow in Eq.~\eqref{eq:softmax_gf}. 
\begin{Lemma}[Variation of the cost along the trajectories of the softmax gradient flow]
\label{lem:equality_Dkl}
Let $\map{\F}{\R^n}{\R}$ be continuously differentiable on $\simplex{n}$. Given $\eps>0$, define $S(\bv{x}; \eps) := \softmax(-\eps^{-1}\nabla \F(\bv{x}))$. Consider a trajectory $\bv{x}(t)$ of the softmax gradient flow in Eq.~\eqref{eq:softmax_gf} starting from the simplex. Then, for every $t>0$, it holds
\begin{equation}
\label{eq:V_dissipation}
\tau \frac{d}{dt}V(\bv{x}) = -\eps \Bigl( \DKL{\bv{x}}{S(\bv{x}; \eps)} + \DKL{S(\bv{x}; \eps)}{\bv{x}} \Bigr).
\end{equation}
\end{Lemma}

The proof is given in Section~\ref{sec:proofs}. 
Eq.~\eqref{eq:V_dissipation} in Lemma~\ref{lem:equality_Dkl} provides an exact dissipation identity for the entropy-regularized objective $V$ along the trajectories of the softmax gradient flow.
Since the symmetric Kullback--Leibler divergence is nonnegative and vanishes only when $\bv{x}=S(\bv{x};\eps)$, the lemma shows that $V$ is strictly decreasing along every nonstationary trajectory (i.e., solution of Eq.~\eqref{eq:softmax_gf}). Hence, $V$ is an energy (Lyapunov-like) function for the dynamics, and the dissipation identity provides the key estimate to establish exponential convergence. We can now prove exponential convergence of the softmax gradient flow.
\begin{Theorem}[Exponential convergence of the softmax gradient flow in the simplex]
\label{thm:exp_convergence}
Let $\map{\F}{\R^n}{\R}$ be convex and continuously differentiable on $\simplex{n}$ and consider the softmax gradient flow in Eq.~\eqref{eq:softmax_gf}. Let $\bv{x}^\star\in\simplex{n}$ be the unique minimizer of problem~\eqref{eq:entropy_regularized_problem} and let $\bv{x}(t)$ be a solution of Eq.~\eqref{eq:softmax_gf} rooted from some initial condition $\bv{x}(0)$.
Then, for every $\bv{x}(0)\in\simplex{n}$, the following holds: 
\begin{enumerate}
\item
\label{item1:proof_conv} the energy $V(\bv{x}(t))$ exponentially converges to the energy minimum $V(\bv{x}^\star)$ with rate $c = {\tau}^{-1}$, that is
\begin{equation}
    \label{eq:exp_cost_bound}
    V(\bv{x}(t))-V(\bv{x}^\star) \le \e^{-t/\tau}\bigl(V(\bv{x}(0))-V(\bv{x}^\star)\bigr),\qquad\text{for all }t\ge 0,
\end{equation}
\item 
\label{item2:proof_conv}
the solution $\bv{x}(t)$ of the dynamics in Eq.~\eqref{eq:softmax_gf} exponentially converges to $\bv{x}^\star$ with rate $c = \frac{1}{2\tau}$, that is
\begin{equation*}
\|\bv{x}(t)-\bv{x}^\star\|_2 \le \sqrt{\frac{2}{\eps}\bigl(V(\bv{x}(0))-V(\bv{x}^\star)\bigr)}\,\e^{-t/(2\tau)},
\qquad\text{for all }t\ge 0.
\end{equation*}
\end{enumerate}
\end{Theorem}

The proof is in Section~\ref{sec:proofs}.
Theorem~\ref{thm:exp_convergence} shows that the softmax gradient flow converges to the solution of the entropy-regularized optimization problem~\eqref{eq:entropy_regularized_problem}. Indeed, trajectories starting in the simplex remain feasible for all time by Lemma~\ref{lem:positive_softmax_gf}, and by Lemma~\ref{lem:softmax_optimization_connection} the equilibrium points of Eq.~\eqref{eq:softmax_gf} coincide with the minimizers of the entropy-regularized objective. Since this minimizer is unique, every trajectory converges to the same point $\bv{x}^{\star}$. In particular, the equilibrium $\bv{x}^{\star}$ is globally exponentially stable on $\simplex{n}$.

\begin{Remark}[A Bregman--PL-type inequality]
The exponential convergence result from Theorem~\ref{thm:exp_convergence} can be interpreted as a
Polyak--\L ojasiewicz (PL)-type inequality in the geometry induced by
the negative entropy.
Specifically, the optimality gap satisfies $V(\bv{x}) - V(\bv{x}^\star) \le \eps\, \DKL{\bv{x}}{S(\bv{x};\eps)}$,
where $S(\bv{x};\eps) = \softmax(-\eps^{-1}\nabla \F(\bv{x}))$.
This inequality plays the role of a PL inequality, where the Euclidean gradient residual is replaced by the Bregman divergence associated with the negative entropy, namely the Kullback--Leibler divergence~\cite{LMB:67,AB-MT:03}.
Since $\bv{x} = S(\bv{x};\eps)$ if and only if $\bv{x}$ is optimal, the quantity $\DKL{\bv{x}}{S(\bv{x};\eps)}$ acts as a measure of stationarity in the simplex.
\end{Remark}

\subsection{Introducing Biases}\label{sec:biases}
The derivations in this section naturally extend to optimization problems of the form:
$$
\ds \min_{\bv{x} \in \simplex{n}} \F(\bv{x}) + \eps\DKL{\bv{x}}{\hat{\bv{x}}},
$$
where $\hat{\bv{x}}\in \simplex{n}$ is a given vector. This formulation introduces a bias in the optimal solution. Essentially, the entropy is replaced by the \KL~divergence, which biases the optimal solution toward the vector $\hat{\bv{x}}$. The formal treatment presented in this section extends to this setting. In fact, one can rewrite the \KL~divergence regularized problem  as $\ds \min_{\bv{x} \in \simplex{n}} (\bar{\F}(\bv{x};\eps) - \eps\entropy{(\bv{x})})$, where $\bar{\F}(\bv{x};\eps):=\F(\bv{x})+ \eps\entropy{(\bv{x},\hat{\bv{x}})}$ is again convex and continuously differentiable in $\bv{x}$.

\section{Model Details}\label{sec:model_details}
We first present the derivations for the normative framework, the softmax gradient flow, and its neural instantiation from the main text. 
Then, we show why, as remarked in the main text, our model yields a prior-weighted softmax gating when a weight bias is introduced in the formulation.
We recall from the main text that $\bv{X}_k$ and $\bv{U}_k$ are the random variables of the {state} and {action} at time step $k$, respectively. The state space is {$\sX\subseteq\R^{n_x}$ and the action space is $\sU\subseteq\R^{n_u}$}. The time-indexing is chosen so that the environment transitions from state $\states{k-1}$ to $\states{k}$ when action $\bv{u}_k$ is applied. The agent selects its action by sampling from the policy $\policy{k}{k-1}$.
{The environment model available to the agent is denoted by $\plant{k}{k-1}$; from the main text, we also recall that $\jointxu{k}{k-1} = \plant{k}{k-1}\policy{k}{k-1}$.}

\subsection{Normative Framework} For convenience, we report here the optimization problem from the main text:
\begin{equation}\label{eqn:gateframe_SI}
\begin{split}
    \min_{\weights{k}\in \simplex{\np}} & \DKL{\jointxu{k}{k-1}}{\refjointxu{k}{k-1}} - \eps \entropy{(\weights{k})}\\
    \text{s.t. }& \policy{k}{k-1} = \sum_{\alpha=1}^{\np}\weight{k}{\alpha} \primitive{k}{k-1}{\alpha}.
\end{split}
\end{equation}
All symbols are defined in the main text and we start with formalizing the  mild assumptions described in the results.
\begin{Assumption}
\label{asn:bounded_cost} 
The set $\sW_k := \setdef{\weights{k}}{\weights{k} \mbox{ is feasible for problem~\eqref{eqn:gateframe_SI} and the cost is finite}}$ is non-empty.
\end{Assumption}

\begin{Assumption}
\label{assum_primitives}
The primitives are uniformly bounded and have full support over the action space $\mathcal{U}$.
Specifically, there exist two finite constants $0 < \subscr{\pi}{min} \leq \subscr{\pi}{max}$ such that, for all $\alpha \in \until{\np}$ and for all time steps $k$, $\subscr{\pi}{min} \leq \primitive{k}{k-1}{\alpha}\leq\subscr{\pi}{max}$.
\end{Assumption}
Intuitively, Assumption~\ref{asn:bounded_cost} makes the optimization meaningful, in the sense that the optimal value of problem~\eqref{eqn:gateframe_SI} is finite. This assumption is standard in the context of, e.g., entropy-regularized optimization. For example, the assumption is 
satisfied when $\plant{k}{k-1}$ is absolutely continuous with respect to $\refplant{k}{k-1}$, or when the primitives are absolutely continuous with respect to $\refpolicy{k}{k-1}$. Assumption~\ref{assum_primitives} formalizes the fact that primitives cover the action space and have bounded moments (this is always satisfied for, e.g., finite discrete variables or truncated Gaussians).

The next result formalizes the unconstrained reformulation of the optimization problem stated under our normative framework. The reformulation from the main text is also reported here for convenience
\begin{equation}
\min_{\weights{k}\in \R^{\np}} \F(\weights{k}) + \eps \entropybar {(\weights{k})},
\label{eq:sub2_stepk_modified_DKL}
\end{equation}
where 
\begin{equation}
\label{eqn:F_SI}
\begin{aligned}
\F(\weights{k}) :=  \sum_{\alpha =1}^{\np}\weight{k}{\alpha} \bigg(\E_{\primitive{k}{k-1}{\alpha}}\bigg[ &\ln\sum_{\beta =1}^{\np}\weight{k}{\beta}\primitive{k}{k-1}{\beta} - \ln\refpolicy{k}{k-1} \\
&\qquad+ \DKL{\plant{k}{k-1}}{\refplant{k}{k-1}}\bigg] \bigg).
\end{aligned}
\end{equation}
The result also derives the gradient of $\F$, characterizing its convexity and smoothness properties.
\begin{Lemma}[Normative framework reformulation]
\label{lem::gradient}
Let Assumptions~\ref{asn:bounded_cost} and~\ref{assum_primitives} hold. Then,
\begin{enumerate}
\item \label{item1:equivalence_ops_DKL}
the problem defined in Eqs.~\eqref{eq:sub2_stepk_modified_DKL} and~\eqref{eqn:F_SI} is a reformulation of the optimization problem in Eq.~\eqref{eqn:gateframe_SI};
\item \label{item2:gradient_op_reduced_DKL}
the gradient of $\F(\weights{k})$ is:
\begin{equation*}
\begin{aligned}
\nabla \F(\weights{k}) = \E_{\primitives{k}{k-1}}\bigl[ &
\ln \weights{k}^\top\primitives{k}{k-1} + \DKL{\plant{k}{k-1}}{\refplant{k}{k-1}}\\
&- \ln\refpolicy{k}{k-1}\bigr] + \1_{\np},
\end{aligned}
\end{equation*}
where we are using the notation $\E_{\primitives{k}{k-1}} = \left(\E_{\primitive{k}{k-1}{\alpha}}\right)_{\alpha \in \until{\np}} \in \R^{\np}$;
\item \label{item3:smooth_constant}
the function $\F$ is convex and 
$\subscr{L}{$\F$}$-smooth in $\weights{k}$ with constant $\subscr{L}{$\F$} = \np \frac{\subscr{\pi}{max}}{\subscr{\pi}{min}}$.
\end{enumerate}
\end{Lemma}
The proof is in Section~\ref{sec:proofs}.
Item~\ref{item1:equivalence_ops_DKL} shows the reformulation in the Results section of the main text and derives the functional expression for $\F$ given in the Materials and Methods section of the main text. Item~\ref{item2:gradient_op_reduced_DKL} computes the gradient of $\F$ and this expression is used in both \results~and \methods~sections of the main text to obtain the softmax flow and its neural instantiation. Item~\ref{item3:smooth_constant} characterizes convexity and smoothness of the problem. 

\subsection{Softmax Gradient Flow}
Building on the formal treatment developed so far, the dynamics of the softmax gradient flow associated to our normative framework optimization is:
\begin{align*}
\tau \dot{\bv{w}}_k {=} - \weights{k} + \softmax\big({-} \eps^{-1} &\E_{\primitives{k}{k-1}}\bigl[\ln \weights{k}^\top\primitives{k}{k-1} \\
&\; + \DKL{\plant{k}{k-1}}{\refplant{k}{k-1}} - \ln\refpolicy{k}{k-1}\bigr] \big),
\end{align*}
or, component-wise:
\begin{align*}
\tau \dot{w}_k^\alpha {=} - \weight{k}{\alpha} + \softmax\big({-} \eps^{-1} &\E_{\primitives{k}{k-1}}\bigl[
\ln \weights{k}^\top\primitives{k}{k-1} \\
& + \DKL{\plant{k}{k-1}}{\refplant{k}{k-1}} - \ln\refpolicy{k}{k-1}\bigr] \big)_{\alpha},
\end{align*}
with $\alpha\in \until{\np}$. This is the dynamics from the main text, where we unpack the exponent of the softmax and we leverage the invariance property (Lemma~\ref{rem:softmax_sum}) to drop the constant term $-\1_{\np}$ in the gradient.

In the main text we explicitly discuss the setting where the action space is discrete so that there are  $\du$ possible actions. In this case, it is possible to define the matrix $\Pi(\states{k-1}) \in \R^{\du \times \np}$
\begin{equation*}
    \Pi(\states{k-1}):=
    \begin{pmatrix}
        \primitive{k,1}{k-1}{1} &\primitive{k,1}{k-1}{2} &\hdots &\primitive{k,1}{k-1}{\np}\\
        \primitive{k,2}{k-1}{1} &\primitive{k,2}{k-1}{2} &\hdots &\primitive{k,2}{k-1}{\np}\\
        \vdots &\vdots &\ddots &\vdots\\
        \primitive{k,\du}{k-1}{1}& \primitive{k,\du}{k-1}{2} &\hdots &\primitive{k,\du}{k-1}{\np},
    \end{pmatrix}.
\end{equation*}
Consequently, the dynamics of the softmax flow can be written as
\begin{equation}\label{eqn:modelflow_discrete}
\tau \dot{\bv{w}}_k = - \weights{k} + \softmax\bigl(- \eps^{-1} \Pi(\bv{x}_{k-1})^\top\bigl(\ln \bigl({\Pi(\bv{x}_{k-1}) \bv{\weights{k}}}\bigr) + \bv{s}(\bv{x}_{k-1} , \bv{u}_{k})\bigr) \bigr),
\end{equation}
where the input $\bv{s}$ is defined as in the main text, i.e., $\bv{s}(\bv{x}_{k-1} , \bv{u}_{k}) := \DKL{\plant{k}{k-1}}{\refplant{k}{k-1}} - \ln\refpolicy{k}{k-1}$. All the desirable properties of the softmax flow from the main text follow from the formal treatment in Section~\ref{sec:general_softmaxflow_properties}.
In particular:
\begin{enumerate}
\item the probability simplex $\simplex{n}$ is a forward invariant set for the softmax flow. This means that, if the dynamics is initialized with initial conditions in the probability simplex, then trajectories of the softmax flow remain in the simplex;
\item an equilibrium exists for the softmax flow. Moreover, a point in the state space is an equilibrium for the softmax flow if and only if it is also an optimal solution of our normative framework optimization;
\item the softmax flow is an energy model and the energy function is the cost in Eq.~\eqref{eq:sub2_stepk_modified_DKL}. Moreover, the energy function decreases monotonically along trajectories and vanishes only at the equilibrium (the optimal solution of Eq.~\eqref{eqn:gateframe_SI});
\item finally, for any feasible initial condition, the trajectories of the softmax flow converge exponentially to the optimal solution with a guaranteed convergence rate of $\frac{1}{2\tau}$.
\end{enumerate}

\subsection{Neural Circuit}
\label{sec:bio_implementation_SI}
We derive separately the fast and slow units (Fig.~2(c) in the main text) and then integrate them together.
We recall that the matrix of the primitives is conditioned on the current state available to the agent. This aspect is highlighted by explicitly writing $\Pi(\states{k-1})$. Here, we denote by $\Pi_i(\states{k-1})$ the $i$-th row of the matrix $\Pi(\states{k-1})$, for all $i \in \until{\du}$,
and by $\Pi^{\alpha}(\states{k-1})$ the $\alpha$-th column of $\Pi(\states{k-1})$, for all $\alpha \in \until{\np}$.

\noindent{\bf Fast unit.}
As detailed in the main text, given the current state of $\weights{k}$, the fast unit computes the vector $\bv{y}$, whose value at steady state is $ \bar{\bv{y}}=- \eps^{-1}\nabla \F(\weights{k})$.
Component-wise:
\begin{align*}
\bar{y}^{\alpha} = - \eps^{-1}\left(\nabla \F(\weights{k})\right)_{\alpha}&= - \eps^{-1}\sum_{i = 1}^{\du}\left(\ln \Pi_i(\states{k-1})\weights{k} + s(\bv{x}_{k-1} , \bv{u}_{k,i})\right)\primitive{k,i}{k-1}{\alpha} \\
& = - \eps^{-1}\sum_{i = 1}^{\du}\left(\ln \bar{a}_i + s(\bv{x}_{k-1}, \bv{u}_{k,i})\right)\primitive{k,i}{k-1}{\alpha}\\
& = - \eps^{-1}\sum_{i = 1}^{\du}\left(\bar{b}_i + s(\bv{x}_{k-1} , \bv{u}_{k,i})\right)\primitive{k,i}{k-1}{\alpha}\\
&\ds = - \eps^{-1}  (\Pi^{\alpha}(\states{k-1}))^\top (\bar{\bv{b}} +  \bv{s} ),
\end{align*}
where, as in the main text, $\bar{a}_i:= \Pi_i(\states{k-1})\weights{k}$, $\bar{b}_i:=\ln(\bar{a}_i)$, $\bar{\bv{a}} := [\bar{a}_1, \ldots,\bar{a}_{\du}]^\top$, $\bar{\bv{b}} := [\bar{b}_1, \ldots,\bar{b}_{\du}]^\top$, and $\bv{s} :=[s_1,\ldots,s_{\du}]^\top$. Let $\bv{a} := [a_1, \ldots,a_{\du}]^\top$, and $\bv{b} :=[b_1,\ldots,b_{\du}]^\top$; then, the vector $\bv{y}$ can be computed via a dynamical system with: (i) two sets of equations, with state variables $\bv{a}$ and $\bv{b}$, that converge to $\bar{\bv{a}}$ and $\bar{\bv{b}}$, respectively; (ii) a slower set of dynamical equations, with state variables $\bv{y}$, that converge to the quantity in the last line of the above derivations. 
This yields:
\begin{align*}
\subscr{\tau}{g} \dot{\bv{a}} &=-\bv{a}+\Pi(\states{k-1}) \bv{w}_k,\\
\subscr{\tau}{g} \dot{\bv{b}} &=-\bv{b}+\ln{(\bv{a})},\\
\tilde{\tau}_g\dot{\bv{y}} &= - \eps\bv{y} -\Pi(\states{k-1})^\top \left(\bv{b} + \bv{s}\right),
\end{align*}
 with $\subscr{\tau}{g} \ll \tilde{\tau}_g$. This is the dynamics from the main text corresponding to the circuit for the fast unit in Fig.~2(c).\\

\noindent{\bf Slow unit.}
The core of this unit in Fig.~2(c) is the biologically plausible softmax neural circuit from the literature~\cite{MS-JO:22}. Complementing the derivations in \methods section of the main text, the unit receives as input the components of $\bv{y}$, $y^\alpha$, and returns the vector $\bv{w}_k$, whose value at steady state is $\bv{\bar{w}}_k$ having as components $\bar{w}_k^\alpha = \softmax(\bv{y})_\alpha$. In the circuit from the main text, the neuron $m$ at steady state computes the softmax normalizer $\bar{m}$. Each of the $\bv{r}$-neurons converges to $y^\alpha - \ln\left(\bar{m}\right)$, so that the $\bv{w}_k$-neurons effectively return the softmax, providing the optimal weights. This yields the dynamics in \results  section of the main text, also reported here, component-wise:
\begin{align*}
\subscr{\tau}{s} \dot m &= - m + \sum_{\alpha=1}^{\np} \exp(y^\alpha),\\
\subscr{\tau}{s} \dot r^\alpha &= - r^\alpha + y^\alpha - \ln(m),  \\
\tau \dot w^\alpha_k &= - w^\alpha_k + \exp(r^\alpha),
\end{align*}
where $\alpha=1,\ldots,\np$ and $\subscr{\tau}{s} \ll \tau$ to ensure that the internal dynamics (of neuron $m$ and $r^\alpha$-neurons) is faster than the output dynamics returning the weights.\\

\noindent{\bf The complete dynamics.} The dynamics in the main text, corresponding to the circuit in Fig.~2(c), is obtained by combining the dynamics for the fast and slow units, with $\tilde{\tau}_g \ll \tau_s$ so that the former dynamics is effectively faster than the latter.

\subsection{Introducing Biases Yields Prior-Weighted Softmax Gating} Following the arguments at the end of Section~\ref{sec:general_results}, 
biases can be introduced in the normative framework optimization by replacing the entropic regularizer with the \KL~divergence. In this case, the optimization becomes
\begin{equation}
\label{eqn:gateframe_biased}
\begin{aligned}
\min_{\weights{k}\in \simplex{\np}} & \DKL{\jointxu{k}{k-1}}{\refjointxu{k}{k-1}} + \eps \DKL{\weights{k}}{\pweights{k}}\\
\text{s.t. }& \policy{k}{k-1} = \sum_{\alpha=1}^{\np}\weight{k}{\alpha} \primitive{k}{k-1}{\alpha},
\end{aligned}
\end{equation}
where $\pweights{k}\in\operatorname{ri}(\simplex{\np})$ is a bias vector.
By Proposition~\ref{prop:cross_entropy}, the identity $\ds \DKL{\weights{k}}{\pweights{k}} = \entropy{(\weights{k}, \pweights{k})} - \entropy{(\weights{k})}$ holds, implying that the above optimization problem remains convex. Formally, the same derivations presented above can be followed in this more general setting. By doing so, when the bias $\pweights{k}$ is introduced, the softmax flow becomes
\begin{equation}\label{eqn:modelflow_discrete_Gumbel}
\tau \dot{\bv{w}}_k = - \weights{k} + \softmax\Bigl( - \eps^{-1} \bv{h}(\weights{k}) + \ln\pweights{k}\Bigr),
\end{equation}
where we have defined the vector $\bv{h}(\weights{k}) =\Pi(\bv{x}_{k-1})^\top\left(\ln \left({\Pi(\bv{x}_{k-1}) \weights{k}}\right) + \bv{s}(\bv{x}_{k-1} , \bv{u}_{k})\right)$.
Essentially, the above equation differs from Eq.~\eqref{eqn:modelflow_discrete} in the presence of the logarithmic term in the softmax that depends on the bias vector. The equilibrium of Eq.~\eqref{eqn:modelflow_discrete_Gumbel}, that is also the optimal solution of Eq.~\eqref{eqn:gateframe_biased}, is
\begin{equation}\label{eqn:gumbel_softmax}
\optimalweights{k} = \softmax\Bigl(- \eps^{-1} \bv{h}(\optimalweights{k}) + \ln\pweights{k}\Bigr).
\end{equation}
Component-wise,
\begin{equation*}
w_{k}^{\alpha,\star}
=
\frac{
\pweight{k}{\alpha}
\exp\bigl(-\eps^{-1}h_\alpha(\optimalweights{k})\bigr)
}{
\displaystyle
\sum_{\beta=1}^{\np}
\pweight{k}{\beta}
\exp\bigl(-\eps^{-1}h_\beta(\optimalweights{k})\bigr)
}.
\end{equation*}
This is a probability mass function that twists the probability vector $\pweights{k}$ with an exponential~\cite{PG-MR-RW:14} kernel appearing in policies that minimize the variational free energy, including KL control. See, e.g., the surveys in~\cite{AS-HJ-KF-GR:25_new,EG-GR:22}.  Moreover, Eq.~\eqref{eqn:modelflow_discrete_Gumbel} also reveals that the neural circuit in the main text does not necessarily need to be modified to account for biases. These can be included in the input to the network.

\section{Supplemental Implementation Details and Experiments}\label{sec:experiments_details}
We provide supplemental implementation details and experiments. In the experiments reported in the main text, the optimal weights are computed via the softmax gradient flow, after verifying on a subset of experiments that its neural implementation returns the same weights. Further implementation details are available at the accompanying repository~\cite{FR:26_repo}.

\subsection{Emerging Flocking Behaviors Experiments}
We refer to Table~\ref{tab:flocking_parameters} for descriptions and values of the parameters described below. We recall from the main text that the state of boid (agent) $i$, $\mathbf{x}_k^i$, is the stack of the boid's position $\mathbf{p}_k^i=\big[p_{x,k}^i,p_{y,k}^i\big]^\top$ and velocity $\mathbf{v}_k^i=\big[v_{x,k}^i,v_{y,k}^i\big]^\top$. The maximum velocity is bounded; more precisely, there exists $v_{\max} > 0$ such that $\|\mathbf{v}_k^i\| \leq v_{\max}$. The action of boid $i$, $\mathbf{u}_k^i=\big[u_{x,k}^i,u_{y,k}^i\big]^\top$, is the boid acceleration. The acceleration is also bounded, that is, there exists $u_{\max} > 0$ such that $\|\mathbf{u}_k^i\| \leq u_{\max}$. 
The environment model for boid $i$, $\plantboid{k}{k-1}{i}$, is a Gaussian with covariance specified in Table~\ref{tab:flocking_parameters} and mean given by the second-order dynamics from~\cite{FC-SS:07, HL-WJR-IC:00,RO:06}:
\begin{equation}
\label{flock_model}
\begin{cases}
\mathbf{p}_k^i=\mathbf{p}_{k-1}^i+\mathbf{v}_{k-1}^i dt\\
\mathbf{v}_k^i=\mathbf{v}_{k-1}^i+\mathbf{u}_k^i dt
\end{cases},\quad i=1,\ldots,N.
\end{equation}
In the simulations from the main text, the dynamics of the boids is the one in Eq.~\eqref{flock_model}. In the simulations, $\bv{u}_k$ is sampled at each $k$ from the policy obtained with our computational model. In all the experiments, at a given $k$, a neighboring boid $j$ is visible to boid $i$ if the relative angle
$
\theta_{ij,k}
:=
\arccos\!\left(
\frac{
(\mathbf v_k^i)^\top(\mathbf p_k^j-\mathbf p_k^i)
}{
\|\mathbf v_k^i\|\,
\|\mathbf p_k^j-\mathbf p_k^i\|
}
\right)
\leq \frac{\alpha}{2},
$
where $\alpha$ is the vision angle (see also Fig.~3(a) in the main text). The three concentric zones are from the literature, see~\cite{IDC-JK-RJ-GDR-NRF:02} and related references in \results.
Following the definition of social forces in~\cite{RL-YL-LEK:10}, the means of the Gaussian primitives are given by:
\begin{itemize}
\item 
for the \emph{separation primitive}
$$
\mathbf{u}_{\mathrm{sep}}^i =
\begin{cases}
\ds \frac{1}{|\mathcal{S}^i_k|} \sum\limits_{j \in \mathcal{S}^i_k} g\left(\|\mathbf{p}_k^i - \mathbf{p}_k^j\|\right) \frac{\mathbf{p}_k^i - \mathbf{p}_k^j}{ \|\mathbf{p}_k^i - \mathbf{p}_k^j\|}, & \text{if } |\mathcal{S}^i_k| > 0, \\
\0_2, & \text{otherwise},
\end{cases}
$$
with $\mathcal{S}^i_k$ being the set of neighbors visible to boid $i$ in the separation zone;
\item 
for the \emph{alignment primitive}
$$
\mathbf{u}_{\mathrm{ali}}^i =
\begin{cases}
\ds\frac{1}{|\mathcal{A}^i_k|} \sum\limits_{j \in \mathcal{A}^i_k} g\left(\|\mathbf{p}_k^i - \mathbf{p}_k^j\|\right)\frac{\mathbf{v}_k^j}{\|\mathbf{v}_k^j\|}, & \text{if } |\mathcal{A}^i_k| > 0, \\
u_{\max}\frac{\mathbf{v}_k^i}{\|\mathbf{v}_k^i\|}, & \text{otherwise},
\end{cases}
$$
where $\mathcal{A}^i_k$ is the set of neighbors visible to boid $i$ in the alignment zone;
\item 
for the \emph{cohesion primitive} 
$$
\mathbf{u}_{\mathrm{coh}}^i =
\begin{cases}
\ds \frac{1}{|\mathcal{C}^i_k|} \sum\limits_{j \in \mathcal{C}^i_k} g\left(\|\mathbf{p}_k^i - \mathbf{p}_k^j\|\right) \frac{\mathbf{p}_k^j - \mathbf{p}_k^i}{\|\mathbf{p}_k^i - \mathbf{p}_k^j\|}, & \text{if } |\mathcal{C}^i_k| > 0, \\
\0_2, & \text{otherwise},
\end{cases}
$$
where $\mathcal{C}^i_k$ is the set of neighbors visible to boid $i$ in the cohesion zone.
\end{itemize}
The forces depend on the function $g$; this is a function of the distance between agents. According to, e.g.,~\cite{RL-YL-LEK:10} this is set to
\begin{equation*}
g(d) =
\begin{cases}
u_{\max} \left(1 - \dfrac{d}{r_{\mathrm{sep}}} \right) & \text{if } d \in [0, r_{\mathrm{sep}}] \\
u_{\max} & \text{if } d \in ]r_{\mathrm{sep}}, r_{\mathrm{ali}}]\\
u_{\max} \left( \dfrac{r_{\mathrm{coh}} - d}{r_{\mathrm{coh}} - r_{\mathrm{ali}}} \right) & \text{if } d \in ]r_{\mathrm{ali}}, r_{\mathrm{coh}}].
\end{cases}
\end{equation*}
The function $g$  modulates the magnitude of the social forces. The function is decreasing inside the separation zone. Then, it becomes constant in the alignment zone and finally decreases again in the cohesion zone. Note that, according to the above social forces, when there is no neighbor in the alignment zone, the boid is encouraged to maintain its current heading~\cite{HL-WJR-IC:00}.

The polarization order parameter in Fig.~\ref{fig:boids_metrics} is defined as in~\cite{IDC-JK-RJ-GDR-NRF:02,KT-YK-CCI-CH-MJL-IDC:13}:
\begin{equation*}
    p_k = \left\| \frac{1}{N} \sum_{i=1}^{N} \frac{\mathbf{v}^i_k}{\|\mathbf{v}^i_k\|} \right\|.
\end{equation*}
For the experiment shown in Fig.~3(b) of the main text, the generative model is inspired by~\cite{CH-BM-LDC-RPM-KJF-IDC:24}. Specifically, $\refplant{k}{k-1}$ is a Gaussian with covariance specified in Table~\ref{tab:flocking_parameters} and centered in
\begin{equation*}
\bar{\mathbf{x}}_{\mathrm{q},k}^i =
\begin{bmatrix}
\displaystyle \frac{1}{|\mathcal{C}_k^i|} \sum_{j \in \mathcal{C}_k^i} \mathbf{p}^j_k\\
\displaystyle \frac{1}{|\mathcal{A}_k^i|} \sum_{j \in \mathcal{A}_k^i} \mathbf{v}^j_k
\end{bmatrix}.
\end{equation*}
This generative model captures the tendency of a boid to stay close to its neighbors and align with their velocity in a free environment~\cite{CWR:87,AMR-GEM-AT-PY-NTO:22,HL-GEM-JW-KVDV-RTV-AT-NTO:19,DJTS:06}. The milling order parameter in Fig.~\ref{fig:boids_metrics}(b) is:
\begin{equation*}
    m_k = \left\| \frac{1}{N} \sum_{i=1}^{N} \left(\frac{\mathbf{p}^i_k - \mathbf{p}_{\mathrm{cm},k}}{\|\mathbf{p}^i_k - \mathbf{p}_{\mathrm{cm},k}\|} \times \frac{\mathbf{v}^i_k}{\|\mathbf{v}^i_k\|} \right) \right\| ,
\end{equation*}
where $\mathbf{p}_{\mathrm{cm},k} = \frac{1}{N} \sum_{i=1}^{N} \mathbf{p}^i_k$ is the position of the flock's center of mass at step $k$.
The experiments {involve $60$ boids} and are obtained by including in the generative model a collision avoidance term. This term induces an increase in the separation weight when boids are too close.
The emergence of milling following an increase in separation weight is consistent with previous zonal models of collective motion, e.g.,~\cite{IDC-JK-RJ-GDR-NRF:02}. In these models individuals follow a hierarchical rule where short-range repulsion has higher priority over alignment and attraction. In regimes where alignment is insufficient to produce global consensus, this interaction structure gives rise to milling (torus) states.

For the soft-control experiments of Fig.~3(d) in  the main text,  $10\%$ of the $40$ boids are goal-informed (leaders). In these experiments, the generative model for the non-informed boids (followers) is left unchanged while informed boids are equipped with a model encoding a goal-directed behavior. Letting $\bv{p}_{\mathrm{g}}$ be the goal position (set to $\mathbf{p}_{\mathrm{g}}=[-15,-15]^\top$ in the experiments) the generative model for the leaders is again a Gaussian with the same covariance matrix used for the generative model in the previous experiments, this time centered in
$$
\bar{\mathbf{x}}_{\mathrm{q},k}^i =
\begin{bmatrix}
\mathbf{p}_{\mathrm{g}} \\
\displaystyle \frac{\mathbf{p}_{\mathrm{g}}-\mathbf{p}_k^i}{\|\mathbf{p}_{\mathrm{g}}-\mathbf{p}_k^i\|}
\end{bmatrix}.
$$
The distance in Fig.~\ref{fig:boids_metrics}(c) is defined as the Euclidean distance between the center of mass of the flock and the goal position:
\begin{equation*}
d_k = \left\| \mathbf{p}_{\mathrm{cm},k} - \mathbf{p}_{\mathrm{g}} \right\|.
\end{equation*}
Finally, in Fig.~\ref{fig:boid_soft_control} the followers evolve according to three different models from the literature~\cite{IDC-JK-RJ-GDR-NRF:02,FC-SS:07,TV-AC-EBJ-IC-OS:95}. 
We refer to the original papers and to the repository~\cite{FR:26_repo} for the description of these models. 

\subsection{Exploration/Exploitation in Human Decision-Making Experiments} In the experiments from the main text we use data and numerical methods from~\cite{SJG:18}.
Beliefs over rewards are updated via a Kalman filter. At each trial $k=1,\ldots,T$, the participant $i$ selects arm $j=u^i_k\in \{1, 2\}$ based on the current belief state $\bv{x}^i_{k-1}$ and subsequently observes reward $r^i_k$. 
To streamline notation, from now on we omit the participant superscript.
The state is the vector of posterior means and variances updated as in~\cite{SJG:18}. Specifically, for the
chosen arm $j=u_k$, the posterior mean and variance are updated
according to
\begin{equation}
\label{eqn:gershman_state}
\begin{aligned}
\mu_{j,k}
&=
\mu_{j,k-1}
+
\alpha_{j,k-1}
\left(r_k-\mu_{j,k-1}\right),
\\
s_{j,k}
&=
s_{j,k-1}
-
\alpha_{j,k-1}s_{j,k-1},
\end{aligned}
\end{equation}
with learning rate
\[
\alpha_{j,k-1}
=
\frac{s_{j,k-1}}{s_{j,k-1}+\nu},
\]
where $\nu=10$ is the observation-noise variance used in the
Kalman-filter implementation.
The posterior parameters of the unchosen arm remain unchanged.
At the beginning of each of the $B$ blocks, the belief means and
variances are reset to the prior values reported in
Table~\ref{tab:mab_parameters}, and the latent reward mean of each stochastic arm $j$ is sampled from a zero-mean Gaussian distribution with variance equal to $s_{j,0}$ and held fixed throughout the block.
In Experiment~1, one arm is stochastic and the other yields a fixed reward of zero. In Experiment~2, both arms are stochastic. 
The Hybrid model in~\cite{SJG:18} is
\begin{equation*}
    p_{\mathrm{h}}(u_k = 1|\bv{x}_{k-1}) = \Phi\left(\gamma \, \left(\sqrt{s_{1,k-1}} - \sqrt{s_{2,k-1}}\right) + \beta \, \frac{\mu_{1,k-1} - \mu_{2,k-1}}{\sqrt{s_{1,k-1} + s_{2,k-1}}}\right),
\end{equation*}
where $\Phi(\cdot)$ denotes the standard normal cumulative distribution function, $\gamma$ is the coefficient for directed exploration, and $\beta$ is the coefficient for random exploration. In the main text, the coefficients $\gamma$ and $\beta$ are inferred from the data via probit regression.
The exploitation, uncertainty-seeking and risk-aversion primitives from the main text are, respectively:
\begin{equation*}
\begin{split}
    \pi^1(u_k = j|\bv{x}_{k-1}) &= 
    \begin{cases}
        0.99 \ & \text{if $j=\arg \max_{i \in \{1,2\}} {\mu_{i,k-1}} $}\\
        0.01 \ & \text{otherwise}
    \end{cases}
    ,\\
    \pi^2(u_k = j|\bv{x}_{k-1}) &= \ds\frac{s_{j,k-1}}{\sum_{i=1}^2 s_{i,k-1}}
    ,\\
    \pi^3(u_k = j|\bv{x}_{k-1}) &= \ds\frac{\exp({-s_{j,k-1}})}{\sum_{i=1}^2 \exp({-s_{i,k-1}})}
    .
\end{split}
\end{equation*}
For the exploitation primitive, if $\mu_{1,k-1} = \mu_{2,k-1}$, both arms are assigned equal probability, i.e., $\pi^1(u_k = j \mid \mathbf{x}_{k-1}) = 1/2$ for $j \in \{1,2\}$.
The parameter values are given in Table~\ref{tab:mab_parameters}.
Since we  do not want to learn additional parameters over~\cite{SJG:18}, leveraging the flexibility offered by our normative framework, we set $q(\bv{x}_{k} \mid \bv{x}_{k-1}, u_k)$ equal to $p(\bv{x}_{k} \mid \bv{x}_{k-1}, u_k)$. In fact, notice from the optimization properties highlighted in the main text that, by doing so, the complexity term $\F$ does not depend on $p(\bv{x}_{k} \mid \bv{x}_{k-1}, u_k)$ and $q(\bv{x}_{k} \mid \bv{x}_{k-1}, u_k)$.

In Fig.~\ref{fig::mab_logit}(a)--(b) we present supplemental evaluation of our model in the same settings considered in Fig.~4(b)--(c) from the main text, but using logistic (softmax) regression rather than the probit regression. Further, Fig.~\ref{fig::mab_logit}(c) shows the Predictive Log-Likelihood (PLL) from a 3-fold cross validation among models in the case of logit regression. These findings confirm the results from the main text. Finally, Fig.~\ref{fig::mab_logit}(d) shows the PXP comparison for $\eps \in [0.001, 1]$, using logit-regression approach. The figure indicates that low temperatures can be beneficial, whereas excessively large values of $\eps$ degrade performance. This deterioration, confirming the result in Fig.~\ref{fig::mab_additional_probit_results}(b), is, however, expected. In fact, as temperature increases, the entropic term in our model dominates, driving the weights toward a uniform distribution.

\subsection{Layered Control Experiments}
Parameter values and descriptions are in Table~\ref{tab:quadrotor_parameters}, unless differently specified.  We recall that the state of the quadrotor is
$\mathbf{x}_k =
\begin{bmatrix}
p_{x,k}, \dot{p}_{x,k}, p_{z,k}, \dot{p}_{z,k}, \theta_k, \dot{\theta}_k
\end{bmatrix}^\top$,
where $p_{x,k}$ and $p_{z,k}$ denote the horizontal and vertical position, $\dot{p}_{x,k}$ and $\dot{p}_{z,k}$ the velocities on the x and z axes, respectively, $\theta$ is the pitch angle and $\dot{\theta}$ is the angular velocity (see Fig.~5(a) in the main text). The control input is $\mathbf{u}_{k} = \left[u_{1,k}, u_{2,k} \right]^\top$, representing the thrust commands applied to the two rotors.
The reference trajectory is defined by the tuple ${(\bv{x}_{\mathrm{ref},k}, \bv{u}_{\mathrm{ref},k})}_{k=1,\ldots,T}$, where $\bv{x}_{\mathrm{ref},k}=\begin{bmatrix}
p_{x,\mathrm{ref},k}, \dot{p}_{x,\mathrm{ref},k}, p_{z,\mathrm{ref},k}, \dot{p}_{z,\mathrm{ref},k}, \theta_{\mathrm{ref},k}, \dot{\theta}_{\mathrm{ref},k} 
\end{bmatrix}^\top$, $\bv{u}_{\mathrm{ref},k}=\left[u_{1,\mathrm{ref},k}, u_{2,\mathrm{ref},k} \right]^\top$. In our experiments, $\bv{x}_{\mathrm{ref},k}$ defines a lemniscate trajectory (see Fig.~5(b) in the main text), with constant zero pitch angle and angular velocity.
The equations of motion of the quadrotor in the main text are
$$
\begin{cases}
\ddot{p}_x = \frac{\sin\theta}{m}(u_1 + u_2) + d_x\\
\ddot{p}_z = \frac{\cos\theta}{m}(u_1 + u_2) - g + d_z\\
\ddot{\theta} = \frac{l}{\sqrt{2} I_{yy}} (u_2 - u_1)
\end{cases},
$$
where $m$ is the mass of the quadrotor, $g$ is the gravitational acceleration, $l$ is the arm length of the quadrotor, i.e., the distance from each motor pair to the center of mass, $I_{yy}$ is the moment of inertia about the $y$-axis, and $d_x$ and $d_z$ are, e.g., environmental disturbances. Both $d_x$ and $d_z$ are drawn from a Gaussian distribution $\mathcal{N}(0,\sigma_{\mathrm{d}}^2)$. 
In the experiments, the cost $c$, also used to compute the total episode costs, is quadratic:
\begin{equation*}
c(\bv{x}_k,\bv{u}_k) = (\bv{x}_k-\bv{x}_{\mathrm{ref},k})^\top Q (\bv{x}_k-\bv{x}_{\mathrm{ref},k}) + (\bv{u}_k-\bv{u}_{\mathrm{ref},k})^\top R (\bv{u}_k-\bv{u}_{\mathrm{ref},k}),
\end{equation*}
where $Q$ and $R$ are the state and action cost matrices, respectively. In accordance with the formulation of the objective from the main text, $\refjointxu{k}{k-1}=\frac{1}{Z} \refjointxutwisted{k}{k-1} \exp\left(-c(\bv{x}_k,\bv{u}_k)\right)$. Here, $Z$ is a normalization constant, the cost $c(\bv{x}_k,\bv{u}_k)$ is computed using the state predicted by the linearized drone dynamics from \texttt{Safe-Control-Gym}~\cite{ZY-AWH-SZ-LB-MG-JP-APS:22}, and $\refjointxutwisted{k}{k-1}=\tilderefplant{k}{k-1}\tilderefpolicy{k}{k-1}$, with $\tilderefplant{k}{k-1}=\plant{k}{k-1}$ and uniform $\tilderefpolicy{k}{k-1}$. 
See the repository~\cite{FR:26_repo} for details.
The two primitives from the main text are Gaussian distributions centered in the output of two Proportional Derivative (PD) controllers: 
\begin{equation}
\begin{split}
\primitive{k}{k-1}{1} &=
\mathcal{N}\left(\bv{u}^1_{k},
\Sigma_\mathrm{u} \right),\\
\primitive{k}{k-1}{2} &=
\mathcal{N}\left(\bv{u}^2_{k},
\Sigma_\mathrm{u} \right).
\end{split}
\label{primitives_quadrotor}
\end{equation}
The covariance matrix $\Sigma_\mathrm{u}$ captures actuation noise. The PD controllers/primitives used in the experiments are obtained from the environment~\cite{ZY-AWH-SZ-LB-MG-JP-APS:22} and are standard in the literature~\cite{JP-HZ-SZ-JX-APS:21,BE-EA:07}. In particular, the controllers follow a planar cascade architecture, where an outer-loop position PD computes the desired thrust and pitch angle, and an inner-loop attitude PD regulates the pitch dynamics.

In the experiments from the main text (using $1000$ different random seeds for disturbances and actuator noise), the action space is discretized into a uniform grid with $\du$ bins. All episodes start from the fixed initial state $\mathbf{x}_0$ (see Fig.~5(b) from the main text) and have length $T$. 
Statistical comparisons between methods are performed using a two-sided Mann--Whitney $U$ test~\cite{HBM-DRW:47} and effect sizes are quantified using Cliff's $\delta$~\cite{NC:93}.
The reduction in total episode cost relative to the PID baseline is statistically significant (as reported in the main text), with $\delta = -0.14$. 

In addition to the findings reported in the main text, we conduct further experiments using different primitives. The parameters for this additional set of experiments are in Table~\ref{tab:quadrotor_lqr_parameters}. The primitives are Gaussian policies as in Eq.~\eqref{primitives_quadrotor}, but in this case $\bv{u}^{1}_{k}$ and $\bv{u}^{2}_{k}$ are the control inputs generated by two Linear Quadratic Regulator (LQR) controllers. As in the previous quadrotor experiments, $\Sigma_{\mathrm{u}}$ models actuator noise. The LQRs have complementary objectives: (i) stabilization and (ii) trajectory tracking. 
{The stabilizing primitive has state and action cost matrices $Q_{\mathrm{s}}$ and $R_{\mathrm{s}}$, respectively}. These matrices are chosen so that the attitude error has the greatest weight in the cost. The stabilizing primitive uses as state reference the current position, with both velocities and pitch angle equal to zero, and as action reference ${(\bv{u}_{\mathrm{ref},k})}_{k=1,\ldots,T}$. Therefore, the drone controlled by this primitive does not pursue the trajectory. In contrast, the tracking primitive uses state and action cost matrices $Q_{\mathrm{t}}$ and $R_{\mathrm{t}}$, respectively, which emphasize position tracking, with reference ${(\bv{x}_{\mathrm{ref},k},\bv{u}_{\mathrm{ref},k})}_{k=1,\ldots,T}$.
The generative model is defined as in the previous quadrotor experiments.
In addition to actuator noise and white Gaussian noise in the dynamics, in this set of experiments we also introduce a single impulsive disturbance of fixed magnitude $m_{\mathrm{i}}$ applied at a randomly selected time step in each episode. For each random seed, all controllers are evaluated under identical disturbance realizations, ensuring a paired comparison across methods.  
As a baseline, we consider a standard discrete-time LQR controller whose state and action cost matrices, $Q_{\mathrm{b}}$ and $R_{\mathrm{b}}$, respectively, coincide with those used in the generative model.
Figure~\ref{fig::quadrotor_lqr} summarizes the results.
Panel~(a) shows representative trajectories obtained using our proposed model and, for comparison, the baseline LQR. The baseline LQR exhibits large transient deviations following the impulse. In contrast, our composed policy recovers faster from the impulsive disturbance due to the use of the stabilization primitive, resulting in a reduced excursion from the reference.
Panel~(b) reports the boxplots of the episode cost over $1000$ random seeds. Our model achieves a lower median episode cost than both individual primitives and the baseline LQR. In particular, it improves upon the stabilization primitive, which provides strong local regulation but does not explicitly track the reference, and upon the tracking primitive, which is optimized for nominal trajectory following but less robust to abrupt perturbations. Moreover, our model achieves a reduction in episode cost relative to the baseline LQR that is statistically significant ($p < 10^{-38}$), with $\delta = -0.34$.

\section{Proving the Statements}\label{sec:proofs}
\subsection*{Proof of Lemma~\ref{lem:prox_softmax}}
For any $\bv{x} \in \R^n$, we have
\begin{align}
\displaystyle \prox{ \entropybar - \frac{\norm{\cdot}^2}{2}}(\bv{x}) &:= \argmin_{\bv{z} \in \R^n} \left( \entropybar(\bv{z}) - \frac{\norm{\bv{z}}^2}{2} + \frac{1}{2}\norm{\bv{x} - \bv{z}}_2^2\right) \nonumber\\
&= \argmin_{\bv{z} \in \simplex{n}} \left( - \entropy(\bv{z}) - \frac{\norm{\bv{z}}^2}{2} + \frac{1}{2}\norm{\bv{x}}^2 - \bv{x}^\top \bv{z} + \frac{1}{2}\norm{\bv{z}}^2\right) \label{eq:lem_prox_entropy1}\\
&= \argmin_{\bv{z} \in \simplex{n}} \left( \sum_{i=1}^n z_i \ln{z_i} - \bv{x}^\top \bv{z}\right), \label{eq:lem_prox_entropy2}
\end{align}
where in the equality in Eq.~\eqref{eq:lem_prox_entropy1} we used the fact that, by definition of entropic barrier function, the minimum is obtained in $\simplex{n}$. Now, to solve the constrained problem~\eqref{eq:lem_prox_entropy2}, consider the Lagrangian function
\[
\mathcal{L}(\bv{z}, \lambda, \bv{\mu}) = -\bv{x}^\top \bv{z} +  \sum_{i=1}^n z_i \ln z_i + \lambda\left( \sum_{i=1}^n z_i - 1 \right) - \sum_{i=1}^n \mu_i z_i,
\]
where $\lambda \in \mathbb{R}$ is the multiplier associated with the equality constraint $\sum_{i=1}^n z_i = 1$, and $\mu_i \geq 0$, $i \in \until{n}$, are the multipliers associated with the inequality constraints $- z_i \leq 0$.
For all $i \in \until{n}$, the optimality KKT conditions are: (i) \emph{stationarity:} $\frac{\partial \mathcal{L}}{\partial z_i} = -x_i + (1 + \ln z_i) + \lambda - \mu_i = 0$; (ii) \emph{primal feasibility:} $\sum_{i=1}^n z_i = 1$; (iii) \emph{dual feasibility:} $\mu_i \geq 0$; and (iv) \emph{complementary slackness:} $\mu_i z_i = 0$.

\noindent
Now, since the term $(1 + \ln z_i)$ tends to $-\infty$ as any $z_i \to 0^+$, for the stationarity conditions to hold, the minimizer must lie strictly in the interior of the simplex. By complementary slackness, this implies that $\mu_i = 0$. Thus, the stationarity condition simplifies to
$$
-x_i + (1 + \ln z_i) + \lambda = 0 \quad\iff\quad z_i = \exp\left(x_i - \lambda - 1\right).
$$
Then, by primal feasibility we have 
\begin{align*}
\sum_{i=1}^n z_i = 1 \quad&\iff\quad \sum_{i=1}^n \exp\left({x_i - \lambda} - 1\right) = 1 \quad\iff
\quad \exp\left(-\lambda - 1\right)= \frac{1}{\sum_{j=1}^n \exp\left( {x_j} \right)}.
\end{align*}
Finally, the unique optimal solution is:
\[
z_i^\star = \frac{\exp\left( {x_i} \right)}{\sum_{j=1}^n \exp\left({x_j} \right)} =: \left(\softmax\left(\bv{x} \right) \right)_i, \quad i \in \until{n}.
\]
This concludes the proof.\qed

We can now give the proof of Lemma \ref{lem:softmax_optimization_connection}.

\subsection*{Proof of Lemma~\ref{lem:softmax_optimization_connection}}
First, note that problem~\eqref{eq:entropy_regularized_problem} is equivalent to the following composite optimization problem
\begin{equation}
\label{eq:composite_f_H_withnorm}
\min_{\bv{x} \in \R^n} f(\bv{x}) + g(\bv{x}),
\end{equation}
where $\ds f(\bv{x}) = \F(\bv{x}) + \eps \frac{\norm{\bv{x}}^2}{2}$ is continuously differentiable and $\ds g(\bv{x}) = \eps\bigg(\entropybar(\bv{x}) - \frac{\norm{\bv{x}}^2}{2}\bigg)$ is CCP.
The fact that $g$ is CCP follows from the following facts: (i) the effective domain (the set of points in the function's domain at which the function's value is finite) of $g$ is $\simplex{n}$; (ii) in the interior of the simplex, the Hessian of $g$ is positive semidefinite; (iii) convexity extends to the closure of the simplex due to $g$ being lower semicontinuous (the negative entropy is lower semicontinuous on the simplex, and the sum of a lower semicontinuous function and a continuous function is lower semicontinuous) and convex on the relative interior of the simplex.
Next, let $\bv{\xstar}$ be the minimizer of problem~\eqref{eq:composite_f_H_withnorm}. Then, by first-order necessary and sufficient optimality conditions for convex optimization problems~\cite{RTR:70} we have $\0_n \in \nabla f(\bv{x^{\star}}) + \partial g(\bv{\xstar})$.
Multiplying by $\eps^{-1}$ and adding and subtracting $\bv{\xstar}$ to the right-hand side of the above inclusion yields
\[
\begin{aligned}
\0_n \in (I_n + \eps^{-1}\partial g)(\bv{\xstar}) + \eps^{-1} \nabla f(\bv{\xstar}) - \bv{\xstar} &  \iff  
\bv{x}^\star-\eps^{-1}\nabla f(\bv{x}^\star) \in (I_n+\eps^{-1}\partial g)(\bv{x}^\star)\\
& \iff  \bv{\xstar} \in (I_n {+} \eps^{-1}\partial g)^{-1}\left(\bv{\xstar} {-} \eps^{-1} \bigl(\nabla \F(\bv{\xstar}) {+}\eps \bv{\xstar} \bigr)\right).
\end{aligned}
\]
Recalling that $\prox{\eps^{-1} g}= (I_n + \eps^{-1}\partial g)^{-1}$ and, being by assumption $g$ CCP, then $\prox{\eps^{-1} g}$ is single-valued~\cite{NP-SB:14} and, by Lemma~\ref{lem:prox_softmax}, we have $\prox{\eps^{-1} g} = \softmax$. Therefore,
$$
\bv{\xstar} = \softmax\left(- \eps^{-1} \nabla \F(\bv{\xstar})\right),
$$
that is, $\bv{\xstar}$ is an equilibrium point of the proximal gradient flow dynamics in Eq.~\eqref{eq:softmax_gf}. This establishes the proof for items~\ref{item1:op_soft} and~\ref{item2:op_soft}.
\qed

\subsection*{Proof of Lemma~\ref{lem:positive_softmax_gf}} First note that $\simplex{n} = \R^n_{\geq 0} \cap \setdef{\bv{x} \in \R^n}{\sum_{i = 1}^{n} x_i=1}$. Then, by Nagumo's Theorem~\cite{MN:1942} (see also~\cite[Exercise 3.12]{FB:24-CTDS}), to prove that $\simplex{n}$ is forward invariant for the dynamics in Eq.~\eqref{eq:softmax_gf} we need to show that (i) the positive orthant is forward invariant, (ii) the dynamics preserve the total mass constraint, that is $\sum_{i = 1}^n \dot x_i = 0$ for all $\bv{x} \in \R^n$ such that $\sum_{i = 1}^{n} x_i=1$.
By Nagumo's Theorem~\cite{MN:1942}, the positive orthant is forward invariant for a vector field $f$ if and only if $f_i(\bv{x}) \geq 0$, for all $\bv{x} \in \R_{\geq 0}^n$ such that $x_i =0$.
Let us consider the softmax gradient flow in Eq.~\eqref{eq:softmax_gf} written in components
\begin{align*}
\tau \dot{x}_i &= - x_i + \softmax(-\eps^{-1}\nabla \F(\bv{x}))_i \\ 
&\ds= - x_i + \frac{\exp\bigl(-\eps^{-1}\nabla_i \F(\bv{x})\bigr)}{\sum_{j=1}^n \exp\bigl(-\eps^{-1}\nabla_j \F(\bv{x})\bigr)}:= \subscr{\F}{{sm,$i$}}(\bv{x}), \quad i \in \until{n}.
\end{align*}
Then, for all $\bv{x} \in \R_{\geq 0}^n$ such that $x_i =0$ we have
$
\ds \subscr{\F}{sm,$i$}(\bv{x}) = \frac{\operatorname{exp}\left(-\eps^{-1}\nabla_i \F(\bv{x})\right)}{\sum_{j=1}^n \operatorname{exp}\left(-\eps^{-1}\nabla_j \F(\bv{x})\right)} > 0,
$
for each $i$.
Next, let $\bv{x} \in \R^n$ such that $\sum_{i = 1}^{n} x_i=1$. We have
$$
\ds \sum_{i=1}^n \tau \dot{x}_i = - \sum_{i=1}^n x_i + \sum_{i=1}^n \softmax\left(- \eps^{-1}\nabla \F(\bv{x})\right)_i = -1 + 1 = 0.
$$
This concludes the proof. \qed

\subsection*{Proof of Lemma~\ref{lem:equality_Dkl}}
Let $\bv{x}(0) \in \simplex{n}$ and denote by $\bv{x}(t) = \phi_t(\bv{x}(0))$ the corresponding trajectory of Eq.~\eqref{eq:softmax_gf}. By Lemma~\ref{lem:positive_softmax_gf}, the simplex $\simplex{n}$ is forward invariant, so $\bv{x}(t)\in\simplex{n}$ for all $t \ge 0$. In particular $V(\bv{x}(t)) = \F(\bv{x}(t)) - \eps \entropy(\bv{x}(t))$. Furthermore, since the softmax map has strictly positive entries, we have $\bv{x}(t)\in\operatorname{ri}(\simplex{n})$ for all $t>0$.
Since $\F$ is continuously differentiable on $\simplex{n}$ and the entropy term is continuously differentiable on $\operatorname{ri}(\simplex{n})$, it follows that $V$ is continuously differentiable on $\operatorname{ri}(\simplex{n})$. Therefore, for $t>0$, the map $t \mapsto V(\bv{x}(t))$ is differentiable and, by the chain rule,
\begin{equation*}
\frac{d}{dt}V(\bv{x}) = \nabla \F(\bv{x})^\top \dot{\bv{x}} + \eps \sum_{i=1}^n (\ln x_i + 1)\,\dot x_i.
\end{equation*}
Substituting the flow dynamics in Eq.~\eqref{eq:softmax_gf} into the above equation, we obtain
\begin{equation*}
\tau \frac{d}{dt}V(\bv{x}) = \nabla \F(\bv{x})^\top (S(\bv{x}; \eps)-\bv{x}) + \eps \sum_{i=1}^n (\ln x_i + 1)(S_i(\bv{x}; \eps)-x_i).
\end{equation*}
Since both $\bv{x}$ and $S(\bv{x}; \eps)$ belong to the simplex, we know $\sum_{i=1}^{n} (x_i - S_i(\bv{x}; \eps)) = 0$, and therefore
\beq
\label{eq:dot_G}
\tau \frac{d}{dt}V(\bv{x}) = \nabla \F(\bv{x})^\top (S(\bv{x}; \eps)-\bv{x}) + \eps \sum_{i=1}^n \ln x_i (S_i(\bv{x}; \eps)-x_i).
\eeq
Since
\begin{equation*}
\ln S_i(\bv{x}; \eps) = \ln \Biggl( \frac{\exp(-\eps^{-1}\nabla_i \F(\bv{x}))}{\sum_{j=1}^n \exp(-\eps^{-1}\nabla_j \F(\bv{x}))}\Biggr)
= -\eps^{-1}\nabla_i \F(\bv{x}) - \ln \Biggl(\sum_{j=1}^n \exp(-\eps^{-1}\nabla_j \F(\bv{x}))\Biggr),
\end{equation*}
we obtain
\begin{equation}
\label{eq:grad_F_i}
\nabla_i \F(\bv{x})
=
-\eps \ln S_i(\bv{x}; \eps) - \eps \ln \biggl(\sum_{j=1}^n \exp(-\eps^{-1}\nabla_j \F(\bv{x})) \biggr),
\qquad i=1,\dots,n.
\end{equation}
Moreover, by definition,
$
\nabla \F(\bv{x})^\top(\bv{x}-S(\bv{x}; \eps))
=\sum_{i=1}^n \nabla_i \F(\bv{x})
(x_i-S_i(\bv{x}; \eps)).
$
Substituting Eq.~\eqref{eq:grad_F_i} into this last equation yields
\begin{align}
\nabla \F(\bv{x})^\top\bigl(\bv{x}-S(\bv{x};\eps)\bigr)
&= -\eps \sum_{i=1}^n
   \bigl(x_i-S_i(\bv{x};\eps)\bigr)
   \ln S_i(\bv{x};\eps)
   \nonumber\\
&\quad
   -\eps
   \ln\left(
      \sum_{j=1}^n
      \exp\bigl(-\eps^{-1}\nabla_j\F(\bv{x})\bigr)
   \right)
   \sum_{i=1}^n
   \bigl(x_i-S_i(\bv{x};\eps)\bigr)
   \nonumber\\
&= -\eps \sum_{i=1}^n
   \bigl(x_i-S_i(\bv{x};\eps)\bigr)
   \ln S_i(\bv{x};\eps).
\label{eq:lemma6}
\end{align}
where in the last equality we used again the fact that since both $\bv{x}$ and $S(\bv{x}; \eps)$ belong to the simplex, the second term vanishes.
Substituting Eq.~\eqref{eq:lemma6} into Eq.~\eqref{eq:dot_G}, we obtain
\begin{align*}
\tau \frac{d}{dt}V(\bv{x})
&= \eps \sum_{i=1}^n (x_i-S_i(\bv{x}; \eps))\ln S_i(\bv{x}; \eps) + \eps \sum_{i=1}^n (S_i(\bv{x}; \eps)-x_i)\ln x_i \\
&= -\eps \sum_{i=1}^n (x_i-S_i(\bv{x}; \eps)) \ln \frac{x_i}{S_i(\bv{x}; \eps)}\\
&= -\eps \Bigl( \DKL{\bv{x}}{S(\bv{x}; \eps)} + \DKL{S(\bv{x}; \eps)}{\bv{x}} \Bigr).
\end{align*}
This concludes the proof. \qed

\subsection*{Proof of Theorem~\ref{thm:exp_convergence}}
Pick $\bv{x}(0) \in \simplex{n}$ and denote by $\bv{x}(t) = \phi_t(\bv{x}(0))$ the solution of the dynamics in Eq.~\eqref{eq:softmax_gf} with initial condition $\bv{x}(0)$. By Lemma~\ref{lem:positive_softmax_gf}, the simplex $\simplex{n}$ is forward invariant, so $\bv{x}(t)\in\simplex{n}$ for all $t\ge 0$. Moreover, since the softmax map has strictly positive entries for each argument, it follows that $\bv{x}(t)\in\operatorname{ri}(\simplex{n})$ for all $t > 0$ and $V(\bv{x}(t)) = \F(\bv{x}(t)) - \eps \entropy(\bv{x}(t))$.

We next characterize the equilibria of the dynamics. By item~\ref{item2:op_soft} of Lemma~\ref{lem:softmax_optimization_connection}, the equilibrium points of Eq.~\eqref{eq:softmax_gf} coincide with the minimizers of the entropy-regularized problem in Eq.~\eqref{eq:entropy_regularized_problem}. By Remark~\ref{remark:exist_uniq}, this problem admits a unique minimizer $\bv{x}^\star \in \simplex{n}$, and hence the dynamics admits a unique equilibrium $\bv{x}^\star$.
To prove convergence of the softmax gradient flow in Eq.~\eqref{eq:softmax_gf} to this equilibrium, we derive a differential inequality for the entropy-regularized objective $V(\bv{x})$, which yields exponential 
decay of the optimality gap and exponential convergence of the trajectories.

By the smooth best response characterization of the softmax map (see Eq.~\eqref{eq:best_response_map}), we have that, for any $\bv{v}\in\R^n$, $\softmax(\bv{v}) = \argmax_{\bv{u}\in\simplex{n}}
\left\{ \bv{v}^\top \bv{u} + \entropy(\bv{u}) \right\}$.
Applying this with $\bv{v}=-\eps^{-1}\nabla \F(\bv{x})$, we obtain
\begin{align*}
S(\bv{x};\eps) =\softmax(-\eps^{-1}\nabla \F(\bv{x}))
& = \argmax_{\bv{u}\in\simplex{n}}
\left\{
-\eps^{-1}\nabla \F(\bv{x})^\top \bv{u}
+ \entropy(\bv{u})
\right\}\\
& =
\argmin_{\bv{u}\in\simplex{n}}
\left\{
\nabla \F(\bv{x})^\top \bv{u} + \eps\sum_{i=1}^n u_i\ln u_i
\right\},
\end{align*}
where in the last step we multiplied the objective by $\eps>0$, which does not change the minimizer.
Choosing $\bv{u}=\bv{x}^\star$ yields
\begin{equation}
\label{eq:ineq_1}
\nabla \F(\bv{x})^\top S(\bv{x}; \eps)+\eps\sum_{i=1}^n S_i(\bv{x}; \eps)\ln S_i(\bv{x}; \eps)
\le
\nabla \F(\bv{x})^\top \bv{x}^\star+\eps\sum_{i=1}^n x_i^\star\ln x_i^\star.
\end{equation}
By convexity of $\F$, it holds $\F(\bv{x})-\F(\bv{x}^\star)\le \nabla \F(\bv{x})^\top(\bv{x}-\bv{x}^\star)$. Adding $\eps\bigl(\sum_{i=1}^n x_i\ln x_i -\sum_{i=1}^n x_i^\star\ln x_i^\star\bigr)$ to both sides of this inequality yields
\begin{align}
V(\bv{x})-V(\bv{x}^\star)
&\le \nabla \F(\bv{x})^\top \bv{x} + \eps \sum_{i=1}^n x_i\ln x_i - \biggl(\nabla \F(\bv{x})^\top \bv{x}^\star + \eps \sum_{i=1}^n x_i^\star\ln x_i^\star\biggr), \nonumber\\
&\le\nabla \F(\bv{x})^\top(\bv{x}-S(\bv{x}; \eps)) + \eps\sum_{i=1}^n \bigl(x_i\ln x_i-S_i(\bv{x}; \eps)\ln S_i(\bv{x}; \eps)\bigr), \label{eq:ineq_3}
\end{align}
where in the last inequality we have used the inequality in Eq.~\eqref{eq:ineq_1}. 
From Eq.~\eqref{eq:lemma6} in the proof of Lemma~\ref{lem:equality_Dkl}, we have
$
\nabla \F(\bv{x})^\top(\bv{x}-S(\bv{x}; \eps)) = -\eps \sum_{i=1}^n (x_i - S_i(\bv{x}; \eps)) \ln S_i(\bv{x}; \eps).
$
Substituting this into the inequality in Eq.~\eqref{eq:ineq_3} yields
\begin{align}
V(\bv{x})-V(\bv{x}^\star)
& \le -\eps \sum_{i=1}^n (x_i - S_i(\bv{x}; \eps)) \ln S_i(\bv{x}; \eps) + \eps\sum_{i=1}^n \bigl(x_i\ln x_i-S_i(\bv{x}; \eps)\ln S_i(\bv{x}; \eps)\bigr)  \nonumber \\ 
&= \eps\sum_{i=1}^n x_i\ln\frac{x_i}{S_i(\bv{x}; \eps)} = \eps\, \DKL{\bv{x}}{S(\bv{x}; \eps)}.
\label{eq:gap_bound_dkl}
\end{align}
By Lemma~\ref{lem:equality_Dkl}, the inequality in Eq.~\eqref{eq:gap_bound_dkl}, and using the nonnegativity of $\DKL{S(\bv{x}; \eps)}{\bv{x}}$, we obtain
\begin{equation*}
V(\bv{x})-V(\bv{x}^\star) \le \eps \DKL{\bv{x}}{S(\bv{x}; \eps)} \le \eps \Bigl( \DKL{\bv{x}}{S(\bv{x}; \eps)} + \DKL{S(\bv{x}; \eps)}{\bv{x}} \Bigr) = - \tau \frac{d}{dt} V(\bv{x}).
\end{equation*}
Therefore, $\frac{d}{dt} V(\bv{x})\le -\frac{1}{\tau} \bigl( V(\bv{x})-V(\bv{x}^\star) \bigr)$ and, by Gr\"onwall's inequality, we have
\begin{equation*}
V(\bv{x}(t))-V(\bv{x}^\star) \le \e^{-t/\tau}\bigl(V(\bv{x}(0))-V(\bv{x}^\star)\bigr)
\qquad\text{for all }t\ge 0.
\end{equation*}
This proves statement~\ref{item1:proof_conv}.

We now prove statement~\ref{item2:proof_conv}. 
By convexity of $\F$,
\begin{equation}
\label{eq:F_convexity}
\F(\bv{x})-\F(\bv{x}^\star)\ge \nabla \F(\bv{x}^\star)^\top(\bv{x}-\bv{x}^\star).
\end{equation}
Since $\bv{x}^\star$ is the minimizer of $V$ on $\simplex{n}$ and lies in $\operatorname{ri}(\simplex{n})$ (see Remark~\ref{remark:exist_uniq}), the first-order optimality condition implies that there exists $\lambda^\star\in\R$ such that
$
\nabla \F(\bv{x}^\star)+\eps(\ln \bv{x}^\star+\1_n)=\lambda^\star \1_n.
$
Rearranging, $\nabla \F(\bv{x}^\star) = \lambda^\star \1_n - \eps (\ln \bv{x}^\star + \1_n)$.
Substituting into the inequality in Eq.~\eqref{eq:F_convexity}, we have
\begin{align*}
\F(\bv{x}) - \F(\bv{x}^\star)
&\ge (\lambda^\star \1_n - \eps (\ln \bv{x}^\star + \1_n))^\top (\bv{x} - \bv{x}^\star) = \lambda^\star \1_n^\top (\bv{x} - \bv{x}^\star)
- \eps (\ln \bv{x}^\star + \1_n)^\top (\bv{x} - \bv{x}^\star)\\
&\ge - \eps (\ln \bv{x}^\star + \1_n)^\top (\bv{x} - \bv{x}^\star),
\end{align*}
where in the last inequality we used the fact that since $\bv{x}, \bv{x}^\star \in \simplex{n}$, we have $\1_n^\top (\bv{x} - \bv{x}^\star) = 1 - 1 = 0$.
Next, we compute
\begin{align*}
V(\bv{x}) - V(\bv{x}^\star)
&= \F(\bv{x}) - \F(\bv{x}^\star)
+ \eps \sum_{i=1}^n \big(x_i \ln x_i - x_i^\star \ln x_i^\star\big) \\
&\ge - \eps (\ln \bv{x}^\star + \1_n)^\top (\bv{x} - \bv{x}^\star)
+ \eps \sum_{i=1}^n \big(x_i \ln x_i - x_i^\star \ln x_i^\star\big)\\
&=- \eps \sum_{i=1}^n \ln x_i^\star (x_i - x_i^\star) + \eps \sum_{i=1}^n \big(x_i \ln x_i - x_i^\star \ln x_i^\star\big),
\end{align*}
where in the last equality we used again $\1_n^\top (\bv{x} - \bv{x}^\star) = 0$.
Rearranging the terms,
\begin{align*}
V(\bv{x}) - V(\bv{x}^\star)
&\ge \eps \sum_{i=1}^n
\Bigl[
x_i \ln x_i - x_i^\star \ln x_i^\star
- \ln x_i^\star (x_i - x_i^\star)
\Bigr] \\
&= \eps \sum_{i=1}^n
\Bigl[
x_i \ln x_i - x_i \ln x_i^\star
\Bigr]
= \eps \sum_{i=1}^n x_i \ln \frac{x_i}{x_i^\star} = \eps \, \DKL{\bv{x}}{\bv{x}^\star}\\
&\ge \frac{\eps}{2}\|\bv{x}-\bv{x}^\star\|_2^2,
\end{align*}
where in the last inequality we used Pinsker's inequality (see, e.g.,~\cite{IC-JK:11,MT:14}), that is $
\DKL{\bv{x}}{\bv{x}^\star}\ge \frac12\|\bv{x}-\bv{x}^\star\|_1^2$, considering that the $D_{\mathrm{KL}}$ is defined using natural logarithms, and the bound $\frac12\|\bv{x}-\bv{x}^\star\|_1^2\ge \frac12\|\bv{x}-\bv{x}^\star\|_2^2$. 
Then, by using the exponential estimate for the cost in Eq.~\eqref{eq:exp_cost_bound}, we obtain
\begin{equation*}
\|\bv{x}(t)-\bv{x}^\star\|_2^2
\le
\frac{2}{\eps}\bigl(V(\bv{x}(t))-V(\bv{x}^\star)\bigr)
\le
\frac{2}{\eps}\e^{-t/\tau}\bigl(V(\bv{x}(0))-V(\bv{x}^\star)\bigr),
\end{equation*}
and therefore
\begin{equation*}
\|\bv{x}(t)-\bv{x}^\star\|_2
\le
\sqrt{\frac{2}{\eps}\bigl(V(\bv{x}(0))-V(\bv{x}^\star)\bigr)}\,\e^{-t/(2\tau)}.
\end{equation*}
This concludes the proof. \qed

\subsection*{Proof of Lemma~\ref{lem::gradient}}
We begin with item~\ref{item1:equivalence_ops_DKL}. By Lemma~\ref{lem:chain_rule_KL}, the objective of the problem in Eq.~\eqref{eqn:gateframe_SI} becomes:
\begin{equation}
\label{eq:cost_reformulation_DKL}
\DKL{{\policy{k}{k-1}}\!}{\!{\refpolicy{k}{k-1}}}  {+} \E_{\policy{k}{k-1}}\big[\!\DKL{{\plant{k}{k-1}}\!}{\!{\refplant{k}{k-1}}}\big] {-} \eps \entropy{(\weights{k})}.
\end{equation}
Substituting the constraint $\displaystyle \policy{k}{k-1} = \sum_{\alpha =1}^{\np}\weight{k}{\alpha}\primitive{k}{k-1}{{\alpha}}$ into Eq.~\eqref{eq:cost_reformulation_DKL} yields, for the first term:
\begin{align*}
\DKL{\policy{k}{k-1}\!}{\!\refpolicy{k}{k-1}} &= \int \sum_{\alpha=1}^{\np} \weight{k}{\alpha}\primitive{k}{k-1}{\alpha} \ln \Bigg(\frac{\sum_{\beta =1}^{\np}\weight{k}{\beta}\primitive{k}{k-1}{\beta}}{\refpolicy{k}{k-1}}\Bigg) d\bv{u}_k\\
& = \!\sum_{\alpha=1}^{\np}\!\weight{k}{\alpha}\E_{\primitive{k}{k-1}{\alpha}}\bigg[\!\ln\!\sum_{\beta =1}^{\np}\!\weight{k}{\beta}\primitive{k}{k-1}{\beta}{-}\ln\refpolicy{k}{k-1}\bigg].
\end{align*}
For the second term we have
\begin{align*}
&\E_{\policy{k}{k-1}} \big[\DKL{\plant{k}{k-1}}{\refplant{k}{k-1}}\big]\\
&\qquad = 
\sum_{\alpha=1}^{\np}\weight{k}{\alpha}\E_{\primitive{k}{k-1}{\alpha}}\left[\DKL{\plant{k}{k-1}}{\refplant{k}{k-1}}\right].
\end{align*}
Embedding the simplex constraint of the problem in Eq.~\eqref{eqn:gateframe_SI} in the cost via the indicator function yields the desired conclusion.

To prove item~\ref{item2:gradient_op_reduced_DKL}, first note that Assumption~\ref{assum_primitives} implies that the function $\F(\weights{k})$ is smooth in $\weights{k}$. Indeed, the only potential source of discontinuity is the term $\ln \sum_{\beta =1}^{\np}\weight{k}{\beta}\primitive{k}{k-1}{\beta}$.  However, since this sum is always non-negative for $\weights{k} \in \simplex{\np}$, the function $\F(\weights{k})$ remains well-defined.
Next, we compute the gradient of $\F(\weights{k})$ w.r.t. a generic weight, say $\weight{k}{\gamma}$. We have
\[
\begin{aligned}
\frac{\partial}{\partial \weight{k}{\gamma}} \F(\weights{k}) {=} & \frac{\partial}{\partial \weight{k}{\gamma}}\! \Bigg( \!\sum_{\alpha=1}^{\np}\weight{k}{\alpha} \E_{\primitive{k}{k-1}{\alpha}}\!\Big[\!\DKL{\plant{k}{k-1}\!\!}{\!\!\refplant{k}{k-1}}{-}\!\ln\refpolicy{k}{k-1}\Big] \!\Bigg) \\
&\ + \frac{\partial}{\partial \weight{k}{\gamma}} \left( \sum_{\alpha =1}^{\np}\weight{k}{\alpha} \E_{\primitive{k}{k-1}{\alpha}}\bigg[ \ln\sum_{\beta =1}^{\np}\weight{k}{\beta}\primitive{k}{k-1}{\beta}\bigg] \right)\\
= & \E_{\primitive{k}{k-1}{\gamma}} \bigg[\DKL{\plant{k}{k-1}}{\refplant{k}{k-1}} - \ln\refpolicy{k}{k-1} \\
&\qquad\qquad\qquad+ \ln \sum_{\beta =1}^{\np}\weight{k}{\beta}\primitive{k}{k-1}{\beta} \bigg] + 1,
\end{aligned}
\]
where we used the following derivations:
\[
\begin{aligned}
&\frac{\partial}{\partial \weight{k}{\gamma}} \Bigg(\sum_{\alpha =1}^{\np}\weight{k}{\alpha} \E_{\primitive{k}{k-1}{\alpha}}\bigg[ \ln \sum_{\beta =1}^{\np}\weight{k}{\beta}\primitive{k}{k-1}{\beta}\bigg]\Bigg)\\
&\quad= \frac{\partial}{\partial \weight{k}{\gamma}} \Bigg( \int \sum_{\alpha=1}^{\np}\weight{k}{\alpha} \primitive{k}{k-1}{\alpha} \ln \sum_{\beta =1}^{\np}\weight{k}{\beta}\primitive{k}{k-1}{\beta} d\bv{u}_k\Bigg)\\
&\quad= \int \frac{\partial}{\partial \weight{k}{\gamma}} \Bigg( \sum_{\alpha =1}^{\np}\weight{k}{\alpha} \primitive{k}{k-1}{\alpha}\ln \sum_{\beta =1}^{\np}\weight{k}{\beta}\primitive{k}{k-1}{\beta}\Bigg) d\bv{u}_k\\
&\quad=\int \primitive{k}{k-1}{\gamma}\ln \sum_{\beta =1}^{\np}\weight{k}{\beta}\primitive{k}{k-1}{\beta} + \sum_{\alpha =1}^{\np}\weight{k}{\alpha} \primitive{k}{k-1}{\alpha} \frac{\primitive{k}{k-1}{\gamma}}{\sum\limits_{\beta=1}^{\np}\weight{k}{\beta}\primitive{k}{k-1}{\beta}} d\bv{u}_k\\
&\quad=\int \primitive{k}{k-1}{\gamma} \Bigg( \ln \sum_{\beta =1}^{\np}\weight{k}{\beta}\primitive{k}{k-1}{\beta} + 1 \Bigg) d\bv{u}_k\\
&\quad=\E_{\primitive{k}{k-1}{\gamma}}\left[\ln \sum_{\beta =1}^{\np}\weight{k}{\beta}\primitive{k}{k-1}{\beta} \right] + 1.
\end{aligned}
\]

To prove item~\ref{item3:smooth_constant}, we compute the Hessian of $\F(\weights{k})$. It suffices to show that, for every
$\weights{k}\in\simplex{\np}$,
$\0_{\np\times\np}
\preceq
\nabla^2\F(\weights{k})
\preceq
\subscr{L}{$\F$}I_{\np}$.
The first inequality establishes convexity of $\F$ on $\simplex{\np}$, whereas the second establishes that $\F$ is $\subscr{L}{$\F$}$-smooth on $\simplex{\np}$.
We compute
\begin{align*}
\left(\nabla^2 \F(\weights{k})\right)_{\alpha \gamma} &= \frac{\partial}{\partial \weight{k}{\alpha}}\left(\E_{\primitive{k}{k-1}{\gamma}}\left[\ln\sum_{\beta = 1}^{\np}\weight{k}{\beta}\primitive{k}{k-1}{\beta} \right]\right) \\
&= \int \primitive{k}{k-1}{\gamma} \frac{\partial}{\partial \weight{k}{\alpha}} \left(\ln\sum_{\beta = 1}^{\np}\weight{k}{\beta}\primitive{k}{k-1}{\beta}\right)d\bv{u}_k \\
&= \int \primitive{k}{k-1}{\gamma} \frac{\primitive{k}{k-1}{\alpha}}{\weights{k}^\top\primitives{k}{k-1}}d\bv{u}_k\\
&=\E_{\primitive{k}{k-1}{\gamma}}\left[\frac{\primitive{k}{k-1}{\alpha}}{\weights{k}^\top\primitives{k}{k-1}}\right].
\end{align*}

Note that, for any $\bv{z}\in \R^{\np}$, $\bv{z} \neq \0_{\np}$, it holds 
\begin{align*}
\bv{z}^\top \nabla^2 \F(\weights{k}) \bv{z} &= \sum_{\alpha=1}^{\np}\sum_{\gamma=1}^{\np} z_\alpha z_\gamma \E_{\primitive{k}{k-1}{\gamma}}\left[\frac{\primitive{k}{k-1}{\alpha}}{\weights{k}^\top\primitives{k}{k-1}}\right]\\
&= \sum_{\alpha=1}^{\np}\sum_{\gamma=1}^{\np} z_\alpha z_\gamma \int \frac{\primitive{k}{k-1}{\gamma}\primitive{k}{k-1}{\alpha}}{\weights{k}^\top\primitives{k}{k-1}}d\bv{u}_k\\
&= \int\sum_{\alpha=1}^{\np}\sum_{\gamma=1}^{\np} z_\alpha z_\gamma \frac{\primitive{k}{k-1}{\gamma}\primitive{k}{k-1}{\alpha}}{\weights{k}^\top\primitives{k}{k-1}}d\bv{u}_k\\
&= \int
\frac{\bv{z}^\top \primitives{k}{k-1}\primitives{k}{k-1}^\top \bv{z}}{\weights{k}^\top\primitives{k}{k-1}}d\bv{u}_k\\
&= \int
\frac{(\bv{z}^\top \primitives{k}{k-1})^2}{\weights{k}^\top\primitives{k}{k-1}}d\bv{u}_k \geq0,
\end{align*}
where the last inequality follows being $\weights{k}^\top\primitives{k}{k-1} >0$ for all $\weights{k} \in \simplex{\np}$.
Therefore, the Hessian $\nabla^2 \F(\weights{k})$ is a positive-semidefinite matrix, which in turn ensures that the map $\F(\weights{k})$ is convex. Next, recall that for any symmetric matrix, the largest eigenvalue is given by the maximum of the Rayleigh quotient, that is
\[
\subscr{\lambda}{max}\left(\nabla^2 \F(\weights{k})\right) = \max_{\norm{\bv{y}}_2 = 1} \bv{y}^\top \nabla^2 \F(\weights{k}) \bv{y}.
\]
Then, to prove our statement, it suffices to show $\bv{y}^\top \nabla^2 \F(\weights{k}) \bv{y} \leq \subscr{L}{$\F$}$, for all $\bv{y} \in \R^{\np}$ such that $\norm{\bv{y}}_2 = 1$.
We have
\begin{align*}
\bv{y}^{\top}\nabla^2\F(\weights{k})\bv{y}
&=
\int
\frac{
\bigl(\bv{y}^{\top}\primitives{k}{k-1}\bigr)^2
}{
\weights{k}^{\top}\primitives{k}{k-1}
}
\,d\bv{u}_k
\\
&\leq
\frac{1}{\subscr{\pi}{min}}
\int
\bigl(\bv{y}^{\top}\primitives{k}{k-1}\bigr)^2
\,d\bv{u}_k
\\
&\leq
\frac{\|\bv{y}\|_2^2}{\subscr{\pi}{min}}
\int
\|\primitives{k}{k-1}\|_2^2
\,d\bv{u}_k
\\
&=
\frac{1}{\subscr{\pi}{min}}
\sum_{\alpha=1}^{\np}
\int
\bigl(\primitive{k}{k-1}{\alpha}\bigr)^2
\,d\bv{u}_k
\\
&\leq
\frac{\subscr{\pi}{max}}{\subscr{\pi}{min}}
\sum_{\alpha=1}^{\np}
\int
\primitive{k}{k-1}{\alpha}
\,d\bv{u}_k
=
\np\frac{\subscr{\pi}{max}}{\subscr{\pi}{min}}.
\end{align*}
where in the first inequality we used the fact that $\weights{k}^{\top}\primitives{k}{k-1}=\sum_{\alpha=1}^{\np}\weight{k}{\alpha}\primitive{k}{k-1}{\alpha}
\geq \subscr{\pi}{min}
\sum_{\alpha=1}^{\np}\weight{k}{\alpha}
= \subscr{\pi}{min}$, while in the second inequality we applied the Cauchy-Schwarz inequality.
This concludes the proof. \qed

\newpage

\section{Supplemental Figures}\label{sec:figures}

\begin{figure}[h]
	\centering
\noindent\includegraphics[width=0.69\linewidth]{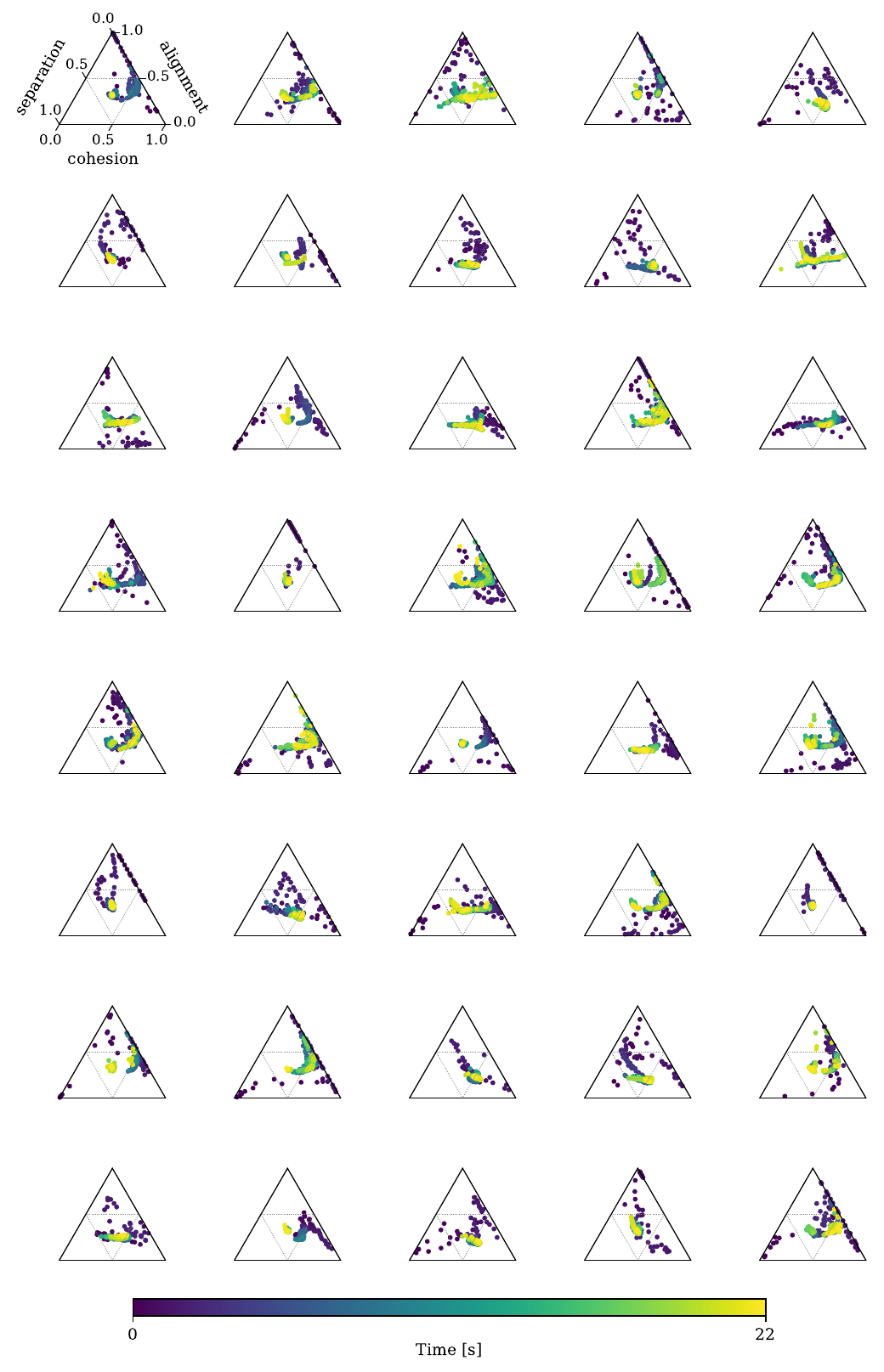}
\caption{Time evolution of the weights of all the boids in Fig.~3(b) (where the top-left simplex is shown). Over time, the weights become approximately equal.
As highlighted in the main text, this suggests a link between balanced (uniform) use of forces and the onset of polarization.
}
\label{fig:boid_polarization_all_weights}
\end{figure}

\clearpage
\begin{figure}
	\centering
\noindent\includegraphics[width=1\linewidth]{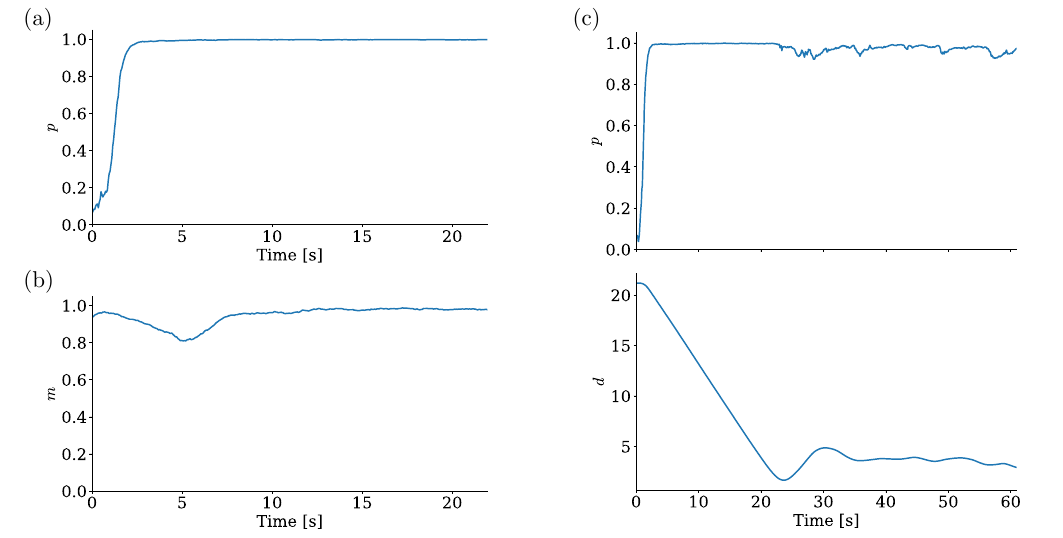}
\caption{Metrics for the collective behavior experiments in Fig.~3. 
See Sec.~\ref{sec:experiments_details} for the definitions of the order parameters.
(a) Evolution of the polarization order parameter metric for the experiment in Fig.~3(b). Over time, the order parameter becomes equal to $1$, indicating the onset of polarization.
(b) Evolution of the milling order parameter metric for the experiment in Fig.~3(c). An order parameter equal to $1$ indicates the onset of milling.
(c) Evolution of the polarization order parameter (top) and distance of the flock from goal (bottom) for the experiment in Fig.~3(d). The panels show that the flock is polarized and soft-controlled to the target destination.
}
\label{fig:boids_metrics}
\end{figure}

\clearpage
\begin{figure}
	\centering
\noindent\includegraphics[width=1\linewidth]{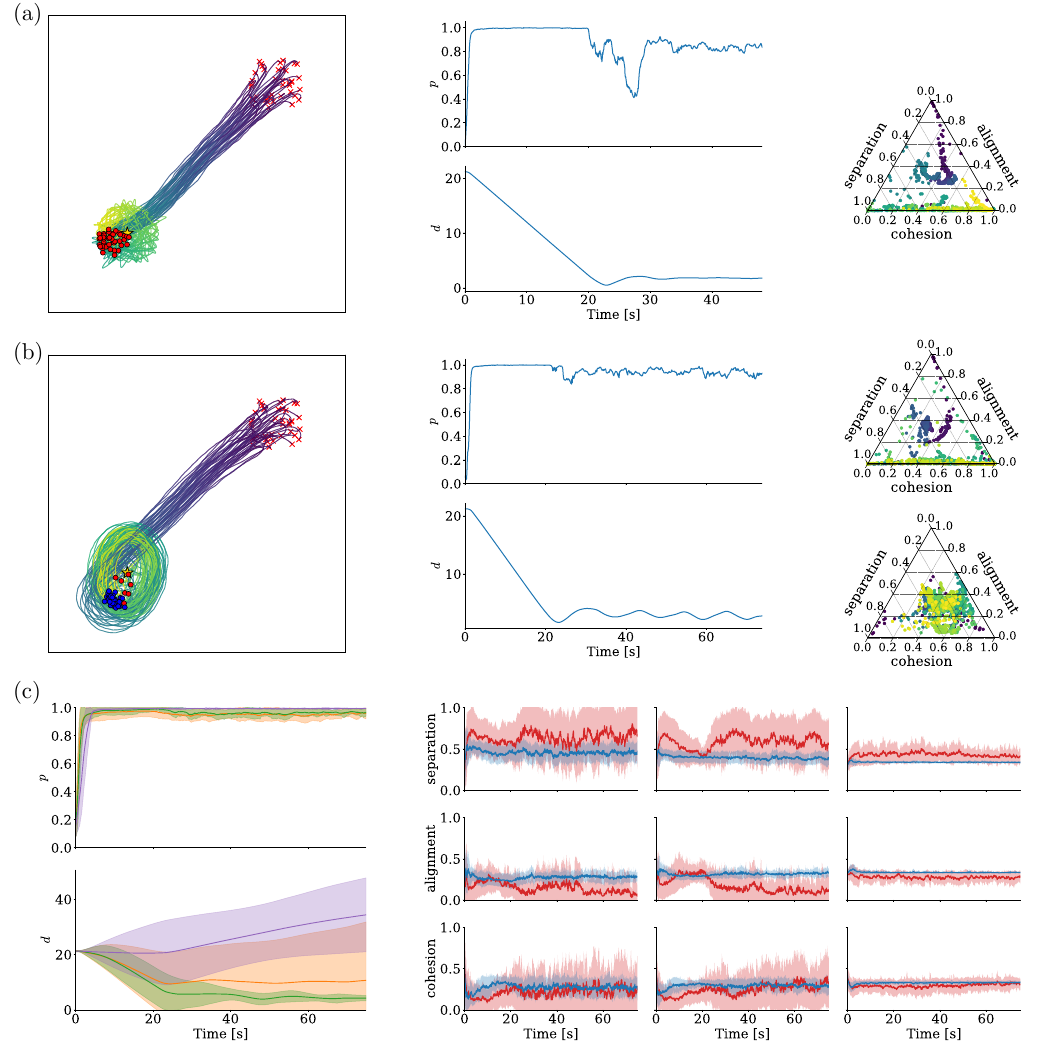}
\caption{Additional collective behavior experiments.
(a)
(Left) Group behavior when $100\%$ of the $N$ boids (red-filled circles) are informed about a goal position (star). Boids start from random initial conditions (crosses).
(Middle) The distance from the target decreases and the group exhibits polarization. This indicates a goal-directed yet polarized emerging behavior.
(Right) Time evolution of the optimal primitives' weights from our model. The weights are shown for one boid representative of all the (informed) boids in the group. 
(b)
(Left) Group behavior when $20\%$ of the $N$ boids are informed (red-filled circles), while the others (blue-filled circles) are not.
(Middle) Confirming the results from the main text, the distance from the target decreases (bottom) and the group exhibits polarization (top). This indicates a goal-directed yet polarized emerging behavior.
(Right) Time evolution of the optimal primitives' weights from our model for one representative informed boid (top) and one representative uninformed boid (bottom). As highlighted in the main text, the weights of the uninformed and informed boids behave differently over time. 
(c)
(Left) Mean and standard deviation of the group metrics across $30$ simulations in the setting of Fig.~3(d) from main text, with temperature $\eps$ equal to $0.1$ (orange), $1$ (green), and $10$ (purple).
While polarization remains high across all temperatures (top), varying $\eps$ may hinder the navigation of boids to the goal destination (bottom). Results suggest the existence of a sweetspot for $\eps$ that enables goal-directed and cohesive behavior.
(Right) Mean and standard deviation of weights across $30$ seeds. Plots show the weights for a representative informed boid (red) and a representative uninformed boid (blue), when $\eps$ is equal to $0.1$ (left column), $1$ (middle column), and $10$ (right column). As expected, increasing the temperature leads to more uniform weights and, as a result, the informed boids orchestrate primitives similarly to uninformed agents.
}
\label{fig:boid_goal_other_cases}
\end{figure}

\clearpage
\begin{figure}
	\centering
\noindent\includegraphics[width=1\linewidth]{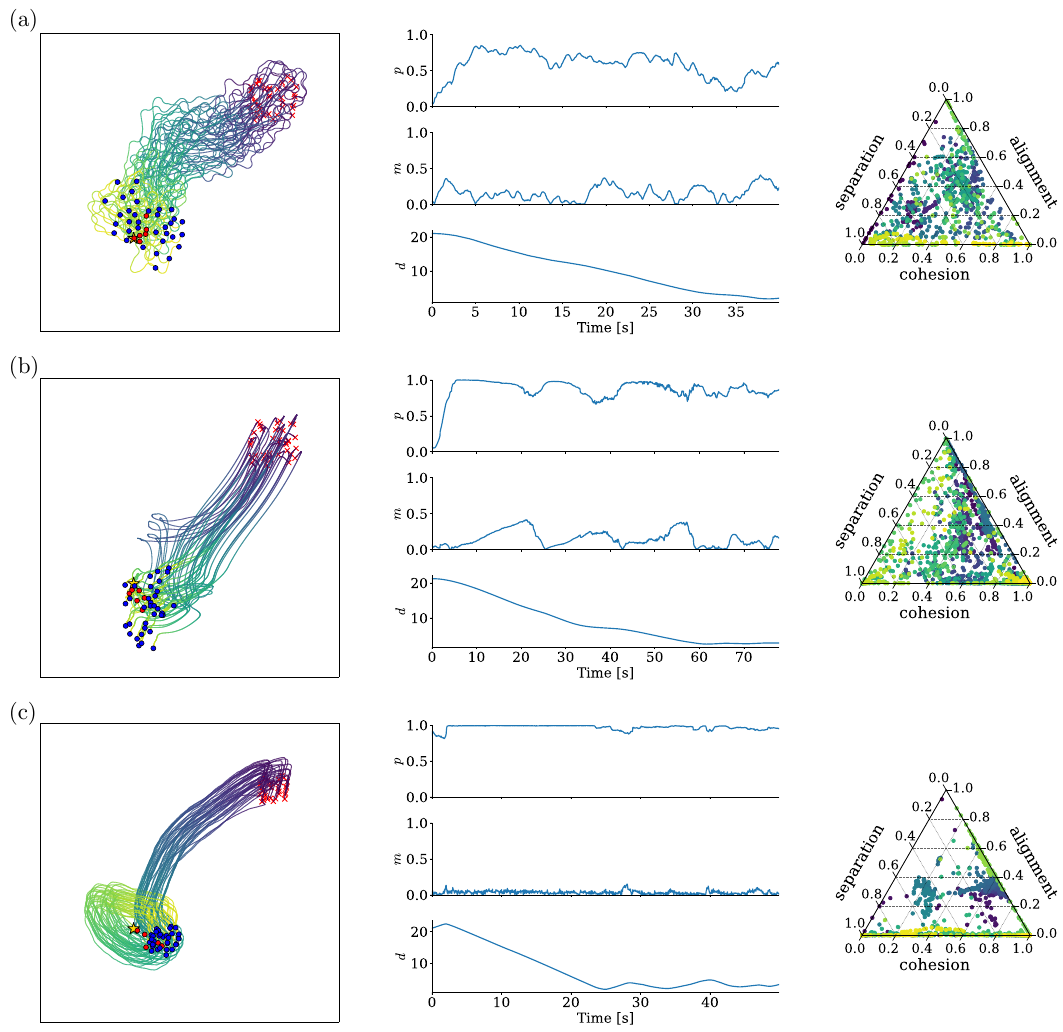}
\caption{Additional soft-control experiments.  Boids start from random initial conditions (crosses); $6$ boids (red-filled circles) informed about the goal position (star) implement our policy composition model, while the uninformed boids (blue-filled circles) evolve according to different models from the literature.
(a)
(Left) The goal-informed boids guide the uninformed boids, which evolve according to the three-zone model from~\cite{IDC-JK-RJ-GDR-NRF:02}, toward the goal.
(Middle) Temporal evolution of group-level metrics. The distance from the goal decreases with time.
(Right) Time evolution of the optimal primitives' weights of a representative informed boid on the simplex.
(b)
(Left) The goal-informed boids implementing our model guide the Cucker-Smale boids~\cite{FC-SS:07} toward the goal.
(Middle) Temporal evolution of group-level metrics.
(Right) Time evolution of the optimal primitives' weights of a representative informed agent on the simplex.
(c)
(Left) The goal-informed boids guide the uninformed Vicsek boids~\cite{TV-AC-EBJ-IC-OS:95} toward the goal.
(Middle) Temporal evolution of group-level metrics.
(Right) Time evolution of the optimal primitives' weights of a representative informed agent on the simplex. Collectively, these experiments confirm the ability of the boids using our model to soft-control the group toward the goal position, while maintaining cohesiveness. Across all models,  the temporal variations in the primitives' weights determined by our model support the tendency of leaders to a flexible use of primitives, without converging to a balanced mixture, consistently with the results in the main text.
}
\label{fig:boid_soft_control}
\end{figure}

\clearpage
\begin{figure}
	\centering
\noindent\includegraphics[width=1\linewidth]{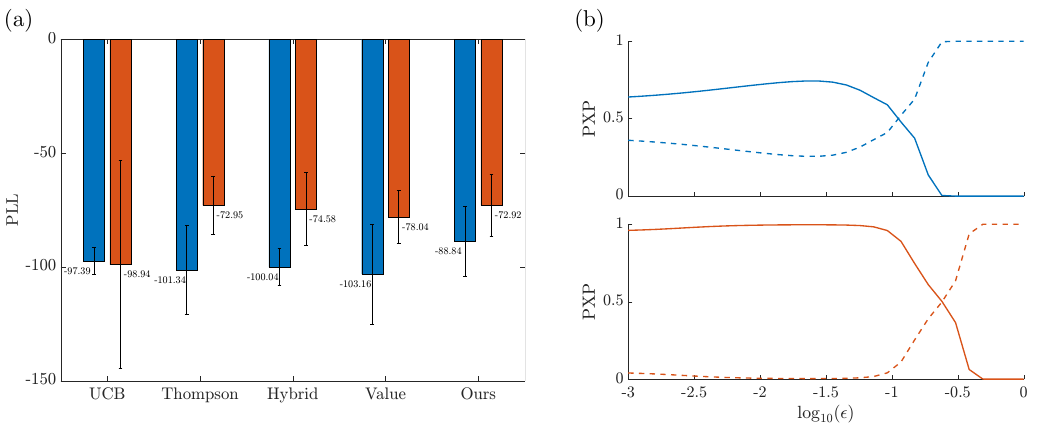}
\caption{Additional results for the exploration/exploitation in human decision-making experiments.
(a) Mean and standard deviation of predictive log-likelihood (PLL) obtained via 3-fold cross-validation. For each fold, models are fitted on two-thirds of the participants and evaluated on the held-out third; results are then averaged across folds. Our model achieves the highest average PLL in Experiment~1 and performs comparably to the best policy (Thompson sampling) in Experiment~2. These results support the findings from the main text. 
(b) PXP for our model (solid line) and the Hybrid model (dashed line) as a function of the regularization parameter $\eps$ in Experiment~1 (top) and Experiment~2 (bottom). These results show that large values of $\eps$ yield lower PXP. This is because, as the regularization parameter increases, the weights tend to become uniform.
}
\label{fig::mab_additional_probit_results}
\end{figure}

\clearpage
\begin{figure}
	\centering
\noindent\includegraphics[width=1\linewidth]{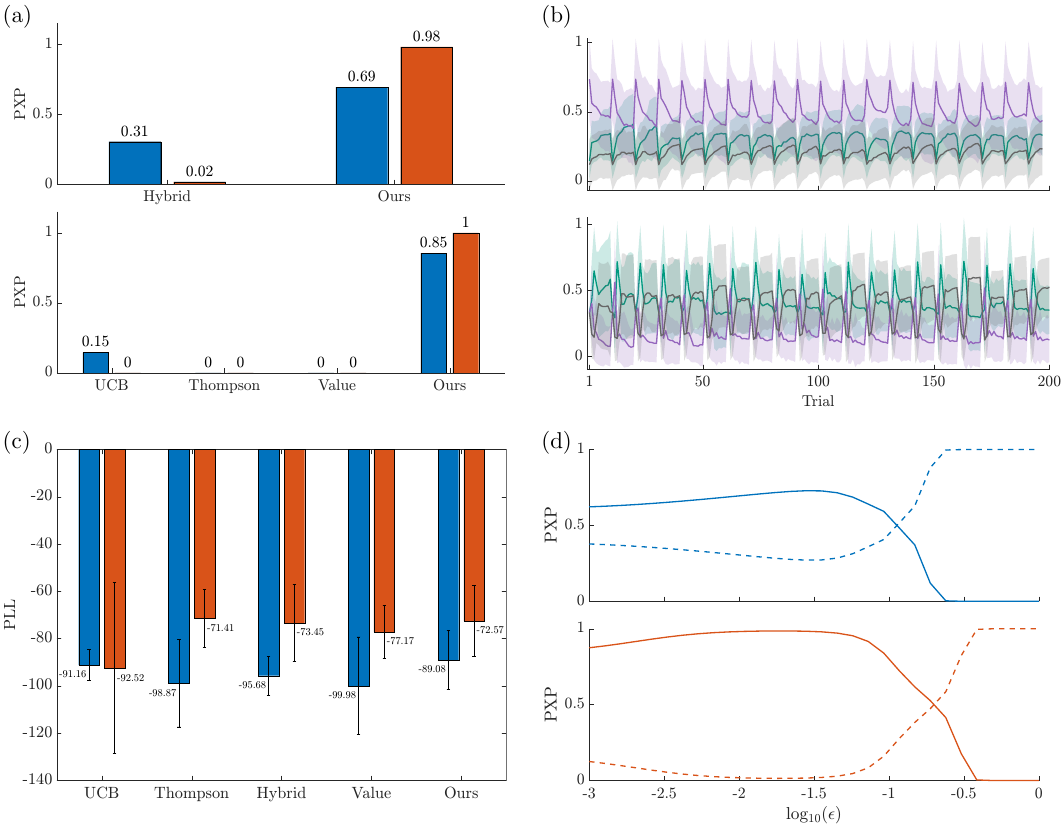}
\caption{Exploration/exploitation in human decision-making experiments using logit regression. (a) Protected exceedance probabilities (PXP) in  Experiment~1 (blue) and Experiment~2 (orange) using logit regression. (Top) Our model outperforms the Hybrid model from~\cite{SJG:18} in both experiments, indicating stronger explanatory power even under the alternative link function. (Bottom) The comparison with UCB, Thompson, and Value shows that our model again achieves the highest PXP in both experiments, consistent with the probit-based results reported in the main text. 
(b) Trial-by-trial evolution of the mean (solid lines) and standard deviation (shaded areas) of the optimal primitives' weights across participants under logit regression. Results for Experiment~1 (top) and 2 (bottom) are shown. The temporal patterns of exploitation (teal), uncertainty-seeking (purple), and risk-aversion (gray) weights closely mirror those obtained with probit regression in the main text, suggesting robustness to the choice of regression link.
(c) Mean and standard deviation of predictive log-likelihood (PLL) obtained via 3-fold cross-validation using logit regression. For each fold, models are fitted on two-thirds of the participants and evaluated on the held-out third; results are then averaged across folds. Our model achieves the highest average PLL in Experiment~1 and performs comparably to the best policy (Thompson sampling) in Experiment~2.
(d) PXP for our model (solid line) and the Hybrid model (dashed line) as a function of the regularization parameter $\eps$, using logit regression. Results in Experiment~1 (top) and Experiment~2 (bottom). 
As in the probit-regression case shown in the main text, large values of $\eps$ yield lower PXP, as the weights tend to become uniform.
Overall, these results confirm that the conclusions drawn in the main text are robust to the choice of link function (probit or logit). 
}
\label{fig::mab_logit}
\end{figure}

\clearpage
\begin{figure}
	\centering
\noindent\includegraphics[width=1\linewidth]{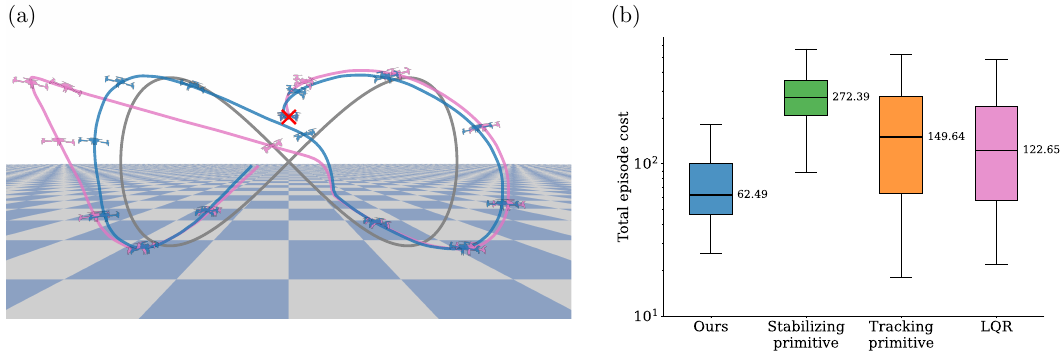}
\caption{Additional layered control experiments combining two linear quadratic regulator (LQR) primitives. 
The environmental settings are as in Fig.~5(b)--(c) with the addition of an impulsive disturbance.
(a) Tracking when our computational model is used (blue). Reference trajectory in grey; trajectory obtained with the baseline LQR in pink. The start position is shown with a red cross and the translucent drone renderings show vehicle configurations during the experiment. The quadrotor equipped with our model recovers faster from the impulsive disturbance.
(b) Comparison of total episode cost between our proposed policy, the individual primitives and the baseline LQR.
Statistics are obtained over $1000$ random seeds.
Renderings generated using \texttt{Safe-Control-Gym}. See Table~\ref{tab:quadrotor_lqr_parameters} for simulation parameters.}
\label{fig::quadrotor_lqr}
\end{figure}

\clearpage

\section{Supplemental Tables}
\label{sec:tables}

\begin{table}[tbph]
\caption{\label{tab:flocking_parameters}%
Parameter values used in the collective-behavior experiments (unless otherwise stated).
Parameters are assumed in accordance with~\cite{CKH-HH:08,RL-YL-LEK:10,IDC-JK-RJ-GDR-NRF:02}.
Initial boid states in the simulations are randomly selected, weights are initialized to the uniform distribution, and the numerical solver \texttt{solve\_ivp} in Python is used.}
\begin{ruledtabular}
\begin{tabular}{lll}
\textrm{Parameter} & \textrm{Description} & \textrm{Value} \\
\colrule
$dt$ & Time step & $0.05$ \\
$N$ & Number of boids & $40$ \\
$d_{\mathrm{u}}$ & Action space discretization & $30^2$ \\
$v_{\max}$ & Maximum velocity norm & $1$ \\
$u_{\max}$ & Maximum acceleration norm & $3$ \\
$\alpha$ & Vision angle & $320^{\circ}$ \\
$r_{\mathrm{sep}}$ & Separation radius & $1$ \\
$r_{\mathrm{ali}}$ & Alignment radius & $3$ \\
$r_{\mathrm{coh}}$ & Cohesion radius & $12$ \\
$\Sigma_{\mathrm{sep}}$ & Covariance of separation primitive & $0.1 I_2$ \\
$\Sigma_{\mathrm{ali}}$ & Covariance of alignment primitive & $0.1 I_2$ \\
$\Sigma_{\mathrm{coh}}$ & Covariance of cohesion primitive & $0.1 I_2$ \\
$\Sigma_{\mathrm{p}}$
& Covariance of $\plantboid{k}{k-1}{i}$, for all $i\in\{1,\ldots,N\}$
& $0.01 {I_4}$ \\
$\Sigma_{\mathrm{q}}$
& Covariance of $\refplantboid{k}{k-1}{i}$, for all $i\in\{1,\ldots,N\}$
& $0.01 I_4$ \\
$\eps$ & Entropy-regularization parameter & $0.5$ \\
$\tau$ & Time-scale parameter & $1$ \\
\end{tabular}
\end{ruledtabular}
\end{table}

\clearpage

\begin{table}[tbp]
\caption{\label{tab:mab_parameters}
Parameter values used in the multi-armed bandit experiments.
Parameters are assumed in accordance with~\cite{SJG:18}.
Weights are initialized to the uniform distribution, and the numerical solver
\texttt{ode23t} in MATLAB is used for weights integration.}
\begin{ruledtabular}
\begin{tabular}{lll}
\textrm{Parameter} & \textrm{Description} & \textrm{Value} \\
\colrule
$\du$ & Action space dimension (number of arms) & $2$ \\
$B$ & Number of blocks per experiment & $20$ \\
$T$ & Trials per block & $10$ \\
$\mu_{1,0}$
& Prior mean of arm 1 (in both experiments)
& $0$ \\
$\mu_{2,0}$
& Prior mean of arm 2 (in both experiments)
& $0$ \\
$s_{1,0}$ & Prior variance of arm 1 & $10$ (Exp.~1), $100$ (Exp.~2)\\
$s_{2,0}$ & Prior variance of arm 2 & $0$ (Exp.~1), $100$ (Exp.~2)\\
$\eps$ & Entropy-regularization parameter & $0.01$ \\
$\tau$ & Time-scale parameter & $1$ \\
\end{tabular}
\end{ruledtabular}
\end{table}

\clearpage
\newpage

\begin{table}[tbp]
\caption{\label{tab:quadrotor_parameters}%
Parameter values used in the quadrotor experiments with composition of PD primitives.
Weights are initialized to the uniform distribution, and the numerical solver
\texttt{solve\_ivp} in Python is used for weights integration.}
\begin{ruledtabular}
\begin{tabular}{lll}
\textrm{Parameter} & \textrm{Description} & \textrm{Value} \\
\colrule
$dt$ & Sampling time & $0.02$ \\
$\mathbf{x}_0$ & Initial state & $[0,0,1.2,0,0,0]^\top$ \\
$T$ & Episode length & $300$ \\
$\mathbf{x}_\mathrm{eq}$ & Hover equilibrium state & $\0_6$ \\
$\mathbf{u}_\mathrm{eq}$ & Hover equilibrium input & $[0.1323,\, 0.1323]^\top$ \\
$\du$ & Action space discretization & $30^2$ \\
$\sigma_{\mathrm{d}}$ & Dynamics disturbance standard deviation & $0.03$ \\
$\Sigma_\mathrm{u}$ & Actuator-noise covariance & $(0.01)^2I_2$ \\
$Q$ & State cost matrix & $\mathrm{diag}(1,0.1,1,0.1,0.1,0.1)$ \\
$R$ & Control cost matrix & $0.1I_2$ \\
$k^{\mathrm{s}}_{\mathrm{d},x}$
& Horizontal derivative gain of stabilization primitive
& $0.05$ \\
$k^{\mathrm{s}}_{\mathrm{d},z}$
& Vertical derivative gain of stabilization primitive
& $0.1$ \\
$k^{\mathrm{s}}_{\mathrm{p},\theta}$
& Attitude proportional gain of stabilization primitive
& $70000$ \\
$k^{\mathrm{s}}_{\mathrm{d},\theta}$
& Attitude derivative gain of stabilization primitive
& $20000$ \\
$k^{\mathrm{t}}_{\mathrm{p},x}$
& Horizontal proportional gain of tracking primitive
& $0.45$ \\
$k^{\mathrm{t}}_{\mathrm{d},x}$
& Horizontal derivative gain of tracking primitive
& $0.2$ \\
$k^{\mathrm{t}}_{\mathrm{p},z}$
& Vertical proportional gain of tracking primitive
& $1.0$ \\
$k^{\mathrm{t}}_{\mathrm{d},z}$
& Vertical derivative gain of tracking primitive
& $0.35$ \\
$k^{\mathrm{t}}_{\mathrm{p},\theta}$
& Attitude proportional gain of tracking primitive
& $14000$ \\
$k^{\mathrm{t}}_{\mathrm{d},\theta}$
& Attitude derivative gain of tracking primitive
& $10000$ \\
$\eps$ & Entropy-regularization parameter & $0.05$ \\
$\tau$ & Time-scale parameter & $1$ \\
\end{tabular}
\end{ruledtabular}
\end{table}

\clearpage

\begin{table}[tbp]
\caption{\label{tab:quadrotor_lqr_parameters}%
Parameter values used in the quadrotor experiments with composition of LQR primitives.
Weights are initialized to the uniform distribution, and the numerical solver
\texttt{solve\_ivp} in Python is used for weights integration.}
\begin{ruledtabular}
\begin{tabular}{lll}
\textrm{Parameter} & \textrm{Description} & \textrm{Value} \\
\colrule
$dt$ & Sampling time & $0.02$ \\
$\mathbf{x}_0$ & Initial state & $[0,0,1.2,0,0,0]^\top$ \\
$T$ & Episode length & $300$ \\
$\du$ & Action space discretization & $30^2$ \\
$\sigma_{\mathrm{d}}$ & Dynamics disturbance standard deviation & $0.03$ \\
$\Sigma_\mathrm{u}$ & Actuator-noise covariance & $(0.01)^2I_2$ \\
$m_{\mathrm{i}}$ & Impulse disturbance magnitude & $1.0$ \\
$Q$ & State cost matrix & $\mathrm{diag}(1,0.1,1,0.1,0.1,0.1)$ \\
$R$ & Action cost matrix & $0.1I_2$ \\
$Q_{\mathrm{s}}$
& Stabilizing primitive state cost matrix
& $\mathrm{diag}(10^{-6},0.1,10^{-6},10^{-5},0.5,0.1)$ \\
$R_{\mathrm{s}}$
& Stabilizing primitive action cost matrix
& $0.01I_2$ \\
$Q_{\mathrm{t}}$
& Tracking primitive state cost matrix
& $\mathrm{diag}(5,0.1,5,0.1,0.1,0.1)$ \\
$R_{\mathrm{t}}$
& Tracking primitive action cost matrix
& $0.1I_2$ \\
$Q_{\mathrm{b}}$
& Baseline LQR state cost matrix
& $\mathrm{diag}(1,0.1,1,0.1,0.1,0.1)$ \\
$R_{\mathrm{b}}$
& Baseline LQR action cost matrix
& $0.1I_2$ \\
$\eps$ & Entropy-regularization parameter & $0.05$ \\
$\tau$ & Time-scale parameter & $1$ \\
\end{tabular}
\end{ruledtabular}
\end{table}

\clearpage
\bibliography{bib_arxiv}

%% file: my_commands.tex
\DeclareMathOperator*{\argmin}{arg\,min}

\newcommand{\R}{\mathbb{R}}
\newcommand{\E}{\mathbb{E}}

\newcommand{\norm}[1]{\left\Vert #1 \right\Vert}

\newcommand{\bv}[1]{\mathbf{#1}}

\newcommand{\refjointxutwisted}[2]{\tilde q(\bv{x}_{#1},\bv{u}_{#1}\mid \bv{x}_{#2} )}

\newcommand{\jointxu}[2]{p(\bv{x}_{#1},\bv{u}_{#1}\mid \bv{x}_{#2} )}
\newcommand{\refjointxu}[2]{q(\bv{x}_{#1},\bv{u}_{#1}\mid \bv{x}_{#2} )}

\newcommand{\plant}[2]{p(\bv{x}_{#1}\mid \bv{x}_{#2}, \bv{u}_{#1} )}
\newcommand{\refplant}[2]{q(\bv{x}_{#1}\mid \bv{x}_{#2}, \bv{u}_{#1} )}

\newcommand{\tilderefplant}[2]{{\tilde q}(\bv{x}_{#1}\mid \bv{x}_{#2}, \bv{u}_{#1} )}

\newcommand{\policy}[2]{p(\bv{u}_{#1}\mid \bv{x}_{#2} )}
\newcommand{\refpolicy}[2]{q(\bv{u}_{#1}\mid \bv{x}_{#2} )}

\newcommand{\tilderefpolicy}[2]{\tilde{q}(\bv{u}_{#1}\mid \bv{x}_{#2} )}
\newcommand{\primitive}[3]{\pi^{{#3}}(\bv{u}_{#1}\mid \bv{x}_{#2} )}
\newcommand{\primitives}[2]{\bv{\pi}(\bv{u}_{#1}\mid \bv{x}_{#2} )}

\newcommand{\plantboid}[3]{p (\bv{x}_{#1}^{#3}\mid \bv{x}_{#2}^{#3}, \bv{u}_{#1}^{#3} )}
\newcommand{\refplantboid}[3]{q (\bv{x}_{#1}^{{#3}}\mid \bv{x}_{#2}^{{#3}}, \bv{u}_{#1}^{{#3}} )}

\newcommand{\refpolicyboid}[3]{q(\bv{u}_{#1}^{{#3}}\mid \bv{x}_{#2}^{{#3}} )}

\newcommand{\weight}[2]{w_{#1}^{#2}}
\newcommand{\weights}[1]{\bv{w}_{#1}}
\newcommand{\optimalweights}[1]{\bv{w}^{\star}_{#1}}

\newcommand{\pweights}[1]{\hat{\bv{w}}_{#1}}

\newcommand{\optimalpolicy}[2]{p^{\star}(\bv{u}_{#1}\mid \bv{x}_{#2} )}

\newcommand{\KL}{\textup{KL}}
\newcommand{\DKL}[2]{{D}_{\KL}(#1 \parallel #2 )}
\newcommand{\DKLb}[2]{{D}_{\KL}\bigl(#1 \parallel #2 \bigr)}

\newcommand{\states}[1]{\mathbf{x}_{#1}}

\usepackage[]{mdframed}

\newcommand{\e}{\mathrm{e}}

\newcommand{\1}{\mbox{\fontencoding{U}\fontfamily{bbold}\selectfont1}}

\newcommand{\prox}[1]{\mathrm{prox}_{#1}}

\newcommand{\softmax}{\operatorname{softmax}}
\newcommand{\entropy}{\mathsf{H}}
\newcommand{\entropybar}{\subscr{\entropy}{barrier}}

\newcommand{\np}{\subscr{n}{$\pi$}}
\newcommand{\du}{\subscr{d}{u}}
\newcommand{\F}{\mathsf{F}}
\newcommand{\simplex}[1]{\Delta_{#1}}
\newtheorem{definition}{Definition}
\newtheorem{thm}{Theorem}
\newtheorem{lem}{Lemma}
\newtheorem{pro}{Proposition}

\newtheorem{remark}{Remark}

\newcommand{\bd}{\begin{definition}} 
\newcommand{\ed}{\end{definition}} 
\newcommand{\bp}{\begin{pro}} 
\newcommand{\ep}{\end{pro}} 
\newcommand{\bt}{\begin{thm}} 
\newcommand{\et}{\end{thm}}
\newcommand{\blm}{\begin{lem}} 
\newcommand{\elm}{\end{lem}}
\let\eps\varepsilon

\newcommand{\bi}{\begin{itemize}}
\newcommand{\ei}{\end{itemize}}
\newcommand{\beq}{\begin{equation}}
\newcommand{\eeq}{\end{equation}}

\newcommand{\subscr}[2]{#1_{\textup{#2}}}

